\documentclass{elsarticle}
\usepackage{xcolor}
\usepackage{lineno,hyperref}
\usepackage{amsthm,amsmath,amscd,amsfonts,mathtools,stmaryrd,subcaption}
\usepackage[vmargin=2.5cm,hmargin=2.5cm,headheight=14.5pt,top=2cm,headsep=.5cm]{geometry}
\usepackage{soul}
\biboptions{sort&compress}
\modulolinenumbers[5]

\newtheoremstyle{ptheorem}{1em}{0em}{\itshape}{}{\bfseries}{.}{.5em}{\thmname{#1}\thmnumber{
		#2}\thmnote{ (\hspace{-.01pt}{#3})}}

\theoremstyle{ptheorem}

\newtheorem{thm}{Theorem}[section]
\newtheorem{pro}[thm]{Proposition}
\newtheorem{lem}[thm]{Lemma}
\newtheorem{cor}[thm]{Corollary}

\newtheoremstyle{hdef}{1em}{0em}{}{}{\bfseries}{.}{.5em}{\thmname{#1}\thmnumber{
		#2}\thmnote{ (\hspace{-.01pt}{#3})}}
\theoremstyle{hdef}

\newtheorem{dfn}[thm]{Definition}
\newtheorem{rem}[thm]{Remark}

\newtheorem{exa}[thm]{Example}

\numberwithin{equation}{section}
\numberwithin{figure}{section}

\journal{Journal of Mathematical Analysis and Applications}

\begin{document}

\begin{frontmatter}

\title{On first and second order linear Stieltjes differential equations}

\author{Francisco J. Fernández, Ignacio Marqu\'ez Alb\'es and F. Adri\'an F. Tojo
}
\address{Departamento de Estat\'{\i}stica, An\'alise Matem\'atica e Optimizaci\'on \\ Instituto de Matem\'aticas \\ Universidade de Santiago de Compostela \\ 15782, Facultade de Matem\'aticas, Campus Vida, Santiago, Spain. \\ e-mail: fjavier.fernandez@usc.es, ignacio.marquez@usc.es,  fernandoadrian.fernandez@usc.es}
%\fntext[myfootnote]{Corresponding author.}


\begin{abstract}
This work deals with the obtaining of solutions of first and second order Stieltjes differential equations. We define the notion of Stieltjes derivative on the whole domain of the functions involved, provide a notion of $n$-times continuously Stieltjes-differentiable functions and prove existence and uniqueness results of Stieltjes differential equations 
in the space of such functions. We also present the Green's functions associated to the different problems and an application to the Stieltjes harmonic oscillator.
\end{abstract}


\end{frontmatter}

\linenumbers

\section{Introduction}

There has been a recent surge in the study of Stieltjes differential equations focused on obtaining applicable results comparable to those available for classical derivatives \cite{LoRo14,FerFraTo2020,FriLo17,FriTo20,LoMa18,LoMa19,LoMa19Resolution,LoMa20,LoMaMon18,LoMaRo20,MonSat19,MarquezTesis,Ma21,MaMon20,MaMulti,SatSmy20}. These works center their attention in the procuring of solutions of first order differential equations and systems. The theory developed starts with the obtaining of simple solutions, like the solution of the first order linear problem \cite{FriLo17,FriTo20}, which is identified with the exponential, in order to, later, prove existence and uniqueness results in more general settings \cite{LoMa19Resolution,Ma21,MaMulti}. Some of these works also provide interesting practical applications \cite{LoMa18,LoMaMon18} and others  generalize the framework in several ways, such as allowing for sign changing derivators \cite{FriTo20}, considering several different derivators \cite{MaMulti} or generalizing the concept of Stieltjes derivative \cite{MaTo19}.

In any case, all of the aforementioned works restrict themselves to the first order case. The reason behind this is that, in order to study higher order problems, the notion of higher order Stieltjes derivative has to be correctly defined, which is not obvious. In fact, the first difficulty lies on the mere definition of the Stieltjes derivative, which, to the best of our knowledge, is nowhere defined in the literature on the whole domain of definition of the function, something which impedes taking a second derivative.

In this work we provide this definition, which enables us to study second order problems. First, we consider the Stieltjes derivative in the whole of the domain of the function, which allows us to talk about the space of continuously Stieltjes-differentiable functions in the same way we speak of the space of continuously differentiable functions with the usual derivative. We can then explore the first order problem in this space, obtaining existence and uniqueness results that mirror those of the previous works. In fact, we profit from the opportunity of revisiting the solution of the first order linear problem to provide a constructive way of obtaining its solution. All of these steps are also taken with a further generalization: our functions are allowed to take real or complex values. Furthermore, we obtain the explicit expression of the Green's function of the first order linear problem with initial conditions and we construct the Stieltjes versions of the sine and cosine functions using the complex version of the Stieltjes exponential.

Once we have studied the first order problem with various degrees of regularity (something the  subsequent spaces of $n$-times continuously Stieltjes-dif\-fer\-en\-ti\-a\-ble functions allow), we move on to study  second order problems. First, we present existence and uniqueness results for the homogeneous second order problem with constant coefficients and then we study the non homogeneous case with varying degrees of regularity. Here we also obtain the explicit expression of the Green's function of the second order linear problem with initial conditions. All this work is then illustrated with an application to the Stieltjes harmonic oscillator for which we also analyze the resonance effect. Finally, in order to validate the explicit solutions obtained, we compare them with the numerical approximation of the corresponding first order linear system using the numerical scheme introduced in \cite{FerFraTo2020}.

The structure of this work is as follows: In Section~\ref{SPr} we present some preliminary concepts and we prove several results related to Lebesgue-Stieltjes integral. In Section~\ref{SBCg} we introduce the space of bounded Stieltjes differentiable functions and analyze some of its properties. We study the first order linear Stieltjes differential equation in Section~\ref{SFO}, including in the complex case. In this Section we also define the complex Stieltjes exponential and the Stieltjes version of the sine and cosine functions. In Section~\ref{SSO} we study the homogeneous Stieltjes second order problem with constant coefficients, the non homogeneous case and we also obtain an explicit solution for both situations. Finally, in Section~\ref{SSHO} we present an application to the Stieltjes harmonic oscillator. We obtain the explicit solution of the overdamped, critically damped and underdamped cases, and provide an example in which the resonance effect appears. In order to validate the explicit solution obtained, we compare it with the numerical solution of the corresponding first order lineal system.

\section{Preliminaries}\label{SPr}
%%%%%%%%%%%%%%%%%%%%%%%%%%%%%%%%%%%%%%%%%%%%%%%%%%%%
%%%%%%%%%%%%%%%%%%%%%%%%%%%%%%%%%%%%%%%%%%%%%%%%%%%%

Let $[a,b]\subset{\mathbb R}$ be an interval, ${\mathbb F}$ the field ${\mathbb R}$ or ${\mathbb C}$ and $g: \mathbb{R} \to{\mathbb R}$ a left-continuous non-decreasing function. We will refer to such functions as \emph{derivators}. For these functions, we define the set $D_g=\{d_n\}_{n \in  \Lambda}$ (where $ \Lambda\subset\mathbb{N}$) as the set of all discontinuity points of $g$, namely, $D_g=\{ t \in \mathbb R  :  \Delta^+ g(t)>0\}$ where $\Delta^+ g(t):=g(t^+)-g(t)$, $t\in\mathbb R$, and $g(t^+)$ denotes the right hand side limit of $g$ at $t$. We also define 
	\begin{equation}
		C_g:=\{ t \in \mathbb R \, : \, \mbox{$g$ is constant on $(t-\varepsilon,t+\varepsilon)$ for
			some $\varepsilon>0$} \}.
	\end{equation}
Observe that $C_g$ is open in the usual topology of $\mathbb R$, so we can write 
\begin{equation}\label{Cgdecomp}
	C_g=\bigcup_{n \in \widetilde\Lambda}\{(a_n,b_n)\}
\end{equation}
 where 
$\widetilde\Lambda\subset \mathbb{N}$ and $(a_k,b_k)\cap (a_j,b_j)=
\emptyset$ for $k\neq j$. With this notation, we denote $N_g^-:=\{a_n\}_{n \in \widetilde\Lambda}\backslash D_g$, $N_g^+:=\{b_n\}_{n \in \widetilde\Lambda}\backslash D_g$ and
$N_g:=N_g^-\cup N_g^+$.
\begin{rem}
	For the aims of this paper, we will assume without loss of generality that $g(a)=0$. Furthermore, we will also assume that $g$ is continuous at $x=a$. As pointed out in \cite[p. 21]{FriLo17} and \cite[Proposition~4.28]{MarquezTesis}, the continuity assumption has no impact in the study of differential equations, which is our final goal. 
	Finally, in order to properly define the Stieltjes derivative in the whole $[a,b]$, we will also ask that $[a,b]\setminus C_g\neq \emptyset$.
\end{rem}
\medskip

We define 
$g^B:\mathbb{R}\to \mathbb{R}$ as:
\begin{equation}\label{gB}
	g^B(t)=\begin{cases}
		\displaystyle
		\sum_{s\in [a,t)\cap D_g}\Delta^+ g(s), & t>a, \\
		\displaystyle
		-\sum_{s\in [t,a)\cap D_g} \Delta^+g(s), & t\leq a.
	\end{cases}
\end{equation}
It is clear that $g^B$ is a left-continuous and non-decreasing function. 
Moreover, the map $g^C:\mathbb{R}\to{\mathbb R}$ given 
by 
\begin{equation}\label{gC}
	g^C(t):=g(t)-g^B(t),
\end{equation} 
is also non-decreasing and continuous. We say $g^C$ that is the \emph{{continuous part}} of $g$ and $g^B$ is the \emph{{jump part}} of $g$. Observe that 
both $g^B$ and $g^C$ are continuous at $x=a$ and $g^C(a)=g^B(a)=0$. 
\\

Throughout this work we consider the Lebesgue--Stieltjes measure space
$(\mathbb{R},\mathcal{M}_g,\mu_g)$, where $\mathcal{M}_g$ and $\mu_g$ are
the $\sigma$-algebra and measure constructed
in an analogous fashion to the classical Lebesgue measure, where the length of $[c,d)$
is given by $\mu_g([c,d))=g(d)-g(c)$. The~interested reader may refer to~\cite{LoRo14}
for details concerning this measure space. We must emphasize that, in the case of considering
$g(t)=t$, we recover the classic Lebesgue measure space that we will denote by
$(\mathbb{R},\mathcal{L},\mu)\equiv (\mathbb{R},\mathcal{M}_{\operatorname{Id}},\mu_{\operatorname{Id}})$ where $\operatorname{Id}$ is the identity function. Furthermore, we can define the measure space associated with the continuous part, 
$(\mathbb{R},\mathcal{M}_{g^C},\mu_{g^C})$, the jump part, 
$(\mathbb{R},\mathcal{M}_{g^B},\mu_{g^B})$, and the one associated 
with the derivator itself, $(\mathbb{R},\mathcal{M}_g,\mu_g)$. If we 
denote by ${\mathcal B}(\tau_u)$ the Borel $\sigma$-algebra associated to $\tau_u$ in $\mathbb{R}$, we have that ${\mathcal B}(\tau_u)\subset 
\mathcal{M}_g$ and also ${\mathcal B}(\tau_u)\subset \mathcal{M}_{g^M}$, 
with $M=C,B$. We must mention that if $E\subset \mathbb{R}$, 
$\mu^*_{g^M}(E)\leq \mu^*_{g}(E)$, for 
$M=C,B$, being $\mu^*_{g^M}$ and $\mu^*_g$ the outer measures 
associated to $g^M$ and $g$ respectively, $M=C,B$. We also have that 
if $E\subset \mathbb{R}$ is a bounded set, then $\mu^*_g(E)<\infty$.
\\

 We have the following lemma that, in particular, provides us with a relationship between the 
$\sigma$-algebras $\mathcal{M}_g$, $\mathcal{M}_{g^C}$ and 
$\mathcal{M}_{g^B}$.

\begin{lem} \label{sigmaalgebras} The following properties hold for the maps $g$, $g^C$ and $g^B$:
\begin{enumerate}
\item Given an element $E\in \mathcal{M}_g$ there exists 
$H\in G_{\delta}$ (that is, $H$ is a countable intersection of 
open sets) and $N\in \mathcal{M}_{g}$ such that $E\subset H$, 
$N\subset H$, $\mu_g(N)=0$ and $E=H\setminus N$.
\item Given an element $E\in \mathcal{M}_g$ there exists 
$F \in F_{\sigma}$ (that is, $F$ is a countable union 
of closed sets) and $N\in \mathcal{M}_{g}$ such that $\mu_g(N)=0$, 
$F\cap N=\emptyset$ and $E=F\cup N$.
\item $\mathcal{M}_g \subset \mathcal{M}_{g^C}$
\item $\mathcal{M}_{g^B}=\mathcal{P}(\mathbb{R})$.
\end{enumerate}
\end{lem}

\begin{proof} Since ${\mathcal B}(\tau_u)\subset \mathcal{M}_g$ and 
$g$ is left-continuous, we have that:
\begin{equation}
\mu_g(E)=\inf \left\{ \sum_{n \in \mathbb{N}} \mu_g([a_n,b_n)):\;
E\subset \bigcup_{n \in \mathbb{N}} [a_n,b_n) 
\right\}  \inf \left\{ \sum_{n \in \mathbb{N}} \mu_g((a_n,b_n)):\;
E\subset \bigcup_{n \in \mathbb{N}} (a_n,b_n) \right\}.
\end{equation}
Indeed on the one hand given $\{(a_n,b_n)\}_{n \in \mathbb{N}}$ such 
that $E\subset \cup_{n \in \mathbb{N}} (a_n,b_n)$, we have that 
\begin{equation}
\mu_g^*(E)\leq \mu_g^*\left( \bigcup_{n \in \mathbb{N}} (a_n,b_n) \right) 
\leq \sum_{n \in \mathbb{N}} \mu_g^*((a_n,b_n)),
\end{equation}
therefore $\mu_g(E)\leq \inf \left\{ \sum_{n \in \mathbb{N}} \mu_g((a_n,b_n)):\;
E\subset \cup_{n \in \mathbb{N}} (a_n,b_n) \right\}$. On the other hand, 
given $\varepsilon>0$ and $\{[a_n,b_n)\}_{n \in \mathbb{N}}$ 
such that $E\subset \cup_{n \in \mathbb{N}} [a_n,b_n)$, we have, 
thanks to the left-continuity of $g$, that 
there exists $\{(\tilde{a}_n,b_n)\}_{n \in \mathbb{N}}$ such that 
$[a_n,b_n)\subset (\tilde{a}_n,b_n)$ and $\mu_g^*((\tilde{a}_n,b_n))\leq 
\mu_g^*([a_n,b_n))+\varepsilon/2^n$, $n \in \mathbb{N}$. Thus 
$\sum_{n \in \mathbb{N}} \mu_g^* (\tilde{a}_n,b_n) \leq 
\sum_{n \in \mathbb{N}} \mu_g^* ([a_n,b_n))+\varepsilon$ and we 
conclude, taking infimum in both sides of inequality, that 
$\inf \left\{ \sum_{n \in \mathbb{N}} \mu_g((a_n,b_n)):\;
E\subset \cup_{n \in \mathbb{N}} (a_n,b_n) \right\}\leq \mu_g(E)$. Now, 
we can proceded as in \cite[Corollary 15.5 and 15.8]{Bartle1995} 
to obtain 1 and 2, respectively.

Now, for 3, given an element $E\in \mathcal{M}_g$, 
there exists $F \in F_{\sigma}$ and $N\in \mathcal{M}_{g}$ 
such that $\mu_g(N)=0$, 
$F\cap N=\emptyset$ and $E=F\cup N$. Now, 
$F\in {\mathcal B}(\tau_u)\subset \mathcal{M}_c$ and 
$\mu_{g^C}^*(N)\leq \mu_g^*(N)=0$ so we have that $N \in \mathcal{M}_{g^C}$ 
since $(\mathbb{R},\mathcal{M}_{g^C},\mu_{g^C})$ is a complete 
measure space. Therefore, $E\subset \mathcal{M}_{g^C}$.

Finally, for $E\in \mathcal{P}(\mathbb{R})$, we have that 
$E=(E\setminus D_{g^B}) \cup (E\cap D_{g^B})$. Now, 
$(E\setminus D_{g^B})\subset C_{g^B}$ and then 
$\mu^*_{g^B}(E\setminus D_{g^B})=0$, so $E\setminus D_{g^B}\in 
\mathcal{M}_{g^B}$. Finally $E\cap D_{g^B}\in \mathcal{B}
(\tau_u)\subset \mathcal{M}_{g^B}$. Therefore 
$E\in \mathcal{M}_{g^B}$,  which finishes the proof of 4.
\end{proof}

We denote by
$\mathcal{L}^1_g([a,b);{\mathbb F})$ the set of functions $f:[a,b)\rightarrow {\mathbb F}$ such that their real and imaginary parts, that is,
$\operatorname{Re}(f)$ and $\operatorname{Im}(f)$ respectively, are measurable and $\int_{[a,b)}|f|\, \operatorname{d} \mu_g <\infty$. For this
class of functions we define
\begin{equation}%
	\int_{[a,b)} f\, \operatorname{d} \mu_g=\int_{[a,b)} \operatorname{Re}(f)\,\operatorname{d} \mu_g + i \, \int_{[a,b)}\operatorname{Im}(f)\, \operatorname{d} \mu_g.\end{equation}

\begin{lem}\label{lemint} Given a function $f\in \mathcal{L}^1_g([a,b);{\mathbb F})$,
	\begin{equation}
		\int_{[a,t)} f\, \operatorname{d}\mu_g=\int_{[a,t)} f \operatorname{d}\mu_{g^C}+
		\sum_{s \in [a,t)\cap D_g}
		f(s)\Delta^+ g(s),\; \forall t \in [a,b].\end{equation}
\end{lem}

\begin{proof} Given a function $f\in \mathcal{L}^1_g([a,b);{\mathbb F})$, thanks to
	Lemma \ref{sigmaalgebras} and the fact that $\mu_{g^C}(E)\leq \mu_g(E)$
	for all $E\in \mathcal{M}_{g}$, we have that $f\in \mathcal{L}^1_{g^C}([a,b);{\mathbb F})$.
	Now thanks to
	\cite[Theorems 6.3.13, 6.12.3 and 6.12.7]{MonSlaTvr18}, and separating the
	real and imaginary part if necessary,  we have the desired result.
\end{proof}

\begin{cor} Given $E\in \mathcal{M}_g$ and taking $f=\chi_{E}$ (the characteristic function associated to $E$) in Lemma~\ref{lemint} we have that
	\begin{equation}
		\mu_g(E)=\mu_{g^C}(E)+\sum_{s\in E\cap D_g} \Delta^+ g(s).\end{equation}
\end{cor}

We now introduce a tool that would allow us to transform Lebesgue-Stieltjes integrals with respect to $g^C$ into the usual Lebesgue ones. In particular, in light of Lemma~\ref{lemint}, this means that we will have a way of transforming any Lebesgue-Stieltjes integral into a Lebesgue one.

\begin{dfn}[Pseudo-inverse of $g^C$]\label{dfnpi} { 
Given an interval $[a,b]$ and a derivator $g:\mathbb{R}\to \mathbb{R}$,
	we define the \emph{pseudo-inverse} of the continuous part $g^C$ in the the interval $[0,g^C(b)]$ by:
	\begin{equation} \label{pseudoinverse}
		\gamma: x \in [0,g^C(b)] \to \gamma(x)=\min\{t\in [a,b]:g^C(t)=x\}\in [a,b].\end{equation}}
\end{dfn}

In \cite[Lemma 4.1]{LoMa20} we can find some of the properties of the pseudo-inverse of a continuous derivator mapping the real line onto the real line. For our context, by extending linearly  the map $g$ outside of the interval $[a,b]$ we can obtain the required property, which leads to the following result.

\begin{pro} We have the following properties for the pseudo-inverse of the continuous part $g^C$ in the interval $[0,g^C(b)]$:
	\begin{itemize}
		\item For all $x\in [0,g^C(b)]$, $g^C(\gamma(x))=x$.
		\item For all $t \in [a,b]$, $\gamma(g^C(t))\leq t$.
		\item For all $ t\in [a,b]$, $t\notin C_{g^C}\cup N_{g^C}^+$, $\gamma(g^C(t))=t$.
		\item The map $\gamma$ is strictly increasing.
		\item The map $\gamma$ is left-continuous everywhere and continuous at
		every $x\in [0,g^C(b)]$, $x\notin g^C(C_g)$.
	\end{itemize}
\end{pro}
Now we are ready to prove the following result.
\begin{pro}  \label{morfismo} Given an interval $[a,b]$ and a derivator $g:\mathbb{R}\to \mathbb{R}$:
	\begin{enumerate}
		\item The continuous part
		\begin{equation}
			g^C:([a,b],\mathcal{M}_{g^C})\to ([0,g^C(b)],\mathcal{L})\end{equation}
		is a measurable morphism.
		\item The pseudo-inverse of the continuous part $g^C$
		\begin{equation}
			\gamma:([0,g^C(b)],\mathcal{L}) \to
			([a,b],\mathcal{M}_{g^C})\end{equation}
		is a measurable morphism.
		%
		%Furthermore,  $\mu_{g^C}^*(C_{g^C} \cup N_{g^C})=0$ and $\mu_{g^C}^*((g^C)^{-1}(N))=0$.
	\end{enumerate}
\end{pro}

\begin{proof} Let us prove the two statements separately.

	1. Let us consider a subset $E\subset [0,g^C(b)]$ such that $E\in \mathcal{L}$. We have that there exists $F\in F_{\sigma}$ and $N\in \mathcal{L}$ with
	$\mu(N)=0$ such that $F\cap N=\emptyset$ and $E=F\cup N$. It is clear
	that $(g^C)^{-1}(F)\in \mathcal{B}(\tau_u)$ so, if we prove that
	$\mu_{g^C}^*((g^C)^{-1}(N))=0$, where $\mu_{g^C}^*$ is the the outer Lebesgue-Stieltjes measure, we will have finished.

	Now, since $\mu(N)=0$, given $\varepsilon>0$, there exists a countable
	disjoint family $\{[\widetilde{c}_n,\widetilde{d}_n)\}_{n \in \mathbb{N}}$ such that
	$N\subset \bigcup_{n \in \mathbb{N}} [\widetilde{c}_n,\widetilde{d}_n)$ and
	$\sum_{n \in \mathbb{N}}(\widetilde{d}_n-\widetilde{c}_n)<\varepsilon$. We have that
	$(g^C)^{-1}([\widetilde{c}_k,\widetilde{d}_k))=[\gamma(\widetilde{c}_k),
	\gamma(\widetilde{d}_k))$, for all $k\in \mathbb{N}$, thus
	$(g^C)^{-1}(N)\subset \bigcup_{n \in \mathbb{N}} [
	\gamma(\widetilde{c}_n),\gamma(\widetilde{d}_n))$. Finally,
	\begin{equation}\mu^*_{g^C}((g^C)^{-1}N)\leq \sum_{n\in \mathbb{N}}\mu^*_{g^C}
		[\gamma(\widetilde{c}_n),\gamma(\widetilde{d}_n))=
		\sum_{n \in \mathbb{N}} \left(g^C(\gamma(\widetilde{d}_n))-
		g^C(\gamma(\widetilde{c}_n))\right)=\sum_{n \in \mathbb{N}}
		(\widetilde{d}_n-\widetilde{c}_n)<\varepsilon.\end{equation}
	Since $\varepsilon>0$ was arbitrarily chosen, we have that $\mu^*_{g^C}((g^C)^{-1}N)=0$, which finishes the proof of 1.
\\

	2. Let us consider a subset $E\subset [a,b]$ such that $E\in \mathcal{M}_{g^C}$.
	We have { that there exists
	$F \in F_{\sigma}$ and $N\in \mathcal{M}_{g^C}$ such that $\mu_g(N)=0$, $F\cap N=
	\emptyset$ and $E=F\cup N$}. Thus, we conclude that $\gamma^{-1}(E)=
	\gamma^{-1}(F)\cup \gamma^{-1}(N)$. Now, since $\gamma$ is strictly increasing, it is a Borel map, so we
	have that $\gamma^{-1}(F)\in \mathcal{B}(\tau_u)\subset \mathcal{L}$. Hence, if we prove
	that $\mu^*(\gamma^{-1}(N))=0$, where $\mu^*$ is the outer Lebesgue measure,
	we are done. The proof in this case is analogous to the previous
	one, the only difference lies in that, given an interval $[c,d)$, we have
	that $\gamma^{-1}([c,d))\subset [g^C(c),g^C(d)]$, thus $\mu^*(\gamma^{-1}([c,d)))\leq
	g^C(d)-g^C(c)=\mu^*_{g^C}([c,d))$.
\end{proof}

The following Corollary is in the line of \cite[Lemma 1]{FerFraTo2020}.
\begin{cor} \label{corocambiovar} Given a function $f\in \mathcal{L}^1_g([a,b);{\mathbb F})$, for every $t \in [a,b]$,
	\begin{equation}
		\int_{[a,t)} f\, \operatorname{d}\mu_g=\int_{[a,t)} f\, \operatorname{d}\mu_{g^C}+\sum_{s \in [a,t)\cap D_g}
		f(s)\Delta^+ g(s)=\int_{[0,g^C(t))} \widehat{f} \operatorname{d}\mu +
		\sum_{s \in [a,t)\cap D_g}
		f(s)\Delta^+ g(s),\end{equation}
	where $\widehat{f}=f\circ \gamma$ and $\gamma:t\in [0,g^C(b)]\to
	\gamma(t)$ is given by \eqref{pseudoinverse}.
\end{cor}

\begin{proof} We write $(X,\Sigma_X)=([a,t),\mathcal{M}_{g^C})$ and
	$(Y,\Sigma_Y)=([0,g^C(t)),\mathcal{L})$. We have, thanks to Proposition~\ref{morfismo},
	that $\widehat{f}:Y\to {\mathbb F}$ is a measurable function and
	$g^C:(X,\Sigma_X) \to (Y,\Sigma_Y)$ is a measurable morphism, so
	(cf. \cite[Exercise 1.4.38]{Tao2011}) ensures that
	\begin{equation}
		\int_Y \widehat{f} \operatorname{d} g^C_* \mu_{g^C}=\int_X (\widehat{f}\circ g^C)\, \operatorname{d} \mu_{g^C},\end{equation}
	where
	\begin{equation}
		g^C_* \mu_{g^C}:E\in \Sigma_Y \to g^C_* \mu_{g^C}(E)=\mu_{g^C}((g^C)^{-1}(E))\end{equation}
	is the \emph{pushforward measure} in $(Y,\Sigma_Y)$. However, given an element
	$(c,d)\subset [0,g^C(t))$, it is clear that $g^C_* \mu_{g^C}(c,d)=\mu_{g^C}((g^C)^{-1}(c,d))=d-c$. In
	particular, (cf. \cite[Theorem~13.8]{Bartle1995}) $g^C_* \mu_{g^C}=\mu$. Therefore,
	\begin{equation}
		\int_Y \widehat{f} \operatorname{d} g^C_* \mu_{g^C}=\int_{[0,g^C(t))} \widehat{f} \operatorname{d} \mu.\end{equation}
	Finally, since $\mu_{g^C}(C_{g^C}\cup N_{g^C}^+)=0$ and $\gamma(g^C(s))=s$ for all
	$ s\in [a,b]\backslash( C_{g^C}\cup N_{g^C}^+)$, we have that
	\begin{equation}
		\int_{[a,t)} (\widehat{f}\circ g^C)\, \operatorname{d}\mu_{g^C} = \int_{[a,t)} f\, \operatorname{d}\mu_{g^C}.\end{equation}
\end{proof}

Finally, we recall a concept of continuity introduced in \cite{FriLo17} as well as some of its properties. To that end, denote by $\tau_u$ the usual topology of ${\mathbb F}$ and define the \emph{$g$-topology}, $\tau_g$, as the family of those sets $U\subset \mathbb{R}$ such that for every $x\in U$ there exists $\delta>0$ such that if $y\in \mathbb{R}$ satisfies $|g(y)-g(x)|<\delta$ then $y\in U$. Then, the following definition can be understood as the continuity of a function $f:(I,\tau_g)\to({\mathbb F},\tau_u)$, see \cite[Lemma~6]{MaTo19}.

\begin{dfn}[$g$-continuous function]\label{dfncont} A function $f:[a,b]\to {\mathbb F}$ is
\emph{$g$-continuous} at a point $t\in [a,b]$,
or \emph{continuous with respect to $g$} at $t$, if for every $\varepsilon>0$, there exists $\delta>0$ such that $|f(t)-f(s)|<\varepsilon$,
for every $s \in [a,b]$ with $|g(t)-g(s)|<\delta$. If $f$ is $g$-continuous at every point $t\in [a,b]$, we say that $f$ is $g$-continuous on $[a,b]$.
\end{dfn}

%It can be checked that a function $f:(I,\tau_g)\to({\mathbb F},\tau_u)$ is continuous if and only if it is $g$-continuous \cite[Lemma 6]{MaTo19}.

\begin{pro}[{\cite[Proposition~3.2]{FriLo17}}]\label{proreg} If $f:[a, b] \to \mathbb{R}$ is $g$-continuous on $[a, b],$ then
\begin{enumerate}
	\item $f$ is continuous from the left at every $t_{0} \in(a, b]$;
	\item if $g$ is continuous at $t_{0} \in[a, b),$ then so is $f$;
	\item if $g$ is constant on some $[\alpha, \beta] \subset[a, b],$ then so is $f$.
\end{enumerate}
In particular, $g$-continuous functions on $[a, b]$ are continuous on $[a, b]$ when $g$ is continuous on $[a, b)$.
\end{pro}

\section{The space of bounded $g$-differentiable functions}\label{SBCg}

In the literature --see, for instance, \cite{FriLo17,FriTo20,LoRo14,MaTo19}-- authors
use the following definition of Stieltjes derivative.
\begin{dfn}\label{Stieltjesderivative}
	 We define the \emph{Stieltjes derivative}, or \emph{$g$--derivative}, of function $f:[a,b]\to\mathbb R$ at a point {$t\in [a,b]\backslash C_g$} as
	\begin{equation}%
		f'_g(t)=\left\{
		\begin{array}{rl}
			\displaystyle \lim_{s \to t}\frac{f(s)-f(t)}{g(s)-g(t)},\quad & t\not\in D_{g},\vspace{0.1cm}\\
			\displaystyle\lim_{s\to t^+}\frac{f(s)-f(t)}{g(s)-g(t)},\quad & t\in D_{g},
		\end{array}
		\right.
	\end{equation}
	provided the corresponding limits exist and, in that case, we say that $f$ is
	\emph{$g$--differentiable at $t$}. In particular, for $t\in N_g^+\cup N_g^-$,
	the $g$-derivative at $t$ must be understood in the following sense:
	\begin{equation} \label{eq:ngpoints}
		f'_g(t)=\left\{\begin{array}{ll}
			\displaystyle \displaystyle \lim_{s \to t^+}\frac{f(s)-f(t)}{g(s)-g(t)}, & t\in N_g^+, \vspace{0.1cm} \\
			\displaystyle \displaystyle \lim_{s \to t^-}\frac{f(s)-f(t)}{g(s)-g(t)}, & t\in N_g^-.
		\end{array}\right.
	\end{equation}
\end{dfn}
\begin{rem}
	Observe that the points of $C_g$ are excluded from the definition of  $g$--derivative.
	This is because the corresponding limit cannot be considered at those points since they are in a neighborhood where the corresponding function is not defined. Observe also that the  previous definition is also valid for functions with values in ${\mathbb C}$.
\end{rem}

%% EXISTENCIA DE L\'IMITE EN LOS EXREMOS DEL INTERVALO
%\begin{rem} { From now on, given an interval 
%$[a,b]$ and a derivator $g:\mathbb{R}\to \mathbb{R}$ continuous 
%at $x=a$ and $g(a)=0$, we will assume that $a\notin N_g^-\cup C_g$ 
%and $b\notin N_g^+\cup C_g\cup D_g$. This will be in order to guarantee that the derivative can be taken at the points $a$ and $b$. 
%In this case, given a function $f:[a,b]\to \mathbb{F}$:
%\begin{displaymath}
%\begin{aligned}
%f^{\prime}_g(a)=& \lim_{s\to a^+} \frac{f(s)-f(a)}{g(s)-g(a)},\\
%f^{\prime}_g(b)=& \lim_{s\to b^-} \frac{f(s)-f(b)}{g(s)-g(b)},
%\end{aligned}
%\end{displaymath}
%provided the corresponding limits exist.
%}
%\end{rem}

\begin{rem} \label{importante} { Taking into account Definition~\ref{Stieltjesderivative} and given a function 
$f:[a,b]\to \mathbb{R}$, the following conditions will be necessary for the existence of the $g$-derivative in all of the points of  $[a,b]\backslash C_g$:
	\begin{itemize}
		\item If $a\in [a,b]\setminus C_g$, then $a\notin N_g^-$. Indeed if $a\in N_g^-$, to calculate the $g$-derivative at $a$ 
		we need to know the values of $f$ to the left of $a$, which are not defined. Observe that $C_g\cap N_g=\emptyset$ therefore the previous condition is equivalent to $a\notin N_g^-$.
		\item If $b\in [a,b]\setminus C_g$, then $b\notin N_g^+\cup D_g$. Indeed if $b\in N_g^+\cup D_g$, to calculate the $g$-derivative at $b$ 
		we need to know the values of $f$ to the right of $b$, which are not defined. Observe that $C_g\cap N_g=C_g\cap D_g=\emptyset$ therefore 
		the previous condition is equivalent to $b \notin N_g^+\cup D_g$.
		\item There exists $f(t^+)$ for every $t\in (a,b)\cap D_g$ (which is also a sufficient condition for the existence of the $g$-derivative at that point).
		\item Given $t\in (a,b)\cap N_g^-$ and $\varepsilon>0$, there exists $\delta>0$ such that, if $s<t$ with $g(t)-g(s)<\delta$ then, $|f(s)-f(t)|<\varepsilon$. We say, in that case, that $f$ is \emph{$g$-continuous from the left at $t$}. To check this fact it is enough to observe that $g$ is left continuous (in the usual sense) at $t$. The function $f$ might not be $g$-continuous at $t$. Indeed, take for instance
		\begin{equation}\label{gremark}
			g:x\in [0,3]\to g(s)=\left\{\begin{array}{ll}
				x, & x\in [0,1], \\
				1, & x\in [1,2], \\
				x-1, & x\in [2,3].
			\end{array}
			\right.\end{equation}
		Then,
		\begin{equation}
			f:x\in [0,3]\to f(x)=\left\{\begin{array}{ll}
				x, &x\in [0,1], \\
				x+1, & x\in (1,3],
			\end{array}\right.\end{equation}
		is $g$-differentiable at $x=1$ since
		\begin{equation}
			\lim_{s\to 1^-} \frac{f(s)-f(1)}{g(s)-g(1)}=\lim_{s\to 1^-} \frac{s-1}{s-1}=1.\end{equation}
		Observe that $g$ is continuous at $x=1$, but $f$ is not, so $f$ cannot be $g$-continuous at that point.
		\item Given $t\in (a,b)\cap N_g^+$, and $\varepsilon>0$, there exists $\delta>0$ such that, if $s>t$ with $g(s)-g(t)<\delta$ then, $|f(s)-f(t)|<\varepsilon$. We say, in that case, that $f$ is \emph{$g$-continuous from the right at $t$}. To check this fact it is enough to observe that $g$ is right continuous (in the usual sense) at $t$. Observe that, once again, the $f$ might not be $g$-continuous at such points. Indeed, take for instance $g$ as in \eqref{gremark} and
		\begin{equation}
			f:x\in [0,3]\to f(x)=\left\{\begin{array}{ll}
				x, &x\in [0,2), \\
				x+1, & x\in [2,3].
			\end{array}\right.\end{equation}
		In this case, $f$ is $g$-differentiable at $x=2$ but $f$ is not $g$-continuous at such point.
		
		\item Given $t\in (a,b)\backslash(C_g \cup D_g \cup N_g)$, $f$ is $g$-continuous at $t$. In particular, $f$ is continuous at $t$ since $g$ is continuous at those points.		
	\end{itemize}
	We conclude that, interestingly enough, the $g$-differentiability of a function  at a point of $N_g$ does not imply the $g$-continuity of the function at the point. The $g$-differentiability of a function only guarantees the $g$-continuity at the points of $(a,b)\backslash(C_g \cup D_g \cup N_g)$.}
\end{rem}

%We wonder if it is possible to define the Stieltjes derivative at $C_g$ so
%that the $g$-derivative function is continuous under certain conditions.

{\begin{dfn}[{$\mathcal{C}^1_g([a,b];{\mathbb F})$ space}] \label{c1g}
	Let $g:\mathbb R\to\mathbb R$ be such that $a\notin N_g^-$ and $b \notin N_g^+\cup D_g$. We say that $f:[a,b]\to {\mathbb F}$ belongs to $\mathcal{C}^1_g([a,b];{\mathbb F})$ if the following conditions are met:
	\begin{enumerate}
		\item $f \in \mathcal{C}_g([a,b];{\mathbb F})$,
		\item $\exists f_g'(x)$, for every $x\in [a,b]\backslash C_g$,
		\item $\exists h \in \mathcal{C}_g([a,b];{\mathbb F})$ such that $h(x)=f_g'(x)$, for every $x\in [a,b]\backslash C_g$.
	\end{enumerate}
	Unless necessary, we will write $\mathcal{C}_g^1([a,b])$ instead of  $\mathcal{C}_g^1([a,b];{\mathbb F})$ for brevity.
\end{dfn}}

{ Let us show now that if we assume that $b\notin C_g$ the previous definition is consistent insofar as the function given by 3, if it exists, it is unique.

\begin{pro} Let $[a,b]\subset \mathbb{R}$ be a closed interval,
$g:\mathbb{R} \rightarrow \mathbb{R}$ a derivator such that 
$a\notin N_g^-$ and $b\notin C_g\cup N_g^+ \cup D_g$ and $f\in \mathcal{C}_g([a,b];\mathbb{F})$ be $g$-differentiable at  every $x\in [a,b]\backslash C_g$.
	If $h_1,\, h_2 \in \mathcal{C}_g([a,b])$ are such that
	$h_1(x)=h_2(x)=f_g'(x)$, for every $x\in [a,b]\backslash C_g$, then $h_1=h_2$.
\end{pro}

\begin{proof} Let us show that $h_1(x)=h_2(x)$
	for every $x\in C_g$. Given $\widetilde{x}\in C_g$, there exists a unique connected component of $C_g$, $(a_n,b_n)$, such that $\widetilde{x}\in (a_n,b_n)$. Let us see that $h_1(\widetilde{x})
	=h_2(\widetilde{x})=f_g'(b_n)$. Indeed, since $h_1$ is $g$-continuous, we have, by Proposition~\ref{proreg}, that $h_1$ is constant on $(a_n,b_n)$ and left-continuous, therefore, $h_1(\widetilde x)=h_1(b_n)=f_g'(b_n)$ for every $x\in (a_n,b_n)$. The case of $h_2$ is proven analogously.
\end{proof}}

{\begin{rem} Observe that if $b\in C_g$, given a 
function $f\in \mathcal{C}_g^1([a,b])$ the function $h\in 
\mathcal{C}_g([a,b])$ such that $f_g'(x)=h(x)$ for every 
$x\in [a,b]\setminus C_g$ is not uniquely defined in 
a neighborhood of point $b$ since we can not compute 
the $g$-derivative at $x=b_n$, with $b\in (a_n,b_n)\subset C_g$.
\end{rem}}

A possible way of defining the $g$-derivative at the points of $C_g$, which is coherent with the definition of the space ${\mathcal C}_g^1$, follows from the previous proof. Indeed, we can generalize Definition~\ref{Stieltjesderivative} in the following terms.

 \begin{dfn} \label{Stieltjesderivative2} We define the \emph{Stieltjes derivative}, or \emph{$g$--derivative}, of a function $f:[a,b]\to\mathbb F$ at a point $t\in [a,b]$ as
	\begin{equation} \label{eq:gder}
		f'_g(t)=\left\{
		\begin{array}{rl}
			\displaystyle \lim_{s \to t}\frac{f(s)-f(t)}{g(s)-g(t)},\quad & t\not\in D_{g}\cup C_g,\vspace{0.1cm}\\
			\displaystyle\lim_{s\to t^+}\frac{f(s)-f(t)}{g(s)-g(t)},\quad & t\in D_{g}, \vspace{0.1cm}\\
			\displaystyle \lim_{s \to b_n^+}\frac{f(s)-f(b_n)}{g(s)-g(b_n)},\quad & t \in (a_n,b_n)\subset C_g,
		\end{array}
		\right.
	\end{equation}
	with $a_n,b_n$ as in \eqref{Cgdecomp}; provided the corresponding limits exist. In that case, we say that 
	$f$ is \emph{$g$--differentiable at $t$}. The $g$-derivative in 
	the points $N_g$ must be understood as in~\eqref{eq:ngpoints}.
\end{dfn}
\medskip

We have the following result as a consequence of the 
previous definition and \cite[Proposition~2.5]{Ma21}.

\begin{pro}\label{reglaproducto} Let $f_1,f_2:[a,b]\rightarrow \mathbb{F}$ be two $g$-differentiable functions at a point $x$. Then:
	\begin{itemize}
		\item The function $\lambda_{1} f_{1}+\lambda_{2} f_{2}$ is $g$-differentiable at $x$ for any $\lambda_{1}, \lambda_{2} \in \mathbb{R}$ and
		\begin{equation}%
			\left(\lambda_{1} f_{1}+\lambda_{2} f_{2}\right)_{g}^{\prime}(x)=\lambda_{1}\left(f_{1}\right)_{g}^{\prime}(x)+\lambda_{2}\left(f_{2}\right)_{g}^{\prime}(x).
		\end{equation}
		\item The product $f_{1} f_{2}$ is $g$-differentiable at $x$ and
		\begin{equation}%
			\left(f_{1} f_{2}\right)_{g}^{\prime}(x)=\left(f_{1}\right)_{g}^{\prime}(x) f_{2}(x)+\left(f_{2}\right)_{g}^{\prime}(x) f_{1}(x)+\left(f_{1}\right)_{g}^{\prime}(x)\left(f_{2}\right)_{g}^{\prime}(x) \Delta^+ g(x).
		\end{equation}
		\item If $f_2(x)\,(f_2(x)+(f_2)^{\prime}_g(x)\, \Delta^+g(x))\neq 0$, the quotient $f_1/f_2$ is $g$-differentiable at $x$ and
		\begin{equation}%
			\left(\frac{f_1}{f_2}\right)'_g(x)=\frac{\left(f_{1}\right)_{g}^{\prime}(x) f_{2}(x)-
				\left(f_{2}\right)_{g}^{\prime}(x) f_{1}(x)}{f_2(x)\,(f_2(x)+(f_2)^{\prime}_g(x)\, \Delta^+g(x))}.
		\end{equation}
	\end{itemize}
\end{pro}

\begin{dfn} Let $g:\mathbb R\to\mathbb R$ be such that $a\notin N_g^-$ and $b\notin C_g\cup N_g^+\cup D_g$. Given $k\in\mathbb N$, we define $\mathcal{C}^0_g([a,b];{\mathbb F})=\mathcal{C}_g([a,b];{\mathbb F})$ and  $\mathcal{C}^k_g([a,b];{\mathbb F})$ recursively as
	\begin{equation}
		\mathcal{C}^k_g([a,b])=\{f \in \mathcal{C}^{k-1}([a,b];{\mathbb F}):\; (f_g^{(k-1)})'_g\in
		\mathcal{C}_g([a,b];{\mathbb F})\},\end{equation}
	where $f^{(0)}_g=f$ and $f^{(k)}_g=(f^{(k-1)}_g)'_g$, $k\in\mathbb N$.
	We also define $\mathcal{C}^\infty_g([a,b];{\mathbb F}):=\bigcap_{{n\in{\mathbb N}}}\mathcal{C}^k_g([a,b];{\mathbb F})$. Unless necessary, we will write $\mathcal{C}_g^k([a,b])$ instead of  $\mathcal{C}_g^k([a,b];{\mathbb F})$ for brevity.
\end{dfn}

Now we endow $\mathcal{C}^k_g([a,b])$ with a normed space structure. First, observe that $g$-continuous functions on $[a,b]$ are not necessarily bounded \cite[Example 3.3]{FriLo17}, so we will restrict ourselves to the space
$\mathcal{BC}_g([a,b])$ of bounded $g$-continuous functions. This is a Banach space \cite[Theorem~3.4]{FriLo17}  with the supremum norm
\begin{equation}
	\|f\|_0=\sup\{|f(x)|:\; x \in [a,b]\}.\end{equation}

\begin{dfn} Let $g:\mathbb R\to\mathbb R$ be such that $a\notin N_g^-$ and $b\notin C_g\cup N_g^+\cup D_g$. We define:
	\begin{equation}
		\mathcal{BC}_g^1([a,b];{\mathbb F}):=\{f\in \mathcal{C}_g^1([a,b];{\mathbb F}):\; f,\, f'_g \in \mathcal{BC}_g([a,b];{\mathbb F})\}.\end{equation}
	Analogously, given $k\in \mathbb{N}$,
	\begin{equation}
		\mathcal{BC}_g^k([a,b];{\mathbb F})=\{f \in \mathcal{C}_g^k([a,b];{\mathbb F}):\; f^{(n)}_g\in \mathcal{BC}_g([a,b];{\mathbb F}),\;
		\forall n=0,\ldots,k\}\end{equation}
		and we will denote by $\mathcal{BC}_g^0([a,b];{\mathbb F})=
		\mathcal{BC}_g([a,b];{\mathbb F}).$
\end{dfn}

{ In the following results we will assume that 
$[a,b]\subset \mathbb{R}$ and $g:\mathbb{R}\to \mathbb{R}$ is 
a derivator such that $a\notin N_g^-$ and 
$b\notin C_g\cup N_g^+\cup D_g$. We have that $\mathcal{BC}_g^k([a,b])\equiv\mathcal{BC}_g^k([a,b];{\mathbb F})$ is a normed vector space with the norm
\begin{equation}\label{eq:norm}
	\begin{array}{rcl}
		\displaystyle
		\|\cdot\|_k: \mathcal{BC}_g^k([a,b]) & \mapsto & {\mathbb R}\\
		f & \mapsto & \displaystyle \|f\|_k=
		\sum_{0\leq i\leq k} \|f_g^{(i)}\|_0
\end{array}\end{equation}}

Before proving it is also a Banach space, we will present the following Lemma.

\begin{lem} \label{lemabs1} We have the continuous embedding $\mathcal{BC}_g^1([a,b])\hookrightarrow
	\mathcal{AC}_g([a,b])$. Furthermore, for every $f\in \mathcal{BC}_g^1([a,b])$,
	\begin{equation}
		f(x)=f(a)+\int_{[a,x)} f'_g(s)\, \operatorname{d}\mu_g,\; \forall x \in [a,b].\end{equation}
\end{lem}

\begin{proof} This is an immediate consequence of
	\cite[Theorem~6.2]{LoRo14} and \cite[Corollary 6.3]{LoRo14}. Indeed, given $f\in \mathcal{BC}_g^1([a,b])$, it is clear that
	$f\in \mathcal{BC}_g([a,b])$. In particular, $f$ is continuous from the left at the points in $(a,b]\cap D_g$ and constant at the intervals where $g$ is. On the other hand,
	$f'_g\in \mathcal{BC}_g([a,b])\subset \mathcal{L}^1([a,b])$. Hence,
	by the aforementioned results, $f\in \mathcal{AC}_g([a,b])$ and,
	furthermore,
	\begin{equation}
		f(x)=f(a)+\int_{[a,x)} f'_g(s)\, \operatorname{d}\mu_g,\; \forall x \in [a,b].\end{equation}
\end{proof}

From the previous Lemma  we derive the following.
\begin{lem} \label{limintdetcont}Let $h\in \mathcal{BC}_g([a,b])$  and consider the function
	\begin{equation}
		H:x\in [a,b]\to H(x)=\int_{[a,x)} h(s)\, \operatorname{d}\mu_g.\end{equation}
	We have that $H'_g(x)=h(x)$, for every $x\in [a,b]$ and, therefore,
	$H\in \mathcal{BC}_g^1([a,b])$.
\end{lem}

\begin{proof} Indeed, on the one hand, given that $h\in \mathcal{BC}_g([a,b])\subset
	\mathcal{L}^1_g([a,b))$, it holds that $H\in \mathcal{AC}_g([a,b])$, so it is enough to prove that $H_g'(x)=h(x)$ for every $x\in [a,b]$ to get the result.  We study three different cases:
	\begin{itemize}
		\item For $x\in D_g$, it is clear that
		\begin{equation}
			\begin{array}{rcl}
				H_g'(x)&=&
				\displaystyle
				\lim_{s\to x^+} \frac{H(s)-H(x)}{g(s)-g(x)}  \\
				&=& \displaystyle
				\lim_{s\to x^+} \frac{1}{g(s)-g(x)} \int_{[x,s)} h(s)\, \operatorname{d}\mu_g  \\
				&=& \displaystyle \lim_{s\to x^+} \frac{1}{g(s)-g(x)}
				\left(  \int_{\{x\}} h(s)\, \operatorname{d}\mu_g + \int_{(x,s)} h(s)\, \operatorname{d}\mu_g
				\right)   \\
				&=& \displaystyle \lim_{s\to x^+} \frac{ h(x) \Delta^+g(x)}{g(s)-g(x)}=h(x).
		\end{array}\end{equation}
		\item For $x \in [a,b]\backslash ({C_g}\cup D_g)$, let us compute the limit
		\begin{equation}
			\lim_{s\to x} \frac{H(s)-H(x)}{g(s)-g(x)},\end{equation}
		on the domain of the function, namely, $D_x=\{s\in[a,b]: g(s)\neq g(t)\}$. Fix $\varepsilon>0$. Since $h$ is $g$-continuous and $g$ is continuous at $x$, there exists $\delta>0$ such that $|h(u)-h(x)|<\varepsilon$ if $|u-x|<\delta$. Define $\llbracket x,s \rrparenthesis:=[\min\{x,s\},\max\{x,s\})$. Now, for $s\in[a,b]\cap D_x$, $|u-s|<\delta$, we have that
		\begin{equation}
			\begin{array}{rcl}
				\displaystyle
				\left|\frac{H(s)-H(x)}{g(s)-g(x)} -h(x)\right|&=&
				\displaystyle \left| \frac{\operatorname{sgn}(s-x)}{g(s)-g(x)}
				\int_{\llbracket x,s \rrparenthesis} h(u)\, \operatorname{d}\mu_g(u)-h(x)\right| \vspace{0.1cm} \\
				&=& \displaystyle
				\frac{1}{|g(s)-g(x)|} \left| \int_{\llbracket x,s \rrparenthesis} \left(h(u)-h(x) \right)\, \operatorname{d}\mu_g(u) \right|
				\vspace{0.1cm} \\
				&\leq & \displaystyle
				\frac{1}{|g(s)-g(x)|} \int_{\llbracket x,s \rrparenthesis} |h(u)-h(x)|\, \operatorname{d}\mu_g(u)  \\
				&\leq & \displaystyle
				\frac{1}{|g(s)-g(x)|} \int_{\llbracket x,s \rrparenthesis}\varepsilon\, \operatorname{d}\mu_g(u) =\varepsilon.
		\end{array}\end{equation}
		Thus,
		\begin{equation}%
			\lim_{s\to x} \frac{H(s)-H(x)}{g(s)-g(x)}=h(x).\end{equation}
		%\item The proof in the case $ x \in N_g=N_g^-\cup N_g^+$ is similar to the previous one.
		\item Finally, for $x\in(a_n,b_n)\subset C_g$, it holds that
		\begin{equation}
			H'_g(x)=H_g'(b_n)=h(b_n)=h(x),\end{equation}
		where the first equality comes from the definition of the $g$-derivative at the points of
		$C_g$ and the last is a consequence of the $g$-continuity of $h$.\qedhere
	\end{itemize}
\end{proof}

\begin{thm} \label{Banach1} $(\mathcal{BC}_g^k([a,b]),\|\cdot\|_k)$ is a Banach space.
\end{thm}

\begin{proof} Let us check the case $k=1$ (the case $k\geq 2$ is analogous). Let $\{f_n\}_{n \in \mathbb{N}}\subset \mathcal{BC}_g^1([a,b])$
	be a Cauchy sequence. Then, $\{f_n\}_{n \in \mathbb{N}}\subset \mathcal{BC}_g([a,b])$ and
	$\{(f_n)'_g\}_{n \in \mathbb{N}}\subset \mathcal{BC}_g([a,b])$ are Cauchy sequences in the Banach space $\mathcal{BC}_g([a,b])$ so there exist $f,\,h\in \mathcal{BC}_g([a,b])$ such that
	$f_n \to f$  and $(f_n)'_g \to h$ in $\mathcal{BC}_g([a,b])$. Let us check that
	$f'_g(x)$ exists for every $x\in [a,b]$ and that, furthermore, $f_g'=h$. Indeed, let $\varepsilon>0$. Since $(f_n)'_g \to h$, there exists $N\in\mathbb N$ such that $\|(f_n)'_g-h\|_0\leq \varepsilon/(g(a)-g(b))$. Now, using Lemma~\ref{lemabs1}, we have that
	\begin{equation}
		f_n(x)-f_n(a)=\int_{[a,x)} (f_n)'_g(s)\, \operatorname{d}\mu_g,\; \forall x \in [a,b],\end{equation}
	whence, for $n\ge N$,
	\begin{equation}
	 \left|\int_{[a,x)} (f_n)'_g(s)\, \operatorname{d}\mu_g- \int_{[a,x)} h(s)\, \operatorname{d}\mu_g
			\right|  \leq \int_{[a,x)} \left|(f_n)'_g(s)-h(s) \right|\, \operatorname{d}\mu_g \le \varepsilon,\quad \forall x\in [a,b].
	\end{equation}
	This means that
	\begin{equation}
		\lim_{n \to \infty} \int_{[a,x)} (f_n)'_g(s)\, \operatorname{d}\mu_g=
		\int_{[a,x)} h(s)\, \operatorname{d}\mu_g\end{equation}
	uniformly on $[a,b]$. Thus,
	\begin{equation}
		\lim_{n \to \infty} (f_n(x)-f_n(a))=\lim_{n \to \infty} \int_{[a,x)} (f_n)'_g(s)\, \operatorname{d}\mu_g=
		\int_{[a,x)} h(s)\, \operatorname{d}\mu_g\end{equation}
	uniformly on $[a,b]$. Hence,
	\begin{equation}
		f(x)=f(a)+\int_{[a,x)} h(s)\, \operatorname{d}\mu_g.\end{equation}
	Since $h \in \mathcal{BC}_g([a,b])$, by Lemma~\ref{limintdetcont}, we get that $f'_g(x)=h(x)$ for all $x\in[a,b]$ as we wanted to show.
%	\begin{equation}
%		f_g'(x)=h(x),\; \forall x \in [a,b],\end{equation}
%	as we wanted to show.
\end{proof}

Let us study now the properties of the functions in 
$\mathcal{BC}_g^{n}([a,b];\mathbb{F})$.

\begin{rem}\label{rempf} Observe that, given $f_1,\,f_2\in \mathcal{BC}^1_g([a,b])$, the product $f_1\, f_2\not\in\mathcal{BC}^1_g([a,b])$ in general. This happens because the product
	$\left(f_{1}\right)_{g}^{\prime}(x)
	\left(f_{2}\right)_{g}^{\prime}(x) \Delta^+ g(x)$ might not be $g$-continuous. Indeed, take the following derivator
	\begin{equation}
		g:t \in [-1,1]\to g(t)=\left\{
		\begin{array}{ll}
			t,& t \in [-1,0], \vspace{0.1cm} \\
			t+2, & t \in (0,1].
		\end{array}\right.\end{equation}
	and the function
	\begin{equation}
		f:t\in [-1,1]\to f(t)=\left\{
		\begin{array}{ll}
			t, & t \in [-1,0], \vspace{0.1cm} \\
			2, & t\in (0,1].
		\end{array}\right.\end{equation}
	It is easy to check that $f\in\mathcal{BC}_g([a,b])$ and its derivative,
	\begin{equation}
		f'_g:t \in [-1,1]\to f'_g(t)=\left\{\begin{array}{ll}
			1, & t\in [-1,0), \vspace{0.1cm} \\
			\displaystyle \frac{f(0^+)-f(0)}{\Delta^+g(0)}=1,& t=0, \vspace{0.1cm} \\
			0, & t \in (0,1],
		\end{array}\right.\end{equation}
	is also $g$-continuous. On the other hand,
	\begin{equation}
		f^2:t\in [-1,1]\to f^2(t)=\left\{
		\begin{array}{ll}
			t^2, & t \in [-1,0], \vspace{0.1cm} \\
			4, & t\in (0,1],
		\end{array}\right.\end{equation}
	is $g$-continuous, but
	\begin{equation}
		(f^2)'_g:t \in [-1,1]\to (f^2)'_g(t)=\left\{\begin{array}{ll}
			2t, & t \in [-1,0), \vspace{0.1cm} \\
			\displaystyle \frac{f^2(0^+)-f^2(0)}{\Delta^+g(0)}=2, & t=0,\vspace{0.1cm} \\
			0, & t \in (0,1],
		\end{array}\right.\end{equation}
	is not $g$-continuous since $\lim_{t\to 0^-} (f^2)'_g(t)=0\neq 2=(f^2)'_g(0)$. The problem, as mentioned before, relies on the  $g$-continuity (or lack thereof) of the term $\left(f\right)_{g}^{\prime}(x)
	\left(f\right)_{g}^{\prime}(x) \Delta^+ g(x)$. Indeed, given an element $t \in [-1,1]$, we have that
	\begin{equation}
		(f^2)'_g(t)=2 f_g'(t)\,f(t)+(f_g'(t))^2 \Delta^+g(t).\end{equation}
	Now, even though
	\begin{equation}
		\lim_{t\to 0^-}2 f_g'(t)\,f(t) = 2 f_g'(0)\, f(0),\end{equation}
	the same does not happen for
	\begin{equation}
		\lim_{t\to 0^-} (f_g'(t))^2 \Delta^+g(t) =0 \neq (f_g'(0))^2 \Delta^+g(0)=2.\end{equation}
	The problem~in the lack of continuity of the previous term can be solved if one of the functions involved in the product is also continuous at the points of the discontinuity of the derivator, that is, if it is also continuous in the usual sense, as the following Proposition~shows.
\end{rem}

\begin{pro}\label{probcmul} 
	Given $f_1 \in \mathcal{BC}_{g}^1([a,b])\cap \mathcal C([a,b])$ and
	$f_2 \in \mathcal{BC}_{g}^1([a,b])$, it holds that
	$f_1 f_2 \in \mathcal{BC}_g^1([a,b])$ and
	\begin{equation}\label{simplifiedproductrule}
		\left(f_{1} f_{2}\right)_{g}^{\prime}(x)=
		\left(f_{1}\right)_{g}^{\prime}(x) f_{2}(x)+
		\left(f_{2}\right)_{g}^{\prime}(x) f_{1}(x)
		,\; \forall x \in [a,b].\end{equation}
	In particular, if $f_1 \in \mathcal{BC}_{g}^1([a,b])\cap \mathcal{BC}_{g^C}([a,b])$ and $f_2 \in \mathcal{BC}_{g}^1([a,b])$, we have that \eqref{simplifiedproductrule} holds.
\end{pro}

\begin{proof} By the continuity of $f_1$, given
	$x\in [a,b]\cap D_g$, we have that $(f_1)^{\prime}_g(x)=0$.
	Hence, Proposition~\ref{reglaproducto}  and the definition of the $g$-derivative at the points of $C_g$ imply that
	\begin{equation}%
		(f_1 f_2)^{\prime}_g(x)=
		\left(f_{1}\right)_{g}^{\prime}(x) f_{2}(x)+
		\left(f_{2}\right)_{g}^{\prime}(x) f_{1}(x)
		,\; \forall x \in [a,b].\end{equation}
	Therefore,
	$(f_1 f_2)^{\prime}_g\in \mathcal{BC}_g([a,b])$ and, thus,
	$f_1 f_2 \in \mathcal{BC}_g^1([a,b])$.
\end{proof}

In the following Corollary, which can be obtained by Proposition~\ref{probcmul} using induction, we provide a generalization of Proposition~\ref{probcmul}.

\begin{cor} \label{regadi} Given $f_1\in \mathcal{BC}_{g}^{n}([a,b])\cap \mathcal{BC}_{g^C}^{n-1}([a,b])$ and
	$f_2 \in \mathcal{BC}_{g}^n([a,b])$, we have that $f_1 f_2 \in \mathcal{BC}_g^n([a,b])$.
\end{cor}

\section{First orden linear Stieltjes differential equations}\label{SFO}

To simplify the notation we will work on the interval $[a,b]=[0,T]$. In this section we will  analyze a first order linear Stieltjes differential
equation where the coefficients and data are complex valued functions. Additionally, we will prove further properties of the solution and we will show that, under some regularity conditions for coefficients and data, it is possible to
obtain solutions in the space $\mathcal{BC}_g^1([0,T];{\mathbb F})$. 
{ In order to correctly define the regular solutions 
in the space $\mathcal{BC}_g^1([0,T];{\mathbb F})$ we will assume 
that $0\notin N_g^-$ and $T\notin N_g^+\cup D_g 
\cup C_g$. This consideration is not necessary when looking for} solutions in the space of the absolutely continuous functions.

\subsection{The homogeneous case}

Let us consider the first order homogeneous linear problem
\begin{equation} \label{eq:expg}
	\left\{
	\begin{array}{l}
		v_g'(t)-\beta(t)\,v(t)=0,\; g-a.e.\, t \in [0,T),  \\
		v(0)=v_0,
	\end{array}
	\right.\end{equation}
where $\beta \in \mathcal{L}^1_g([0,T];{\mathbb F})$ and $v_0\in{\mathbb F}$. The solution of problem
\eqref{eq:expg} was given, for the first time, in \cite{FriLo17} for the real case.
In this section we will analyze the existence of solution in the complex case and see how to recover the particular
cases studied in \cite{FriLo17}.  Apart from the generalization proposed here for the complex case, we present a
constructive proof of the expression of the solution which brings light to the nature of the structure of the
solutions of problem~\eqref{eq:expg}.

For the work ahead, we will need to use the chain rule for the Stieltjes derivative. In \cite[Proposition~3.15]{MarquezTesis} we can find a version of the result for the derivative of real valued functions at a continuity point of the derivator. Here, we introduce the following more general version.

\begin{pro} \label{propcompo} Let $x_{0} \in[0, T]$,  $f:[0, T] \to \mathbb{R}$ and $h:\mathbb R\to\mathbb F$. Then, the following hold:
	\begin{enumerate}
		\item If $x_{0} \in[0, T] \backslash (C_g\cup D_{g})$ and there exist $h^{\prime}(f(x_{0}))$ and $f_{g}^{\prime}(x_{0})$, then $h \circ f$ is $g$-differentiable at $x_{0}$ and
		\begin{equation}\label{chainrule}
			(h \circ f)_{g}^{\prime}(x_{0})=h^{\prime}(f(x_{0})) f_{g}^{\prime}(x_{0}).
		\end{equation}
	\item If $x_0\in C_g$ and there exist $h^{\prime}(f(x_{0}))$ and $f_{g}^{\prime}(x_{0})$, then $h \circ f$ is $g$-differentiable at $x_{0}$ and
	\begin{equation}
		(h \circ f)_{g}^{\prime}(x_{0})=h^{\prime}(f(b_n)) f_{g}^{\prime}(x_{0}),
	\end{equation} 
	with $b_n$ the corresponding element of \eqref{Cgdecomp} in \eqref{eq:gder}.
	\item If $x_0\in D_g$ and 
	\begin{equation}\label{condchainrule}
		f(s)=f(x_0), \quad s \in (x_0,x_0+\delta) \mbox{ for some }\delta>0,
	\end{equation}
	then 
	$f_g'(x_0)=(h\circ f)'_g(x_0)=0$. In 
	particular, \eqref{chainrule} holds provided $h'(f(x_0))$ exists.
	\item Suppose that $x_0\in D_g$ and condition \eqref{condchainrule} does not hold. If $f(x_0^+)$ exists, $h$ is 
%	defined on a neighborhood of $f(x_0^+)$ and 
	continuous at $f(x_0^+)$ and the limit
	\begin{equation}\label{limitchain}
		\lim_{s\to x_0^+} \frac{h(f(s))-h(f(x_0))}{f(s)-f(x_0)}
	\end{equation}
	exists, then there exist $f_g'(x_0)$ and $(h\circ f)'_g(x_0)$ and
	\begin{equation}%
		(h \circ f)_{g}^{\prime}(x_{0})=\frac{h(f(x_{0}^{+}))-h(f(x_{0}))}{f(x_{0}^{+})-f(x_{0})} f_{g}^{\prime}(x_{0}).
	\end{equation}
	\end{enumerate}
\end{pro}
\begin{proof}
	First, observe that 1 follows directly from \cite[Proposition~3.15]{MarquezTesis}. Now, 2 is a consequence of the definition of the $g$-derivative in $C_g$ together with 1. Noting  that \eqref{condchainrule} guarantees that $f(x_0^+)=f(x_0)$ and $(h\circ f)(x_0^+)=(h\circ f)(x_0)$ is enough to obtain 3. Finally, for 4, the hypotheses ensure that $f'_g(x_0)$ exists and
	\begin{equation}
	\lim_{s\to x_0^+} h\circ f(s)=h(f(x_0^+)),
	\end{equation}
	so $(h\circ f)'_g(x_0)$ also exists. On the other hand, given that \eqref{condchainrule} does not hold, we can find $\{x_n\}_{n\in\mathbb N}\subset[0,T]$ such that $x_n\to x_0$, $x_n>x_0$,  and $f(x_n)\not=f(x_0)$ for all $n\in\mathbb N$. Hence, given \eqref{limitchain}, we have that
	\begin{align*}
		(h\circ f)'_g(x_0)&=\lim_{n\to\infty} \frac{h(f(x_n))-h(f(x))}{g(x_n)-g(x)}\\
		&=\lim_{n\to\infty} \frac{h(f(x_n))-h(f(x))}{f(x_n)-f(x)}\frac{f(x_n)-f(x)}{g(x_n)-g(x)}=\frac{h(f(x_{0}^{+}))-h(f(x_{0}))}{f(x_{0}^{+})-f(x_{0})} f_{g}^{\prime}(x_{0}).
	\end{align*}
\end{proof}

For the following theorem, we will denote by $\ln(z):=\ln |z|+ i \operatorname{Arg}(z)$ for $z\in \mathbb{C}$ the principal branch of the complex logarithm where $\operatorname{Arg}$ is the principal argument.
\begin{thm}\label{solel} Assume $\mu_g(\overline{D_g}\backslash D_g)=0$. Let $\beta \in \mathcal{L}_g^1([0,T),{\mathbb F})$ be such
	that $1+\beta(t)\Delta^+g(t) \neq 0$
	for every $t\in [0,T)\cap D_g$. Then there exists a unique solution $v\in \mathcal{AC}_g([0,T];{\mathbb F})$ of problem~\eqref{eq:expg} which, furthermore, is of the form
	\begin{equation}%
		v(t)=v^B(t)\, v^C(t),
		\end{equation}
	where $v^B\in \mathcal{AC}_{g^B}([0,T];{\mathbb F})$ is the unique solution of the problem
	\begin{equation} \label{eq:expgB}
		\left\{
		\begin{array}{l}
			v_{g^B}'(t)-\beta(t)\,v(t)=0,\; g^B-a.e.\, t \in [0,T),  \\
			v(0)=1,
		\end{array}
		\right.\end{equation}
	given by
	\begin{equation} \label{eq:solexpgB}
		v^B(t)=\displaystyle \prod_{s\in [0,t)\cap D_g}
		\left(1+\beta(s)\Delta^+g(s) \right);\end{equation}
	and $v^C\in \mathcal{AC}_{g^C}([0,T];{\mathbb F})$ is the unique solution of
	\begin{equation} \label{eq:expgC}
		\left\{
		\begin{array}{l}
			v_{g^C}'(t)-\beta(t)\,v(t)=0,\; g^C-a.e.\, t \in [0,T),  \\
			v(0)=v_0,
		\end{array}
		\right.\end{equation}
	given by
	\begin{equation} \label{eq:solexpgC}
		v^C(t)= \displaystyle u(g^C(t)),\end{equation}
	where $u \in \mathcal{AC}([0,T];{\mathbb F})$ is the unique  solution of
	\begin{equation} \label{eq:original2}
		\left\{\begin{array}{l}
			\displaystyle u'(t)=\widehat{\beta}(t)\,u(t),\; a.e.\, t \in [0,g^C(T)),  \\
			\displaystyle u(0)=v_0,
		\end{array} \right.\end{equation}
	where $\widehat{\beta}=\beta\circ \gamma$ and $\gamma$ is provided by Definition~\ref{dfnpi}.

	Furthermore, $v$ can be written as
	\begin{equation} \label{eq:solexpg}
		v(t)=v_0\exp\left(\int_{[0,t)} \widetilde{\beta}(s)\, \operatorname{d}\mu_g\right),\end{equation}
	with
	\begin{equation}\label{eqbtil}
		\widetilde{\beta}(t)=\left\{
		\begin{array}{rl}
			\displaystyle \beta(t),& t \in [0,T)\backslash D_g,  \\
			\displaystyle \frac{\ln\left(1+\beta(t)\Delta^+g(t) \right)}{\Delta^+g(t)}, &
			t \in [0,T)\cap D_g.
		\end{array}\right.\end{equation}
\end{thm}

\begin{proof}
	\emph{Existence and uniqueness:} If $v$ solves \eqref{eq:expg}, then $(x,y)$ where $x:=\operatorname{Re} v$ and $y:=\operatorname{Im} v$ solves the  real system
	\begin{equation}\label{prore}
		\left\{
		\begin{array}{l}
			x_g'(t)-\operatorname{Re}\beta(t)\,x(t)+\operatorname{Im}\beta(t)\,y(t)=0,\; g-a.e.\, t \in [0,T),  \\
			y_g'(t)-\operatorname{Im}\beta(t)\,x(t)-\operatorname{Re}\beta(t)\,y(t)=0,\; g-a.e.\, t \in [0,T),  \\
			x(0)=\operatorname{Re} v_0,\ y(0)=\operatorname{Im}v_0,
		\end{array}
		\right.\end{equation}
	and vice-versa, that is, the function $x+iy$, where $(x,y)$ is a solution of \eqref{prore}, solves \eqref{eq:expg}. Now, it is easy to see that \eqref{prore} satifies the conditions of \cite[Theorem 4.3]{LoMa19} with $L=|\operatorname{Re}(\beta)|+|\operatorname{Im}(\beta)|$, so it has a unique solution on $[0,T]$. Hence, \eqref{eq:expg} has a unique solution there as well.

	\emph{Expression of the solution:} Given the nature of problem~\eqref{eq:expg}
	where the $g$-derivative has to be a multiple of itself it is only natural to use an \emph{antsatz} of the form
	\begin{equation}
		v(t)=v_0\exp\left( \int_{[0,t)}\widetilde{\beta}(s)\, \operatorname{d} \mu_g
		\right),\end{equation}
	with $\widetilde{\beta}\in \mathcal{L}^1_g([0,T);{\mathbb F})$. Given that
	$\mathcal{M}_g \subset \mathcal{M}_{g^C}$ and $\mathcal{M}_g\subset \mathcal{M}_{g^B}$,
	it is clear that if $\widetilde{\beta}\in \mathcal{L}^1_g([0,T);{\mathbb F})$, then $\widetilde{\beta}\in \mathcal{L}_{g^C}^1([0,T);{\mathbb F})$ and
	$\widetilde{\beta}\in \mathcal{L}_{g^B}^1([0,T);{\mathbb F})$. Furthermore,
	\begin{equation}
		\begin{array}{rcl}
			v(t)&=& \displaystyle v_0\exp\left(
			\int_{[0,t)}\widetilde{\beta}(s)\, \operatorname{d} \mu_{g^B}+
			\int_{[0,t)}\widetilde{\beta}(s)\, \operatorname{d} \mu_{g^C}
			\right)  \\
			&=& \displaystyle
			v_0\exp\left(
			\int_{[0,t)}\widetilde{\beta}(s)\, \operatorname{d} \mu_{g^B} \right)
			\exp\left(
			\int_{[0,t)}\widetilde{\beta}(s)\, \operatorname{d} \mu_{g^C} \right) = v^B(t) \, v^C(t),
	\end{array}\end{equation}
	where $v^B(t):=\exp\left(
	\int_{[0,t)}\widetilde{\beta}(s)\, \operatorname{d} \mu_g^B \right)$ and $v^C(t):=v_0\exp\left(
	\int_{[0,t)}\widetilde{\beta}(s)\, \operatorname{d} \mu_g^C \right)$.
	From the definition we deduce that
	$v^B\in \mathcal{AC}_{g^B}([0,T];{\mathbb F})$ and
	$v^C\in \mathcal{AC}_{g^C}([0,T];{\mathbb F})$. Hence,
	given the $g^B$-continuity of $v^B$, we have that $(v^B)'_g(t)=0$
	for every $t\in [0,T)\backslash (\overline{D_g} \cup C_g)$ and, thanks to the $g^C$-continuity of $v^C$, it holds that $(v^C)'_g(t)=0$
	for every $t\in [0,T)\cap D_g$. Thus, by Proposition~\ref{reglaproducto},
	\begin{equation}
		v'_g(t)=\left\{\begin{array}{rl}
			(v^B)'_g(t)\, v^C(t), & t \in [0,T)\cap D_g, \vspace{0.1cm} \\
			v^B(t)\, (v^C)'_g(t), & g-a.e.\, t \in [0,T) 
			\backslash( \overline{D_g}\cup C_g).
		\end{array}\right.\end{equation}
	This implies that we will have a different equation for each of the components of the solution:
	\begin{align} \label{eq:expgcases}
		(v^B)'_g(t)=&\beta(t)\,v^B(t),\; t\in [0,T)\cap D_g, \\
		\label{eq:expgcases2} (v^C)'_g(t)=&\beta(t)\,v^C(t),\; g-a.e.\, t \in [0,T) \backslash( \overline{D_g}\cup C_g).
	\end{align}

	We will start studying equation \eqref{eq:expgcases}. For $t\in[0,T)\cap D_g$ we have that
	\begin{equation}
		(v^B)'_g(t)=
		\frac{v^B(t^+)-v^B(t)}{\Delta^+ g(t)}=(v^B)'_{g^B}(t).\end{equation}
	Now, if we develop equation
	\eqref{eq:expgcases}:
	\begin{equation} \label{eq:expgB1}
		\frac{v^B(t^+)-v^B(t)}{\Delta^+ g(t)} =
		\beta(t)\, v^B(t)\end{equation}
	and we get that
	\begin{equation}\label{eqvb}
		v^B(t^+)=v^B(t)(1+\beta(t)\, \Delta^+ g(t)),\;  t\in [0,T)\cap D_g.\end{equation}
	%We propose, as a candidate function satisfying equation~\ref{eqvb}
	%\begin{equation}
	%v^B(t):=\prod_{s \in [0,t)}(1+\beta(s)\Delta^+g(s)) \in \mathcal{AC}_{g^B}([0,T];{\mathbb F}).
	%\end{equation}
	%We have to check that $v^B$ is well defined (the product is convergent), that it is absolutely $g^B$ continuous and that it satisfies equation~\ref{eqvb}.
	%
	%\textbullet\ \emph{$v^B$ is well defined:}
	In order to get a solution candidate for equation \eqref{eq:expgcases}, define $h(t)={\ln(1+\beta(t)\Delta^+g(t))}/{\Delta^+g(t)}$ if $t\in D_g$, $h(t)=0$ if $t\not\in D_g$. Then, taking into account that $\mu_{g^B}(t)=0$ for every $t\notin D_g$, we define
	\begin{equation}
		\begin{array}{rcl}
			\displaystyle
			H(t)&:=&\displaystyle \exp\left(
			\int_{[0,t)}h(t)\, \operatorname{d} \mu_{g^B} \right)
			\\&=&\displaystyle
			\exp \left( \int_{[0,t)} \frac{\ln(1+\beta(s)\Delta^+g(s))}{\Delta^+g(s)}\, \operatorname{d} \mu_{g^B}
			\right)  \\
			&=&\displaystyle
			\exp\left( \sum_{s\in [0,t)\cap D_g} \ln(1+\beta(s)\Delta^+g(s))
			\right) \vspace{0.1cm} \\
			&=&
			\displaystyle
			\exp \left( \sum_{s\in [0,t)\cap D_g} \left( \ln\left|1+\beta(s)\Delta^+g(s)\right| +
			i\, \operatorname{Arg}(1+\beta(s)\Delta^+g(s)) \right) \right).
	\end{array}\end{equation}
	To show that $H$ is well defined, let us check that the series \begin{equation}\sum_{s\in [0,t)\cap D_g}  \ln\left|1+\beta(s)\Delta^+g(s)\right|\text{\quad and\quad}\sum_{s\in [0,t)\cap D_g}\operatorname{Arg}(1+\beta(s)\Delta^+g(s))\end{equation} are absolutely convergent. We have that
	\begin{equation}
		\sum_{s\in [0,T)\cap D_g}\left|
		\ln\left|1+\beta(s)\Delta^+g(s)\right|
		\right|  \\
		\displaystyle
		= \sum_{s\in A}\left|
		\ln\left|1+\beta(s)\Delta^+g(s)\right|
		\right|+\sum_{s\in B}\left|
		\ln\left|1+\beta(s)\Delta^+g(s)\right|
		\right|,\end{equation}
	where
	\begin{equation}
		\begin{aligned}
			A= & \displaystyle
			\left\{ s\in [0,T)\cap D_g:\; \left|1+\beta(s)\Delta^+g(s)\right|\geq1\right\}=\left\{ s\in [0,T)\cap D_g:\; \ln\left|1+\beta(s)\Delta^+g(s)\right|\geq0\right\},  \\
			B= & \displaystyle
			\left\{ s\in [0,T)\cap D_g:\; \left|1+\beta(s)\Delta^+g(s)\right|<1\right\}=\left\{ s\in [0,T)\cap D_g:\; \ln\left|1+\beta(s)\Delta^+g(s)\right|<0\right\}.
	\end{aligned}\end{equation}
	In order to bound the sum on $A$ it is enough to take into account that
	$0\le\ln(1+x)\leq x$ for every $x\in [0,\infty)$:
	\begin{displaymath}
	\begin{aligned}
		\sum_{s\in A}\left|
		\ln\left|1+\beta(s)\Delta^+g(s)\right|
		\right| =& \sum_{s\in A}
		\ln\left|1+\beta(s)\Delta^+g(s)\right|\\
		\leq &
		\sum_{s \in A}  \ln \left( 1+ |\beta(s)|
		\Delta^+g(s)|\right)  \leq
		 \sum_{s \in A} |\beta(s)| \Delta^+g(s)<\infty,
		 \end{aligned}
		 \end{displaymath}
	because $\beta \in \mathcal{L}_{g^B}^1([0,T),{\mathbb F})$.

{Now, let us focus on the sum on $B$. For any $s\in B$, taking into account that $1+\beta(s)\Delta^+g(s)\ne 0$, we have that
	\begin{displaymath}
		\begin{array}{rcl}
			\displaystyle
			0<  \left|1+\beta(s)\,\Delta^+g(s)\right|^2&=&[    1+\operatorname{Re}(\beta(s))\,\Delta^+g(s)]^2+ 
			[\operatorname{Im}(\beta(s))\,\Delta^+g(s)]^2  \\
			&=& \displaystyle 1+2\, \operatorname{Re}(\beta(s))\, \Delta^+g(s)+ |\beta(s)\,\Delta^+g(s)|^2<1.
	\end{array}\end{displaymath}
In particular, $2\, \operatorname{Re}(\beta(s))\, \Delta^+g(s)+ |\beta(s)\,\Delta^+g(s)|^2<0$ which yields $\operatorname{Re}(\beta(s))<0$. Now, we can consider the following sets:
	\begin{displaymath}
		\begin{array}{rcl}
			B_1&=&\displaystyle \left\{s\in B:\; 0<1+2\, \operatorname{Re}(\beta(s)) \Delta^+g(s)+ |\beta(s)\Delta^+g(s)|^2<\frac{1}{2}\right\},
			\\
			B_2&=&\displaystyle \left\{s\in B:\;
			\frac{1}{2}\leq
			1+2\, \operatorname{Re}(\beta(s)) \Delta^+g(s)+ |\beta(s)\Delta^+g(s)|^2 <1 \right\}.
	\end{array}\end{displaymath}
	Observe that $B=B_1\cup B_2$. The definition of $B_1$ implies that
	\begin{equation}
		1>2\, |\operatorname{Re}(\beta(s))| \Delta^+g(s)-|\beta(s)\Delta^+g(s)|^2>\frac{1}{2},\; \forall s \in B_1.\end{equation}
	Therefore,
	\begin{equation}
		|\operatorname{Re}(\beta(s))| \Delta^+g(s) > \frac{1}{4},\; \forall s \in B_1.\end{equation}
	Hence, we have that $B_1$ is finite since, otherwise, we would have that $\beta \not \in \mathcal{L}_{g^B}^1([0,T),{\mathbb F})$,  which is a contradiction.
	For the elements in the set $B_2$ we have that:
	\begin{equation}
		\frac{1}{2}\geq 2\, |\operatorname{Re}(\beta(s))| \Delta^+g(s)-|\beta(s)\Delta^+g(s)|^2 >0,\; \forall s \in B_2.\end{equation}
	Thus, if we take into account that $\ln(1/(1-x))\leq 2x$, for every $x\in [0,1/2]$,
	\begin{equation}
		\begin{aligned}
			\displaystyle
			\left|\ln\left|1+\beta(s)\Delta^+g(s)\right| \right|&=
			\displaystyle
			\frac{1}{2}\left| \ln\left(
			1+2\, \operatorname{Re}(\beta(s)) \Delta^+g(s)+ |\beta(s)\Delta^+g(s)|^2
			\right)\right|  \\
			&= \displaystyle \frac{1}{2} \ln \left({1}/\left(1+2\, \operatorname{Re}(\beta(s)) \Delta^+g(s)+ |\beta(s)\Delta^+g(s)|^2\right)\right)  \\
			&= \displaystyle \frac{1}{2} \ln \left({1}/\left(1-\left(2\, |\operatorname{Re}(\beta(s))| \Delta^+g(s)- |\beta(s)\Delta^+g(s)|^2\right)\right)
			\right) \\
			& \leq \displaystyle 2 \left(2\, |\operatorname{Re}(\beta(s))| \Delta^+g(s)- |\beta(s)\Delta^+g(s)|^2\right)  \\
			&\leq \displaystyle 4 |\operatorname{Re}(\beta(s))| \Delta^+g(s).
	\end{aligned}\end{equation}}
	Hence,
	\begin{equation}
		\sum_{s\in B}\left|
		\ln\left|1+\beta(s)\Delta^+g(s)\right| \right|<\infty.\end{equation}
	Let us now bound the term associated with the argument. Taking into account that $|\operatorname{atan}(x)|\leq |x|$ for every $x\in \mathbb{R}$, we have that
	\begin{equation}
		\sum_{s\in [0,T)\cap D_g}\left|
		\operatorname{Arg}(1+\beta(s)\Delta^+g(s))\right| \leq
		\sum_{s\in [0,T)\cap D_g}
		\frac{\left| \operatorname{Im}(\beta(s))\Delta^+g(s) \right|}{\left|1+\operatorname{Re}(\beta(s))\Delta^+g(s)\right|}.\end{equation}
	Let us divide the set $[0,T)\cap D_g$ into the subsets
	\begin{equation}
		\begin{array}{rcl}
			\widetilde B_1&=&\displaystyle \left\{s \in [0,T)\cap D_g\ :\  |\operatorname{Re}(\beta(s))| \Delta^+g(s) > {1}/{2}\right\},
			\\
			\widetilde B_2&=&([0,T)\cap D_g)\backslash \widetilde B_1.
	\end{array}\end{equation}
	Observe that $\widetilde B_1$ must be of finite cardinality. On the other hand, given $t\in \widetilde B_2$,
	\begin{equation}
		\left|1+\operatorname{Re}(\beta(s))\Delta^+g(s)\right|\geq \frac{1}{2}.\end{equation}
	Thus,
	\begin{equation}
		\sum_{s\in B_2}\frac{\left| \operatorname{Im}(\beta(s))\Delta^+g(s) \right|}{\left|1+\operatorname{Re}(\beta(s))\Delta^+g(s)\right|}
		\leq 2\, \sum_{s\in B_2} \left| \operatorname{Im}(\beta(s))\Delta^+g(s) \right| <\infty.\end{equation}
	Hence, we conclude that $H$ is well defined. In order to prove that $H$ is a solution of \eqref{eq:expgcases}, we observe that, given $t\in [0,T)\cap D_g$,
	\begin{align*}
		H(t^+)= &
		\lim_{s \to t^+}
		\exp \left( \int_{[0,s)} \frac{\ln(1+\beta(s)\Delta^+g(s))}{\Delta^+g(s)}\, \operatorname{d} \mu_{g^B}
		\right)
		\\
		\displaystyle
		= &
		\lim_{s\to t^+}
		\exp \left( \int_{[0,t)} \frac{\ln(1+\beta(s)\Delta^+g(s))}{\Delta^+g(s)}\, \operatorname{d} \mu_{g^B}
		+\ln(1+\beta(t)\Delta^+g(t))\right.
		\\ & \left.+
		\int_{(t,s)} \frac{\ln(1+\beta(s)\Delta^+g(s))}{\Delta^+g(s)}\, \operatorname{d} \mu_{g^B}
		\right) \\
		\displaystyle
		= & (1+\beta(t)\Delta^+g(t))
		\exp \left( \int_{[0,t)} \frac{\ln(1+\beta(s)\Delta^+g(s))}{\Delta^+g(s)}\, \operatorname{d} \mu_{g^B}
		\right),
		\\
		\displaystyle
		= & (1+\beta(t)\Delta^+g(t))
		H(t),
	\end{align*}
	so equation \eqref{eqvb} holds and $v^B:=H$ is a solution of \eqref{eq:expgcases}.
	Observe that, given any set $A\subset [0,T)\backslash D_g$, we have
	$A=(A\backslash \overline{D_g})\cup (A\cap (\overline{D_g}\backslash D_g) )$ thus $\mu^*_{g^B}(A)\leq
	\mu^*_{g^B}(A\backslash \overline{D_g})+\mu^*_{g^B}(A\cap (\overline{D_g}\backslash D_g)) \leq
	\mu^*_{g^B}(A\backslash \overline{D_g}) + \mu^*_g(\overline{D_g}\backslash D_g)=0$.
	Therefore $v^B$ satisfies \eqref{eq:expgB} and, moreover,
	\begin{equation}%
		v^B(t)=\exp\left( \sum_{s\in [0,t)\cap D_g} \ln(1+\beta(s)\Delta^+g(s))\right)= \prod_{s\in [0,t)\cap D_g}
		\left(1+\beta(s)\Delta^+g(s) \right).\end{equation}

	Let us now study  equation \eqref{eq:expgcases2}. First, observe that, given an element $t\in [0,T)\backslash (\overline{D_g}\cup C_g)$,
	there exists $\delta>0$ such that $g$ is continuous on $(t-\delta,t+\delta)$. In the case $t\in N_g^-$ we further know
	that $g$ is strictly increasing on the interval $(t-\delta,t]$, and constant on
	$(t,t+\delta)$. In the case $t\in N_g^+$, $g$ would be constant on $(t-\delta,t)$ and strictly increasing on $[t,t+\delta)$. In any case (observe that, if  $t\in N_g^-$ we have to take the limit from the left and in the case
	$t\in N_g^+$ the limit from the right, respectively):
	\begin{equation}
		(v^C)'_g(t)=\lim_{s\to t} \frac{v^C(s)-v^C(t)}{g(s)-g(t)}=
		\lim_{s\to t}  \frac{v^C(s)-v^C(t)}{g^C(s)-g^C(t)}=(v^C)'_{g^C}(t).\end{equation}
	Hence, taking into account that  $\mu_g^*(A)=0\Leftrightarrow \mu_{g^C}^*(A)=0$ for any $A\subset [0,T)\backslash \overline{D_g}$, together with the fact that
	$C_g=C_{g^C}$, we see that equation  \eqref{eq:expgcases2} is equivalent to
	\begin{equation}\label{previouseq}
		\displaystyle (v^C)'_{g^C}(t)=\beta(t)\,v^C(t),\; {g^C}-a.e.\, t \in [0,T) \backslash( \overline{D_g}\cup C_{g^C}).\end{equation}
	Let us observe that $\mu_{g^C}( \overline{D_g}\cup C_{g^C})\leq \mu_{g^C}(\overline{D_g}\backslash D_g)+
	\mu_{g^C}(D_g)+\mu_{g^C}(C_{g^C})=0$, since $\mu_{g^C}(\overline{D_g}\backslash D_g)\leq \mu_{g}(\overline{D_g}\backslash D_g)=0$
	by hypothesis. Therefore, \eqref{previouseq} is equivalent to:
	\begin{equation}\label{previouseq2}
		\displaystyle (v^C)'_{g^C}(t)=\beta(t)\,v^C(t),\; {g^C}-a.e.\, t \in [0,T).\end{equation}
	Now we will see that $v^C(t):=u(g^C(t))$, with $u\in \mathcal{AC}([0,g^C(T)];{\mathbb F})$ the solution of \eqref{eq:original2} satisfies equation~\eqref{previouseq2}. On the one hand, we have that $\widehat{\beta}=\beta \circ \gamma
	\in \mathcal{L}^1([0,g^C(T)];{\mathbb F})$. Indeed, the measurability is a consequence of Proposition~\ref{morfismo}. Now, using a similar argument as the one in the proof of Corollary~\ref{corocambiovar}:
	\begin{equation}
		\int_{[0,g^C(T))} |\widehat{\beta}| \, \operatorname{d} \mu =\int_{[0,T)} |\beta|\, \operatorname{d} \mu_{g^C} \leq
		\int_{[0,T)} |\beta|\, \operatorname{d} \mu_{g}<\infty.\end{equation}
	Thus, \eqref{eq:original2} admits a unique solution
	\begin{equation}
		u(t)=v_0\exp\left(\int_{[0,t)} \widehat{\beta}(s)\, \operatorname{d} \mu\right)\in \mathcal{AC}([0,g^C(T)];{\mathbb F}).\end{equation}
	In particular, $v^C(t)=u(g^C(t))$ is such that $(v^C)'_g(t)=0$ for every $t\in D_g$. Indeed,
	\begin{equation}
		(v^C)'_g(t) = \lim_{s\to t^+} \frac{v^C(s)-v^C(t)}{g(s)-g(t)}=
		\lim_{s\to t^+} \frac{u(g^C(s))-u(g^C(t))}{g(s)-g(t)}=0,\end{equation}
	thanks to the continuity of the composition $u\circ g^C$. On the other hand, since $u \in \mathcal{AC}([0,g^C(T)];{\mathbb F})$ is the solution of \eqref{eq:original2},
	there exists a Lebesgue-null set $N \subset [0,g^C(T)]$ such that
	\begin{equation}
		u'(t)=\widehat{\beta}(t)\, u(t), \; \forall t \in [0,g^C(T)]\backslash N.\end{equation}
	In particular,
	\begin{equation}
		u'(g^C(t))=\widehat{\beta}(g^C(t))\, u(g^C(t)), \; \forall t \in [0,T]\backslash (g^C)^{-1}(N),\end{equation}
	whence, by Proposition~\ref{propcompo},
	\begin{equation}
		(v^C)'_{g^C}(t)=u'(g^C(t))=\widehat{\beta}(g^C(t))\, u(g^C(t)) ,
		\; \forall t \in [0,T]\backslash (g^C)^{-1}(N).\end{equation}
	Taking into account that $\gamma(g^C(t))=t$ for every $t \in [0,T]\backslash
	(C_{g^C}\cup N_{g^C}^+)$, that $\mu_{g^C}(C_{g^C}\cup N_{g^C}^+)=0$ and that $\mu_{g^C}((g^C)^{-1}(N))=0$ (see the proof of Proposition~\ref{morfismo}), we deduce that
	\begin{equation}
		(v^C)'_{g^C}(t)=\beta(t)\, v^C(t),\; g^C\text{-a.e. } t \in [0,T).\end{equation}
	Last, in which respects the solution $v^C$, we have that, using a reasoning similar to the one used in the proof of the Corollary~\ref{corocambiovar}, we have that
	\begin{equation}
		v^C(t)=u(g^C(t))=v_0\exp\left(\int_{[0,g^C(t))} \widehat{\beta}(s)\, \operatorname{d} \mu\right)=
		v_0\exp\left(\int_{[0,t)} \beta(s)\, \operatorname{d} \mu_{g^C}\right).\end{equation}

	Finally, let us check that $v:=v^C\, v^B$ is in the space $\mathcal{AC}_g([0,T];
	{\mathbb F})$. To show this, let us define
	\begin{equation}
		\widetilde{\beta}(t)=\begin{dcases}
			\displaystyle \beta(t),& t \in [0,T)\backslash D_g,  \\
			\displaystyle \frac{\ln\left(1+\beta(t)\Delta^+g(t) \right)}{\Delta^+g(t)}, &
			t \in [0,T)\cap D_g,
	\end{dcases}\end{equation}
	and check that
	\begin{equation}
		v(t)=
		\displaystyle
		v_0\exp\left(\int_{[0,t)} \widetilde{\beta}(s)\, \operatorname{d}\mu_g\right) \\
		=\displaystyle
		v_0\exp\left(\int_{[0,t)\backslash D_g} \widetilde{\beta}(s)\, \operatorname{d} \mu_g +
		\sum_{s \in [0,t)\cap D_g} \widetilde{\beta}(s) \Delta^+g(s) \right).\end{equation}
	Indeed, on the one hand,
	\begin{equation}
		\begin{array}{rcl}
			v(t)&=&\displaystyle v_0\left[\prod_{s\in [0,t)\cap D_g}
			\left(1+\beta(s)\Delta^+g(s) \right)\right]
			\exp\left(\int_{[0,g^C(t))} \widehat{\beta}(s)\, \operatorname{d} \mu\right)
			\\
			&=& \displaystyle
			v_0\exp\left(
			\int_{[0,g^C(t))} \widehat{\beta}(s)\, \operatorname{d} \mu +
			\sum_{s\in [0,t)\cap D_g} \frac{\ln\left(1+\beta(s)\Delta^+g(s) \right)}{\Delta^+g(s)}
			\Delta^+g(s)
			\right).
	\end{array}\end{equation}
	Now, thanks to the fact that  $\mu(D_g)=0$ as it is a countable set, we see that
	\begin{equation}
		\int_{[0,g^C(t))} (\beta \circ \gamma)(s)\, \operatorname{d} \mu =
		\int_{[0,g^C(t))} (\widetilde{\beta} \circ \gamma)(s)\, \operatorname{d} \mu.\end{equation}
	Thus, by Corollary~\eqref{corocambiovar},
	\begin{equation}
		\int_{[0,t)} \widetilde{\beta}(s)\, \operatorname{d} \mu_g=
		\int_{[0,g^C(t))} \widehat{\beta}(s)\, \operatorname{d} \mu +
		\sum_{s\in [0,t)\cap D_g} \frac{\ln\left(1+\beta(s)\Delta^+g(s) \right)}{\Delta^+g(s)}
		\Delta^+g(s).\end{equation}
Finally, it is clear that $\widetilde{\beta}\in \mathcal{L}^1_g([0,T);\mathbb{F})$, therefore $v\in \mathcal{AC}_g([0,T);{\mathbb F})$.		
\end{proof}

\begin{rem} We must take into account the following remarks:

	\emph{1.} If $g^C$ is constant, then the solution of \eqref{eq:expg} is reduced to $v_0\,v^B$ and the hypothesis $\mu_g (\overline{D_g} \backslash D_g) = 0$ is not necessary.

	\emph{2.} The hypothesis $\mu_g (\overline{D_g} \backslash D_g) = 0 $ that appears in the
	statement of Theorem~\ref{solel} has been used to express the solution of \eqref{eq:expg} as the product of the solutions of the
	problems \eqref{eq:expgB} and \eqref{eq:expgC}. This hypothesis is not
	essential to guarantee the existence of a solution of problem
	\eqref{eq:expg}. Even in the case
	$\mu_g (\overline{D_g} \backslash D_g) \ne 0$, we 
	will have \eqref{eq:solexpg}
	is well defined and a valid solution of problem~\eqref{eq:expg}. 
	 Indeed, since $\widetilde{\beta}\in 
	\mathcal{L}^1_g([0,T);\mathbb{F})$, we have that
	\begin{displaymath}
	\left(\int_{[0,t)} \widetilde{\beta}(s)\, \operatorname{d} \mu_g
	\right)^{\prime}_g(t)=\widetilde{\beta}(t),
	\;g\text{-a.e. } t \in [0,T).
	\end{displaymath}
	Therefore, \eqref{chainrule} ensures that
	\begin{equation}\label{expkk1}
	\left(\exp\left(\int_{[0,t)} \widetilde{\beta}(s)\, \operatorname{d}\mu_g\right)\right)^{\prime}_g(t)={\beta}(t) \, 
	\exp\left(\int_{[0,t)} \beta(s)\, \operatorname{d}\mu_g\right),\; g\text{-a.e. } t \in [0,T)\setminus D_g.
	\end{equation}
	Now, given $t\in [0,T)\cap D_g$,
	\begin{displaymath}
	\begin{aligned}
	\lim_{s\to t^+} \exp\left(\int_{[0,s)} \widetilde{\beta}(s)\, \operatorname{d}\mu_g\right)=&\lim_{s\to t^+}
	\exp\left(  \int_{[0,t)} \widetilde{\beta}(s)\, \operatorname{d}\mu_g  + \ln(1+\beta(t)\Delta^+g(t)) +\int_{(t,s)}\widetilde{\beta}(s)\, \operatorname{d}\mu_g
	\right) \\
	=& (1+\beta(t)\Delta^+g(t))\, 
	\exp\left(\int_{[0,t)} \widetilde{\beta}(s)\, \operatorname{d}\mu_g\right),
	\end{aligned}
	\end{displaymath}
	so equation~\eqref{expkk1} is also satisfied for the points of 
	$D_g$.
\end{rem}

\begin{rem} The previous result is a generalization of the results in \cite[Section 6]{FriLo17} for several reasons.

	\emph{1.} The solution obtained is valid in the complex case, whereas in \cite{FriLo17} it is only applied to the real case. The generalization to the complex case is immediate considering the complex exponential and the principal branch of the complex logarithm.

	\emph{2.} We have proven that the hypothesis
	\begin{equation}
		\sum_{s\in [0,T)\cap D_g} \left|\ln\left|1+\beta(s)\Delta^+g(s)\right|\right|<\infty\end{equation}
	occurring in
	\cite[Definition 6.1 and Lemma 6.5]{FriLo17} is not necessary, it being a direct consequence of $\beta \in
	\mathcal{L}_g^1([0,T);{\mathbb F})$ and $g(T)<\infty$. This was also proven in \cite[Lemma 3.1]{Ma21} for the real case.

	\emph{3.} The solution obtained generalizes that in  \cite[Lemma 6.5]{FriLo17}. Indeed, in the particular case $\beta \in \mathcal{L}_g^1([0,T);\mathbb{R})$ and given that $1+\beta(t)\Delta^+g(t)\neq 0$ for every $t\in [0,T)\cap D_g$, we have that, for every $t\in [0,T)\cap D_g$,
	\begin{equation}
		\operatorname{Arg}(1+\beta(t)\Delta^+g(t)) = \left\{
		\begin{array}{rl}
			\pi, & 1+\beta(t)\Delta^+g(t)<0,  \\
			0, & 1+\beta(t)\Delta^+g(t)>0.
		\end{array}\right.\end{equation}
	Hence, if we write $T_{\beta}^-:=\{t\in [0,T)\cap D_g:\;
	1+\beta(t)\Delta^+g(t)<0\}$ and
	$T_{\beta}^+=\{t\in [0,T)\cap D_g:\;
	1+\beta(t)\Delta^+g(t)>0\}$ (observe that $T_{\beta}^-$ is of finite cardinality),
	we have that
	\begin{equation}
		\ln\left( 1+\beta(t)\Delta^+g(t) \right)=
		\begin{dcases}
			\displaystyle \ln\left|1+\beta(t)\Delta^+g(t) \right|, & t\in T_{\beta}^+, \\
			\displaystyle \ln\left|1+\beta(t)\Delta^+g(t)\right| + i\pi, & t\in T_{\beta}^-.
	\end{dcases}\end{equation}
	Taking into account the previous observations,
	\begin{equation}
		\begin{array}{rcl}
			v(t)&=&\displaystyle
			\exp \left(
			\int_{[0,t)\backslash D_g} {\beta}(s)\, \operatorname{d} \mu_g +
			\sum_{t\in [0,t)\cap D_g} \ln\left|1+\beta(t)\Delta^+g(t) \right| + i
			\sum_{s \in [0,t)\cap T_{\beta}^-} \pi
			\right)  \\
			&=& \displaystyle
			\cos\left(\sum_{s \in [0,t)\cap T_{\beta}^-} \pi  \right)\exp \left(
			\int_{[0,t)\backslash D_g} {\beta}(s)\, \operatorname{d} \mu_g +
			\sum_{t\in [0,t)\cap D_g} \ln\left|1+\beta(t)\Delta^+g(t) \right|
			\right).
	\end{array}\end{equation}
	Hence, if $T_{\beta}^-=\{t_1,\ldots,t_k\}$ and $t_{k+1}:=T$, we get
	\begin{equation}
		v(t)=\begin{dcases}
			\displaystyle
			\exp \left(
			\int_{[0,t)\backslash D_g} {\beta}(s)\, \operatorname{d} \mu_g +
			\sum_{t\in [0,t)\cap D_g} \ln\left|1+\beta(t)\Delta^+g(t) \right|
			\right), & t\in [0,t_1],  \\
			\displaystyle
			\cos(j\, \pi) \exp \left(
			\int_{[0,t)\backslash D_g} {\beta}(s)\, \operatorname{d} \mu_g +
			\sum_{t\in [0,t)\cap D_g} \ln\left|1+\beta(t)\Delta^+g(t) \right|
			\right), & \substack{\displaystyle t\in (t_j,t_{j+1}],\\\displaystyle j=1,\ldots,k,}
	\end{dcases}\end{equation}
	which is precisely the solution in \cite[Lemma 6.5]{FriLo17}.

	\emph{4.} In the case there exists some element $t\in [0,T)\cap D_g$ such that $1+\beta(t)\Delta^+g(t)=0$, the set
\[T_{\beta}^0:=\{t\in [0,T)\cap D_g:\; 1+\beta(t)\Delta^+g(t)=0\}\] is of finite cardinality and, therefore, if we denote by $t^0_\beta:=\min T_{\beta}^0$ if $T_{\beta}^0\ne\emptyset$, $t^0_\beta:=T$ otherwise,
	we have that
	\begin{equation}
		v(t)=\left\{\begin{array}{ll}
			\displaystyle u(g^C(t)) \prod_{s\in [0,t)\cap D_g}
			\left(1+\beta(s)\Delta^+g(s) \right) , & t \in [0,t^0_\beta],  \\
			0, & t\in (t^0_\beta,T].
		\end{array}\right.\end{equation}
	Taking into account that we are assuming that $g$ is continuous at
	$t=0$, we have that $t^0_\beta=\min T_{\beta}^0>0$. Thus,
	$v(t)\neq 0$ for every $t\in [0,t^0_\beta]$.
\end{rem}

\begin{dfn} Given an element $\beta \in \mathcal{L}^1_g([0,T);{\mathbb F})$ and $v_0=1$,
	we denote the solution of problem~\eqref{eq:expg} constructed in  Theorem~\ref{solel} by $\exp_g(\beta;0,t)\in \mathcal{AC}_g([0,T];{\mathbb F})$
	and call it the \emph{complex $g$-exponential} map 
	or just \emph{$g$-exponential} map.
\end{dfn}
\medskip

In the following result we present some important properties of the complex $g$-exponential function.

\begin{pro} \label{gexpprog}Let be $\beta\in \mathcal{L}^1_g([0,T);{\mathbb F})$. The following properties hold:
	\begin{enumerate}
		\item If $a=\operatorname{Re}\beta$ and $b=\operatorname{Im}\beta$ then
		\begin{equation} \label{eq:gcomplex}
			\begin{aligned}
				\exp_g(\beta;0,t)= &
				\prod_{u \in [0,t)\cap D_g} \left( 1+a(u)\Delta^+g(u)+ib(u)\Delta^+g(u) \right)
				\exp\left(\int_{[0,g^C(t))} (a \circ \gamma)\, \operatorname{d}\mu\right)  \\
				&  \cdot
				\left[
				\cos\left(  \int_{[0,g^C(t))}
				(b\circ \gamma) \, \operatorname{d}\mu
				\right)+ i \sin\left( \int_{[0,g^C(t))}
				(b\circ \gamma) \, \operatorname{d}\mu\right)
				\right].
		\end{aligned}\end{equation}
		\item $ \overline{\exp_g(\beta;0,t)}=\exp_g(\overline{\beta};0,t)$, for every $t\in [0,T]$.
		\item Given $n \in \mathbb{N}$, $n\ge 2$, $ \exp_g(\beta;0,t)^n=
		\exp_g(p_n(\beta);0,t)\in \mathcal{AC}_g([0,T];{\mathbb F})$, where
		\begin{equation}
			p_n(\beta)(t)=n\, \beta(t) + \sum_{k=2}^n \binom{n}{k}\, \beta(t)^k \Delta^+ g(t)^{k-1},\quad n \in \mathbb{N},\ n\ge 2.\end{equation}
		\item Given $n \in \mathbb{N}$, $n\ge 2$,
%		and taking $S=\min\{T,t^0_{\beta}\}$, then 
		$\exp_g(\beta;0,t)^{-n} =
		\exp_g(q_n(\beta);0,t)\in \mathcal{AC}_g([0,t^0_{\beta}];{\mathbb F})$, where
		\begin{equation} \label{eq:poweminusn}
			q_n(\beta)(t)=-\frac{p_n(\beta)(t)}{1+p_n(\beta)(t)\,\Delta^+ g(t)},\quad n \in \mathbb{N},\ n\ge 2.\end{equation}
		Observe that $\exp_g(\beta;0,t)^{-n}$ is not well defined 
		in $(t^0_{\beta},T]$ since $\exp_g(\beta;0,\cdot)=0$ in that set.
		\item Given $\beta_1,\beta_2\in \mathcal{L}_g^1([0,T),{\mathbb F})$ such that $\beta_1\beta_2\Delta^+g\in\mathcal{L}_g^1([0,T),{\mathbb F})$, \begin{equation}\exp_g(\beta_1;0,t)\exp_g(\beta_2;0,t)=\exp_g(\beta_1+\beta_2+\beta_1\beta_2\Delta^+g;0,t).\end{equation}
	\end{enumerate}
\end{pro}

\begin{proof}
	1.  Indeed,
	\begin{equation}
		\begin{array}{rcl}
			&&\exp_g(a+b i; 0,t) \\&=&
			\displaystyle \exp\left( \int_{[0,t)\backslash D_g} a(s)\, \operatorname{d}\mu_g + i \int_{[0,t)\backslash D_g}
			b(s) \, \operatorname{d}\mu_g\right) \\ && \displaystyle
			\cdot \exp\left( \sum_{u \in [0,t)\cap D_g} \ln(1+(a(u)+ib(u))\Delta^+g(u))
			\right)  \\
			&=& \displaystyle
			\exp\left( \int_{[0,t)\backslash D_g} a(s)\, \operatorname{d}\mu_g \right)
			\prod_{u \in [0,t)\cap D_g} \left( 1+a(u)\Delta^+g(u)+ib(u)\Delta^+g(u) \right)
			\\ && \displaystyle
			\cdot  \left[
			\cos\left(  \int_{[0,t)\backslash D_g}
			b(s) \, \operatorname{d}\mu_g
			\right)+ i \sin\left( \int_{[0,t)\backslash D_g}
			b(s) \, \operatorname{d}\mu_g\right)
			\right].
	\end{array}\end{equation}
	Now the formula is obtained reasoning as in
	Corollary~\ref{corocambiovar}.

	2. This property is clear from the definition of the complex conjugate.

	3. Observe that $p_n(\beta)\in \mathcal{L}_g^1([0,T);{\mathbb F})$ since
	\begin{equation}%
		\|p_n(\beta) \|_{\mathcal{L}_g^1([0,T);{\mathbb F})} \leq n\, \|\beta\|_{\mathcal{L}^1_g([0,T)}+
		\sum_{k=2}^n \binom{n}{k}\, \sum_{t\in [0,T)\cap D_g}\left(|\beta(t)| \Delta^+ g(t)\right)^k<\infty.\end{equation}

	Thus the solution of problem~\eqref{eq:expg} where we consider $ p_n(\beta)$ instead of $\beta$ is given by $v=v^Bv^C$ where
	\begin{align*}v^B(t)=& \prod_{s\in [0,t)\cap D_g}
		\left(1+\left(n\, \beta(s) + \sum_{k=2}^n \binom{n}{k}\, \beta(s)^k \Delta^+ g(s)^{k-1}\right)\Delta^+g(s) \right),\\
		=& \prod_{s\in [0,t)\cap D_g}
		\left(1+n\, \beta(s)\Delta^+g(s) + \sum_{k=2}^n \binom{n}{k}\, \beta(s)^k \Delta^+ g(s)^{k} \right),\\
		=& \prod_{s\in [0,t)\cap D_g}
		\left(\sum_{k=0}^n \binom{n}{k}\, \beta(s)^k \Delta^+ g(s)^{k} \right)=\prod_{s\in [0,t)\cap D_g}
		\left(1+ \beta(s)\Delta^+g(s)\right)^n\\ = & \left(\prod_{s\in [0,t)\cap D_g}
		\left(1+ \beta(s)\Delta^+g(s)\right)\right)^n,\\
		v^C(t)=&
		\exp\left(\int_{[0,t)} \left(n\, \beta(s) + \sum_{k=2}^n \binom{n}{k}\, \beta(s)^k \Delta^+ g(s)^{k-1}\right) \operatorname{d} \mu_{g^C}\right)=\exp\left(n\int_{[0,t)} \beta(s)\operatorname{d} \mu_{g^C}\right)\\ = & \left[\exp\left(\int_{[0,t)} \beta(s)\operatorname{d} \mu_{g^C}\right)\right]^n.
	\end{align*}
	Hence, $ \exp_g(\beta;0,t)^n=
	\exp_g(p_n(\beta);0,t)$.

	4.  Observe that $q_n(\beta)\in \mathcal{L}_g^1([0,T);{\mathbb F})$ since
	\begin{equation}%
		\|q_n\|_{\mathcal{L}_g^1([0,t^0_{\beta}T);{\mathbb F})}  \leq n\, \|\beta\|_{\mathcal{L}^1_g([0,T)} +
		\sum_{t\in [0,T)\cap D_g} \frac{|p_n(\beta)(t)\, \Delta^+g(t)|}{|1+p_n(\beta)(t)\, \Delta^+ g(t)|}< \infty,\end{equation}
	because $p_n(\beta)\in \mathcal{L}_g^1([0,T);{\mathbb F})$. Therefore, the solution of problem~\eqref{eq:expg}, where we consider $ q_n(\beta)$ instead of $\beta$, is given by $v=v^Bv^C$ where
	\begin{align*}v^B(t)=& \prod_{s\in [0,t)\cap D_g}
		\left(1-\frac{n\, \beta(s) + \sum_{k=2}^n \binom{n}{k}\, \beta(s)^k \Delta^+ g(s)^{k-1}}{1+\left(n\, \beta(s) + \sum_{k=2}^n \binom{n}{k}\, \beta(s)^k \Delta^+ g(s)^{k-1}\right)\Delta^+ g(s)}\Delta^+g(s) \right),\\
		=& \prod_{s\in [0,t)\cap D_g}
		\left(\frac{1}{1+\left(n\, \beta(s) + \sum_{k=2}^n \binom{n}{k}\, \beta(s)^k \Delta^+ g(s)^{k-1}\right)\Delta^+ g(s)}\right) \\ = & \left(\prod_{s\in [0,t)\cap D_g}
		\left(1+ \beta(s)\Delta^+g(s)\right)\right)^{-n},\\
		v^C(t)=&
		\exp\left(\int_{[0,t)}-\frac{n\, \beta(s) + \sum_{k=2}^n \binom{n}{k}\, \beta(s)^k \Delta^+ g(s)^{k-1}}{1+\left(n\, \beta(s) + \sum_{k=2}^n \binom{n}{k}\, \beta(s)^k \Delta^+ g(s)^{k-1}\right)\Delta^+ g(s)} \operatorname{d} \mu_{g^C}\right)\\ = & \exp\left(\int_{[0,t)}-n\, \beta(s)\operatorname{d} \mu_{g^C}\right)=\left( \exp\left(\int_{[0,t)} \beta(s)\operatorname{d} \mu_{g^C}\right)\right) ^{-n}.
	\end{align*}
	Hence, $ \exp_g(\beta;0,t)^{-n}=
	\exp_g(q_n(\beta);0,t)$.

	5. Let $v(t)=\exp_g(\beta_1;0,t)\exp_g(\beta_2;0,t)$. It is easy to check that $v(0)=1$ and
	$$v_g'=\left( \beta_1+\beta_2+\beta_1\beta_2\Delta^+g\right) v,$$ so $v$ solves
	\begin{equation}
		\left\{
		\begin{array}{l}
			v_g'(t)=\left( \beta_1(t)+\beta_2(t)+
			\beta_1(t)\beta_2(t)\Delta^+g(t)\right) v(t),
			\; g-a.e.\; t\in[0,T),  \\
			v(0)=1,
		\end{array}
		\right.\end{equation}
	Thus, by Theorem~\ref{solel}, in order to see that $v(t)=\exp_g(\beta_1+\beta_2+\beta_1\beta_2\Delta^+g;0,t)$ it  is necessary and sufficient to check that $v$ is absolutely continuous, which is true thanks to \cite[Proposition~3.29]{MarquezTesis}.
\end{proof}

\subsection{{$g$-sine} and {$g$-cosine}}
Let us see now how to use the complex $g$-exponential map in order to define the \emph{$g$-sine} and \emph{$g$-cosine} functions. We observe that the presence of jumps in the derivator prevents us from expressing the the exponential as the product of its real and imaginary parts. Indeed, thanks to~\eqref{gexpprog}:
\begin{equation}\label{eq:expgab}
	\begin{aligned}
		\exp_g(a+bi;0,t)=& \exp_g(a;0,t)\, \exp_g\left( \frac{i\, b}{1+a\, \Delta^+g};0,t\right)  \\
		\neq & \exp_g(a;0,t)\, \exp_g(i\,b; 0,t).
\end{aligned}\end{equation}
This fact will have its repercussion when we analyze the case of second order linear equations.
In view of expression \eqref{eq:gcomplex}, it might be interesting to analyze the case
$a=0$, in order to define the $g$-sine and $g$-cosine.

\begin{dfn}[$g$-sine and $g$-cosine] Let $b\in \mathcal{L}^1_g([0,T];\mathbb{F})$. We define $\sin_g(b;0,t)$ and $\cos_g(b;0,t)$, as the first and second components, respectively, of the unique solution in $\mathcal{AC}_g([0,T];\mathbb{F}^2)$
	of the following linear system:
	\begin{equation} \label{eq:gsincos}
		\left\{\begin{array}{l}
			\left(\begin{array}{c}
				\sin_g(b;0,t) \\
				\cos_g(b;0,t)
			\end{array} \right)'_g(t)=\left( \begin{array}{cc}
				0 & b(t) \\
				-b(t) & 0
			\end{array}\right) \left(\begin{array}{c}
				\sin_g(b;0,t) \\
				\cos_g(b;0,t)
			\end{array} \right), \; g-a.e. \, t \in [0,T),  \\
			\sin_g(b;0,0)=0, \; \cos_g(b;0,0)=1.
		\end{array}
		\right.\end{equation}
\end{dfn}

\begin{rem} Observe that \eqref{eq:gsincos} has, indeed, a unique solution in $\mathcal{AC}_g([0,T];\mathbb{F}^2)$ as it satisfies the conditions of \cite[Theorem 7.3]{FriLo17} with $L=|b|$.
\end{rem}

\begin{pro} Given $b\in \mathcal{L}^1_g([0,T];\mathbb{R})$, we have that
	\begin{equation}\label{sincos}
		\begin{array}{rcl}
			\displaystyle \sin_g(b;0,t)&=&\displaystyle  \frac{\exp_g(b i;0,t)-\exp_g(-bi;0,t)}{2i},  \\
			\displaystyle \cos_g(b;0,t)&=&\displaystyle  \frac{\exp_g(b i;0,t)+\exp_g(-bi;0,t)}{2}.
	\end{array}\end{equation}
	Furthermore, developing the previous expressions,
	\begin{equation}\begin{aligned}
			\displaystyle & \sin_g(b;0,t)  \\
			\displaystyle = & \prod_{u \in [0,t)\cap D_g} \left|1+b(u)\, i \, \Delta^+ g(u) \right|
			\sin \left(\sum_{u \in [0,t)\cap D_g} \operatorname{atan}(b(u)\, \Delta^+ g(u)) \right)
			\cos \left( \int_{[0,g^C(t))} (b\circ \gamma) \, \operatorname{d} \mu \right)
			\\ & \displaystyle +
			\prod_{u \in [0,t)\cap D_g} \left|1+b(u)\, i \, \Delta^+ g(u) \right|
			\cos \left(\sum_{u \in [0,t)\cap D_g} \operatorname{atan}(b(u)\, \Delta^+ g(u)) \right)
			\sin \left( \int_{[0,g^C(t))} (b\circ \gamma) \, \operatorname{d} \mu \right),
	\end{aligned}\end{equation}
	\begin{equation}
		\begin{aligned}
			\displaystyle & \cos_g(b;0,t) \\
			\displaystyle = & \prod_{u \in [0,t)\cap D_g} \left|1+b(u)\, i \, \Delta^+ g(u) \right|
			\cos \left(\sum_{u \in [0,t)\cap D_g} \operatorname{atan}(b(u)\, \Delta^+ g(u)) \right)
			\cos \left( \int_{[0,g^C(t))} (b\circ \gamma) \, \operatorname{d} \mu \right)
			\\ \displaystyle  & -
			\prod_{u \in [0,t)\cap D_g} \left|1+b(u)\, i \, \Delta^+ g(u) \right|
			\sin \left(\sum_{u \in [0,t)\cap D_g} \operatorname{atan}(b(u)\, \Delta^+ g(u)) \right)
			\sin \left( \int_{[0,g^C(t))} (b\circ \gamma) \, \operatorname{d} \mu \right).
	\end{aligned}\end{equation}
\end{pro}

\begin{proof}
	Indeed, differentiating the equations in \eqref{sincos},
	\begin{equation}
		\begin{array}{rcl}
			\displaystyle \left( \sin_g(b;0,t)\right)'_g(t)&=&
			\displaystyle \frac{1}{2i} \left(ib(t) \exp_g(bi;0,t)+ib(t)\exp_g(-bi;0,t) \right)  \\
			&=& \displaystyle b(t) \frac{1}{2} \left(\exp_g(bi;0,t)+\exp_g(-bi;0,t)\right)  \\
			&=& b(t) \cos_g(b;0,t),\; g-a.e.\, t \in [0,T).
	\end{array}\end{equation}
	\begin{equation}
		\begin{array}{rcl}
			\displaystyle \left( \cos(b;0,t)\right)'_g(t)&=&
			\displaystyle \frac{1}{2} \left( i b(t) \exp_g(bi;0,t) - ib(t) \exp_g(-bi;0,t)\right) \\
			&=& \displaystyle -b(t) \frac{1}{2i} \left(\exp_g(bi;0,t)-\exp_g(-bi;0,t)\right) \\
			&=& -b(t) \sin_g(b;0,t),\; g-a.e.\, t \in [0,T).
	\end{array}\end{equation}
	Observe also that $\overline{\exp_g(bi;0,t)}=\exp_g(-bi;0,t)$. Hence,
	\begin{equation}
		\begin{array}{rcl}
			\displaystyle
			\cos_g(b;0,t)&=&
			\displaystyle
			\operatorname{Re}(\exp_g(bi;0,t)), \\
			\displaystyle
			\sin_g(b;0,t)&=&
			\displaystyle \operatorname{Im}(\exp_g(bi;0,t)).
	\end{array}\end{equation}
	In particular, we can obtain the  explicit expression of the $g$-sine and the $g$-cosine separating the real and imaginary parts of $\exp_g(bi;0,t)$. We have that
	\begin{equation}
		\begin{array}{rcr}
			\displaystyle \sin_g(b;0,t)&=&\displaystyle \cos \left( \int_{[0,g^C(t))} (b\circ \gamma) \, \operatorname{d}\mu \right)
			\, \operatorname{Im} \left( \prod_{u \in [0,t)\cap D_g} \left(1+b(u)\, i \, \Delta^+ g(u) \right) \right)  \\
			&& \displaystyle +  \sin \left( \int_{[0,g^C(t))} (b\circ \gamma) \, \operatorname{d}\mu \right)
			\, \operatorname{Re} \left( \prod_{u \in [0,t)\cap D_g} \left(1+b(u)\, i \, \Delta^+ g(u) \right) \right),
	\end{array}\end{equation}
	\begin{equation}
		\begin{array}{rcr}
			\displaystyle \cos_g(b;0,t)&=&\displaystyle \cos \left( \int_{[0,g^C(t))} (b\circ \gamma) \, \operatorname{d}\mu \right)
			\, \operatorname{Re} \left( \prod_{u \in [0,t)\cap D_g} \left(1+b(u)\, i \, \Delta^+ g(u) \right) \right)  \\
			&& \displaystyle -  \sin \left( \int_{[0,g^C(t))} (b\circ \gamma) \, \operatorname{d}\mu \right)
			\, \operatorname{Im} \left( \prod_{u \in [0,t)\cap D_g} \left(1+b(u)\, i \, \Delta^+ g(u) \right) \right).
	\end{array}\end{equation}
	Hence, in order to obtain the result it is enough to observe that
	\begin{equation}
		\begin{aligned}
			& \prod_{u \in [0,t)\cap D_g} \left(1+b(u)\, i \, \Delta^+ g(u) \right) =
			\exp\left( \sum_{u \in [0,t)\cap D_g} \log\left(1+b(u)\, i \, \Delta^+ g(u) \right) \right)
			\\
			=&  \exp\left(  \sum_{u\in [0,t)\cap D_g}
			\log  \left|1+b(u)\, i \, \Delta^+ g(u) \right| + i \sum_{u\in [0,t)\cap D_g}
			\operatorname{Arg}\left( 1+b(u)\, i \, \Delta^+ g(u)\right) \right)
			\\
			=&
			\exp\left( \sum_{u\in [0,t)\cap D_g}
			\log  \left|1+b(u)\, i \, \Delta^+ g(u) \right| \right) \,
			\exp\left(  i \sum_{t\in [0,t)\cap D_g} \operatorname{atan}\left( b(u)\, \Delta^+g(u)\right) \right)
			\\
			=&
			\prod_{u \in [0,t)\cap D_g} \left|1+b(u)\, i \, \Delta^+ g(u) \right|   \\
			&  \cdot \left[
			\cos \left(\sum_{u \in [0,t)\cap D_g} \operatorname{atan}(b(u)\, \Delta^+ g(u)) \right)+ i \,
			\sin \left(\sum_{u \in [0,t)\cap D_g} \operatorname{atan}(b(u)\, \Delta^+ g(u)) \right]
			\right).
	\end{aligned}\end{equation}
\end{proof}
\subsection{The non homogeneous case}

In this section we will study the linear non homogeneous problem:
\begin{equation} \label{eq:nohomoexp}
	\left\{ \begin{array}{l}
		v_g'(t)=\beta(t)\,v(t)+f(t),\; g-a.e.\, t \in [0,T),  \\
		v(0)=v_0,
	\end{array}\right.\end{equation}
where $\beta,\, f \in \mathcal{L}^1_g([0,T);{\mathbb F})$ and $v_0\in{\mathbb F}$. We have the
following result whose proof can be achieved using the techniques employed
in \cite[Theorems 3.5 and 4.6]{Ma21}.

\begin{pro} \label{nonhomoexpprop} Let 
%	$S=\min\{T,t_{\beta}^0\}$ and 
	$\beta,\, f \in \mathcal{L}^1_g([0,T);{\mathbb F})$. Then
	the map $v:[0,t^0_{\beta}]\to {\mathbb F}$ defined as:
	\begin{equation} \label{eq:solnohomoexp}
		v(t)=v_0\exp_g(\beta;0,t) + \exp_g(\beta;0,t)\, \int_{[0,t)} \exp_g(\beta;0,s)^{-1}\, \frac{f(s)}{1+\beta(s) \Delta^+g(s)}\, \operatorname{d} \mu_g,\end{equation}
	is the unique solution in $\mathcal{AC}_g([0,t^0_{\beta}];{\mathbb F})$ of problem~\eqref{eq:nohomoexp} in the interval $[0,t^0_{\beta})$. 
%	for {\color{red}$t\in(t^0_{\beta},T]$}, we have that 
%	$\exp_g(\beta;0,t)=0$, so $v(t)=0$, and we cannot ensure that $v$ is a valid solution of~\eqref{eq:solnohomoexp} {\color{red} on that set.}
%	for $t>\min\{T,t_{\beta}^0\}$.
\end{pro}

\begin{rem} Considering \eqref{eq:nohomoexp}, observe that the set $A=\{s\in [0,T)\cap D_g:\; |1+\beta(s)\Delta^+g(s)|<1/2\}$ has at most finite cardinality
	since $\beta \in \mathcal{L}^1_g([0,T];{\mathbb F})$. 
	{ Indeed, given $s\in [0,T)\cap D_g$ such that 
	$1+\beta(s)\Delta^+g(s)\neq 0$,
\begin{displaymath}
\begin{aligned}
0< |1+\beta(s)\Delta^+g(s)|<\frac{1}{2} & \Leftrightarrow  
 0< 1+ 2\operatorname{Re}(\beta(s)) \Delta^+g(s) +
 |\beta(s)\Delta^+g(s)|^2<\frac{1}{4} \\
  & \Leftrightarrow  \frac{3}{4} < 2 
  |\operatorname{Re}(\beta(s))| \Delta^+g(s)-|\beta(s)\Delta^+g(s)|^2<1,
\end{aligned}
\end{displaymath}
so $|\operatorname{Re}(\beta(s))| \Delta^+g(s)>3/8$ and then 
$A$ has finite cardinality.} Therefore:
	\begin{equation}
		\sum_{s\in[0,T)\cap D_g} \frac{|f(s)|\Delta^+g(s)}{|1+\beta(s)\Delta^+g(s)|} \leq
		\sum_{s\in A} \frac{|f(s)|\Delta^+g(s)}{|1+\beta(s)\Delta^+g(s)|} + 2 \sum_{s\in ([0,T)\cap D_g)\backslash A}
		|f(s)|\Delta^+g(s)<\infty.\end{equation}
\end{rem}

\begin{rem}
	Observe that, by Proposition~\ref{nonhomoexpprop}, if we define
	\begin{equation} G(t,s)=\exp_g(\beta;0,t) \frac{\exp_g(\beta;0,s)^{-1}}{1+\beta(s) \Delta^+g(s)}\chi_{[0,t)}(s),\quad t,s\in[0,T),\end{equation}
	$G$ is the \emph{Green's function} associated to problem~\eqref{eq:nohomoexp}, that is, its solution can be expressed as
	\begin{equation} \label{eq:solnohomoexp2}
		v(t)=v_0\exp_g(\beta;0,t) + \int_{[0,T]} G(t,s)f(s)\operatorname{d} \mu_g(s).\end{equation}
\end{rem}

Now, following the same idea as in the previous section, we will see that, under certain hypotheses about the set
$D_g$, it is possible to decompose the solution of \eqref{eq:nohomoexp}  in terms of the solution of
two problems associated with the continuous and discrete part of the derivator.

\begin{cor} Assume $\mu_g(\overline{D_g}\backslash D_g)=0$  and let  
$\beta,\, f \in \mathcal{L}^1_g([0,T);{\mathbb F})$. Then the unique solution 
$v\in \mathcal{AC}_g([0,t_{\beta}^0];{\mathbb F})$~\eqref{eq:solnohomoexp} of
	\eqref{eq:nohomoexp} in the interval $[0,t_{\beta}^0)$ can be expressed in the following terms:
	\begin{equation}
		v(t)=v^C(t)\,\widetilde{v}^B(t)+v^B(t)\,\widetilde{v}^C(t),\; \forall t \in [0,t_{\beta}^0],\end{equation}
	where
	\begin{itemize}
		\item $v^C\in \mathcal{AC}_{g^C}([0,t_{\beta}^0];{\mathbb F})$ is the unique solution of \eqref{eq:expgC} in the interval $[0,t_{\beta}^0)$ given by \eqref{eq:solexpgC},
		\item $v^B\in \mathcal{AC}_{g^B}([0,t_{\beta}^0];{\mathbb F})$ is the unique solution of \eqref{eq:expgB} in the interval $[0,t_{\beta}^0)$ given by \eqref{eq:solexpgB},
		\item $\widetilde{v}^B\in \mathcal{AC}_{g^B}([0,t_{\beta}^0];{\mathbb F})$ is the unique solution of:
		\begin{equation} \label{eq:soltildevb}
			\left\{\begin{array}{l}
				\displaystyle (\widetilde{v}^B)'_{g^B}(t)=\beta(t)\, \widetilde{v}^B(t) + \frac{f(t)}{v^C(t)(1+\beta(t) \Delta^+g(t))},\; g^B-a.e.\, t \in [0,t_{\beta}^0),  \\
				\displaystyle
				\widetilde{v}^B(0)=\frac{1}{2},
			\end{array}\right.\end{equation}
		given by:
		\begin{equation}
			\widetilde{v}^B(t)=v^B(t) \bigg[\frac{1}{2}+
			\sum_{s\in [0,t)\cap D_g} v^B(s)^{-1} \frac{f(s)\Delta^+g(s)}{v^C(s) (1+\beta(s)\Delta^+g(s))}
			\bigg],\end{equation}
		\item $\widetilde{v}^C \in \mathcal{AC}_{g^C}([0,t_{\beta}^0];{\mathbb F})$ is the unique solution of:
		\begin{equation} \label{eq:soltildevc}
			\left\{\begin{array}{l}
				\displaystyle (\widetilde{v}^C)'_{g^C}(t)=\beta(t) \,\widetilde{v}^C(t) +\frac{f(t)}{v^B(t)},\; g^C-a.e.\, t \in [0,t_{\beta}^0), \\
				\displaystyle \widetilde{v}^C(0)=\frac{v_0}{2},
			\end{array}\right.\end{equation}
		given by:
		\begin{equation}
			\widetilde{v}^C(t)=v^C(t)\, \bigg[\frac{v_0}{2}+\int_{[0,t)} v^C(s)^{-1}\frac{f(s)}{v^B(s)}\, \operatorname{d} \mu_{g^C}
			\bigg].\end{equation}
	\end{itemize}
\end{cor}

%{ La unicidad del problema lineal homog\'eneo est\'a pendiente de analizar ojo!!!!}

\begin{proof} Let us consider the solution of \eqref{eq:nohomoexp} given by \eqref{eq:solnohomoexp} and
	the decomposition $\exp_g(\beta;0,t)=v^B(t) v^C(t)$ given by Proposition~\ref{solel}. We have that:
	\begin{equation}
		\begin{array}{rcl}
			v(t)&=&\displaystyle v^B(t)\, v^C(t)\, \bigg[ v_0+ \int_{[0,t)} v^C(s)^{-1} \, \frac{f(s)}{v^B(s)} \, \operatorname{d} \mu_{g^C} \\
			&&\displaystyle
			+
			\sum_{s\in [0,t)\cap D_g} v^B(s)^{-1}\, \frac{f(s)\Delta^+g(s)}{v^B(s)(1+\beta(s)\Delta^+g(s))}
			\bigg]  \\
			&=& \displaystyle v^B(t) \, v^C(t)\, \bigg[\frac{v_0}{2} + \int_{[0,t)} v^C(s)^{-1} \, \frac{f(s)}{v^B(s)} \, \operatorname{d} \mu_{g^C}(s)\bigg] \vspace{0,2cm} \\
			&& \displaystyle +v^C(t) \, v^B(t) \, \bigg[\frac{v_0}{2}+ \sum_{s\in [0,t)\cap D_g} v^B(s)^{-1}\, \frac{f(s)\Delta^+g(s)}{v^B(s)(1+\beta(s)\Delta^+g(s))}
			\bigg]  \\
			&=& v^B(t)\, \widetilde{v}^C(t)+v^C(t)\,\widetilde{v}^B(t).
	\end{array}\end{equation}
	Finally, by Proposition~\ref{nonhomoexpprop} we have that $\widetilde{v}^C\in \mathcal{AC}_{g^C}([0,t_{\beta}^0];{\mathbb F})$ is the unique solution
	of \eqref{eq:soltildevc} and $\widetilde{v}^B \in \mathcal{AC}_{g^B}([0,t_{\beta}^0];{\mathbb F})$ is the unique solution of \eqref{eq:soltildevb}.
\end{proof}

\subsection{Additional regularity}

In order to correctly define regular solutions, throughout 
this section we will assume that $g:\mathbb{R}\to \mathbb{R}$ is a derivator 
such that $0\notin N_g^-$ and $T\notin N_g^+\cup D_g 
\cup C_g$. We also assume that $t_{\beta}^0=T$, 
otherwise, we redefine $T$ by taking $\min\{T,t_{\beta}^0\}$.

Let us check now that we can obtain solutions of the problem~
\eqref{eq:expg} with greater regularity in the case $\beta \in
\mathcal{BC}_{g}([0,T];{\mathbb F})$. We need the following result, which we state for scalar equations.

{
\begin{pro}[{\cite[Proposition 7.6]{FriLo17}}]\label{props} Let 
$x\in \mathcal{AC}_g([0,T];\mathbb{R})$ be a solution of 
\begin{equation}%
x_{g}^{\prime}(t)=f(t, x(t)), \quad g\text {-a.a. } t \in [0,T).\end{equation}
If $f(\cdot, x(\cdot))$ is $g$-continuous on $[0,T]$, then		
\begin{equation}
\left(x\right)_{g}^{\prime}(t)=f(t, x(t)) \quad 
\text { for all } t \in [0,T) \backslash C_{g}.\end{equation}
\end{pro}}

We have the following corollary.

\begin{cor} \label{regkk1} Let $\beta \in \mathcal{BC}_{g}([0,T];{\mathbb F})$, then the problem
	\begin{equation} \label{eq:expg2}
		\left\{
		\begin{array}{l}
			v_g'(t)-\beta(t)\,v(t)=0,\; \forall \, t \in [0,T],  \\
			v(0)=v_0,
		\end{array}
		\right.\end{equation}
	admits a unique solution in the space $\mathcal{BC}_g^{1}([0,T];{\mathbb F})$.
\end{cor}
\begin{proof} Indeed, on the one hand, we have that
	$\mathcal{BC}_g([0,T];{\mathbb F})\subset \mathcal{L}_g^1([0,T);{\mathbb F})$,
	so there exists a unique solution $v\in \mathcal{AC}_g([0,T];{\mathbb F})$ given by \eqref{eq:solexpg}.
	Let us see that $v\in \mathcal{BC}_g^{1}([0,T];{\mathbb F})$ and that $v$ satisfies equation
	\eqref{eq:expg2} on all of the interval $[0,T]$.  First observe that
	$\beta\, v \in \mathcal{BC}_g([0,T];{\mathbb F})$, so, thanks to Proposition~\ref{props}
	we have that $v_g'(t)=\beta(t)\,v(t)$ for all $t\in [0,T)\setminus C_g$. Observe that
	we can extend the result to $t=T$ thanks to the fact that $T\notin N_g^+\cup D_g$.
	Finally, thanks to Definition~\ref{Stieltjesderivative2}, we have the desired result since
	$v_g'=\beta\, v \in \mathcal{BC}_g([0,T];{\mathbb F})$.
\end{proof}

\begin{rem} In the case where $\beta \in \mathcal{BC}_g^1([0,T])$ we cannot ensure that
	$v\in \mathcal{BC}_g^2([0,T];{\mathbb F})$ since the product of two $\mathcal{BC}_g^1([0,T];{\mathbb F})$
	functions is not, in general, a $\mathcal{BC}_g^1([0,T];{\mathbb F})$ function. However, if $\beta$ and its $g$-derivatives are also continuous,
	we can recover the desired regularity as a consequence of Corollary \ref{regadi}.
\end{rem}

\begin{cor} \label{redadicional} Let $\beta \in
	\mathcal{BC}_{g}^n([0,T];{\mathbb F})\cap \mathcal{BC}_{g^C}^{n-1}([0,T];{\mathbb F})$,
	with $n\in \mathbb{N}$, then the problem~\eqref{eq:expg2} admits a unique solution in the space $\mathcal{BC}_g^{n+1}([0,T];{\mathbb F})$.
\end{cor}

Let us now analyze the non homogeneous case.

\begin{cor} Let $\beta,\, f \in \mathcal{BC}_g([0,T];{\mathbb F})$, then the problem
	\begin{equation} \label{eq:nohomoexp2}
		\left\{ \begin{array}{l}
			v_g'(t)=\beta(t)\,v(t)+f(t),\; \forall t \in [0,T],  \\
			v(0)=v_0,
		\end{array}\right.\end{equation}
	admits a unique solution in the space $\mathcal{BC}^1_g([0,T];{\mathbb F})$.
\end{cor}

\begin{cor}\label{corsolregnh}  Let $\beta \in \mathcal{BC}_{g}^n([0,T];{\mathbb F})
	\cap \mathcal{BC}_{g^C}^{n-1}([0,T];{\mathbb F})$
	and $f\in \mathcal{BC}_{g}^n([0,T];{\mathbb F})$, with $n\in \mathbb{N}$, then problem~\eqref{eq:nohomoexp2}
	admits a unique solution in the space $\mathcal{BC}_g^{n+1}([0,T];{\mathbb F})$.
\end{cor}
\begin{exa} Consider any derivator $g$ and the equation
	\begin{equation} \label{eq:nohomoexp4}
		\left\{ \begin{array}{l}
			v_g'(t)=x\, v(t)+\exp_g(z;0,t),\; \forall t \in [0,T],  \\
			v(0)=1,
		\end{array}\right.
	\end{equation}
	where $x,z\in{\mathbb F}$ are constants.   Defining $\beta(t):=x$, $f(t):=\exp_g(z;0,t)$ we have that $\beta$, $f \in \mathcal{BC}_{g}^{\infty}([0,T];{\mathbb F})$. By Corollary \ref{redadicional}, problem~\eqref{eq:nohomoexp4} has a unique solution $v\in\mathcal{C}^\infty_g([0,T];{\mathbb F})$, which, by Proposition~\ref{nonhomoexpprop}, is provided by  expression \eqref{eq:solnohomoexp} as
	\begin{equation}%
		v(t)=  \exp_g(x;0,t) + \exp_g(x;0,t)\, \int_{[0,t)} \exp_g(x;0,s)^{-1}\, \frac{\exp_g(z;0,s)}{1+x\,\Delta^+g(s)}\, \operatorname{d} \mu_g(s).\end{equation}
	Now, by Proposition~\ref{gexpprog},
	\begin{equation}%
		\begin{aligned}
			\exp_g(x;0,s)^{-1}\, \exp_g(z;0,s)=& \exp_g\left( -\frac{x}{1+x\, \Delta^+g};0,s\right) \,  \exp_g(z;0,s)\\
			=& \exp_g\left( -\frac{x}{1+x\, \Delta^+g}+z-\frac{x\,z\, \Delta^+g}{1+x\, \Delta^+g};0,s\right) \\
			=& \exp_g\left( \frac{z-x}{1+x\, \Delta^+g};0,s\right) .
	\end{aligned}\end{equation}
	Therefore,
	\begin{equation}%
		\begin{aligned}
			&\int_{[0,t)} \exp_g(x;0,s)^{-1}\, \frac{\exp_g(z;0,s)}{1+x\,\Delta^+g(s)}\, \operatorname{d} \mu_g(s)\\=&
			(z-x)^{-1} \, \bigg[ \exp_g\left( \frac{z-x}{1+x\, \Delta^+g};0,t\right) -
			\exp_g\left( \frac{z-x}{1+x\, \Delta^+g};0,0\right) \bigg] \\
			=&(z-x)^{-1}\, \big[\exp_g(x;0,t)^{-1}\, \exp_g(z;0,t)-1\big].
	\end{aligned}\end{equation}
	Finally,
	\begin{equation}%
		\begin{aligned}
			v_p(t)=&  \exp_g(x;0,t) + \exp_g(x;0,t)\, (z-x)^{-1}\, \big[\exp_g(x;0,t)^{-1}\, \exp_g(z;0,t)-1\big]\\
			=&  \exp_g(x;0,t) + (z-x)^{-1}\,\big[ \exp_g(z;0,t)- \exp_g(x;0,t)\big].
	\end{aligned}\end{equation}
	Observe that, differentiating $v_g'$ again, we obtain that $v_g^{\prime \prime}-
	(x+z)\,v_g^{\prime}+x\,z\,v=0$, so, for any values $P,\,Q\in{\mathbb C}$, taking $x=(-P+\sqrt{P^2-4Q})/2$, $z=P-x$, $v$ solves the equation $v_g^{\prime \prime}+P\,v_g^{\prime}+Q\,v=0$.

	This fact illustrates how we can obtain a solution of a second order problem~from a first order problem. In the next section we will enter into the detail of  second order problems.
\end{exa}

\section{Linear $g$-differential problems of second order with constant coefficients}\label{SSO}

In this section we consider $g$-differential problems of second order
with constant coefficients. Since we will assume that the coefficients are constant, we will look for  solutions in the space $\mathcal{BC}_g^2([0,T];{\mathbb F})$. Once again, we assume that $0\notin N_g^-$ and $T\notin N_g^+\cup D_g 
\cup C_g$.

\subsection{The homogeneous case}

Let us consider the second order homogeneous linear Cauchy problem
\begin{equation} \label{eq:secondorder}
	\left\{\begin{array}{l}
		\displaystyle v_g''(t)+ P\, v_g'(t)+Q\,v(t)=0,\; \forall t \in [0,T],  \\
		\displaystyle v(0)=x_0,  \\
		\displaystyle v_g'(0)=v_0,
	\end{array}\right.\end{equation}
where $P,\, Q,\, x_0,\, v_0 \in {\mathbb F}$.
We start by defining what we understand as a solution of problem~\eqref{eq:secondorder}.

\begin{dfn}We say $v\in \mathcal{BC}_g^2([0,T];{\mathbb F})$ is a solution of
	\eqref{eq:secondorder} if it satisfies the equation
	\begin{equation}
		v_g''(t)+ P\, v_g'(t)+Q\,v(t)=0,\; \forall t \in [0,T]\end{equation}
	and the initial conditions $v(0)=x_0$ and $v_g'(0)=v_0$.
\end{dfn}
We have the following Lemma, whose proof is straightforward from the linearity of the $g$-derivative.

\begin{lem} \label{2orderexis} Given
	$v_1,\,v_2 \in \mathcal{BC}_g^2([0,T];{\mathbb F})$ such that
	\begin{equation}\label{eq:secondorder2}
		(v_k)_g''(t)+P\, (v_k)_g'(t) + Q\, v_k(t)=0,\; \forall t \in [0,T], \; k=1,2,\end{equation}
	there exist constants $c_1,\, c_2 \in {\mathbb F}$ such that $v=c_1\, v_1+c_2\,v_2$ is a solution of \eqref{eq:secondorder} if $v_1(0)\, (v_2)'_g(0)-v_2(0)\,(v_1)'_g(0)\neq 0$.
\end{lem}

\begin{thm}\label{thmsol2o} For \eqref{eq:secondorder}, the following hold:
	\begin{itemize}
		\item If $P^2-4\,Q\neq 0$, then, defining $\lambda_1=
		(-P+ \sqrt{P^2-4\,Q})/2$ and $\lambda_2=(-P- \sqrt{P^2-4\,Q})/2$, we have that
		\begin{equation}%
			v(t)=\left(\frac{v_{0}-\lambda _{2}\,x_{0}}{\lambda _{1}-\lambda _{2}}\right)\, \exp_g(\lambda_1;0,t)
			-\left(\frac{v_{0}-\lambda _{1}\,x_{0}}{\lambda _{1}-\lambda _{2}}\right) \, \exp_g(\lambda_2;0,t)\end{equation}
		is a solution of \eqref{eq:secondorder}. Furthermore, $v\in
		\mathcal{BC}_g^{\infty}([0,T];{\mathbb F})$ and it is the unique solution in that space.
		\item If $P^2-4\,Q=0$, then, taking $\lambda=-P/2$,
		\begin{equation}%
			v(t)=x_0\, \exp_g(\lambda;0,t)+(v_0-\lambda\, x_0)\, \exp_g(\lambda;0,t)\, \int_{[0,t)}
			\frac{1}{1+\lambda\, \Delta^+g(s)}\, \operatorname{d} \mu_g(s)\end{equation}
		is a solution of \eqref{eq:secondorder}. Furthermore, $v\in
		\mathcal{BC}_g^{\infty}([0,T];{\mathbb F})$ and it is the unique solution in that space.
	\end{itemize}
\end{thm}

\begin{proof}We consider the \emph{characteristic equation} of problem
	\eqref{eq:secondorder},
	\begin{equation} \label{eq:characteristic}
		\lambda^2+P \, \lambda + Q=0.
	\end{equation}

	If $P^2-4\,Q\neq 0$, let
	$v_1=\exp_g(\lambda_1;0,t)$ and $v_2(t)=\exp_g(\lambda_2;0,t)$. By Corollary~\ref{redadicional} we have that
	$v_1,\, v_2 \in \mathcal{BC}_g^{\infty}([0,T];{\mathbb F})\subset \mathcal{BC}_g^2([0,T];{\mathbb F})$. Furthermore, it can be checked that both functions satisfy \eqref{eq:secondorder2}. On the other hand,
	\begin{equation}
		v_1(0)\, (v_2)'_g(0)-v_2(0)\,(v_1)'_g(0)=\lambda_2-\lambda_1 \neq 0.\end{equation}
	Hence, by Lemma~\ref{2orderexis}, there exists a solution of problem~
	\eqref{eq:secondorder} given by
	\begin{equation}
		v(t)=\left(\frac{v_{0}-\lambda _{2}\,x_{0}}{\lambda _{1}-\lambda _{2}}\right)\, \exp_g(\lambda_1;0,t)
		-\left(\frac{v_{0}-\lambda _{1}\,x_{0}}{\lambda _{1}-\lambda _{2}}\right) \, \exp_g(\lambda_2;0,t).\end{equation}

	If $P^2-4\,Q=0$ we get the double root $\lambda=-P/2$ of the characteristic equation.
	Observe that the left hand side of the equation occurring in
	\eqref{eq:secondorder} can be written as $(\partial_g+P/2)^2v$ where $\partial_g$ denotes the $g$-derivative operator. Hence, we define $v_1(t):=\exp_g(\lambda;0,t))$, which is a solution of
	$(\partial_g+P/2)v=0$ and consider the unique solution of problem
	\begin{equation} \label{eq:nohomoexp3}
		\left\{ \begin{array}{l}
			(v_2)_g'(t)=\lambda\, v_2(t)+v_1(t),\; g-a.e.\, t \in [0,T),  \\
			v_2(0)=0,
		\end{array}\right.\end{equation}
	Since $(\partial_g+P/2)v_1=0$, it is clear that $(\partial_g+P/2)^2v_2=0$. Furthermore, $v_2(0)=0$
	and $(v_2)_g'(0)=1$, so $v_2$ is the solution we are looking for. By Corollary~\ref{redadicional}, $v_1\in \mathcal{BC}_g^{\infty}([0,T];{\mathbb F})$ and, applying Corollary~\ref{corsolregnh}, $v_2\in\mathcal{BC}_g^{\infty}([0,T];{\mathbb F})$ as well. Now, thanks to Proposition~\ref{nonhomoexpprop}, we have that:
	\begin{equation}%
		v_2(t)=\exp_g(\lambda;0,t) \, \int_{[0,t)} \frac{1}{1+\lambda \, \Delta^+ g(s)}\, \operatorname{d} \mu_g.\end{equation}
	Since $v_2(0)=0$, $(v_2)'_g(0)=1$ we have that:
	\begin{equation}%
		v_1(0)\, (v_2)'_g(0)-v_2(0)\,(v_1)'_g(0)=1 \neq 0,\end{equation}
	Thus, by Lemma~\ref{2orderexis}, there exists a solution of problem~\eqref{eq:secondorder} given by
	\begin{equation}%
		v(t)=x_0\, \exp_g(\lambda;0,t)+(v_0-\lambda\, x_0)\, \exp_g(\lambda;0,t)\, \int_{[0,t)} \exp_g(\lambda \, \chi_{D_g} ;0,s)\, \operatorname{d} \mu_g.\end{equation}

	Finally if we define $u(t)=v_g'(t)\in \mathcal{BC}_g^{\infty}([0,T])$ we have that the pair of functions
	$(u,v)\in [\mathcal{AC}_g([0,T];{\mathbb F} )]^2$ satisfies the following system of differential
	equations:

	\begin{equation}%
		\left\{\begin{aligned}
			\begin{pmatrix}
				v  \\
				u
			\end{pmatrix}'_g(t)= &
			\begin{pmatrix}
				0 & 1 \\
				-Q & -P
			\end{pmatrix} \begin{pmatrix}
				v (t)\\
				u (t)
			\end{pmatrix},  \\
			v(0)= & x_0,\; u(0)=v_0.
		\end{aligned}\right.\end{equation}
	Thanks to \cite[Theorem 7.3]{FriLo17} we have that the previous system has a
	unique solution in $[\mathcal{AC}_g([0,T];{\mathbb F} )]^2$, therefore $v$ is the unique
	solution in the space $\mathcal{BC}_g^{\infty}([0,T];{\mathbb F})$ of
	problem~\eqref{eq:secondorder2}.
\end{proof}

\subsection{The non homogeneous case}

In this section we focus on the non homogeneous version of the second order linear problem, namely,
\begin{equation} \label{eq:secondordernh}
	\left\{\begin{aligned}
		& v_g''(t)+ P\, v_g'(t)+Q\,v(t)=f(t),\; \forall t \in [0,T], \\
		& v(0)=x_0, \\
		& v_g'(0)=v_0,
	\end{aligned}\right.\end{equation}
where $P,Q,x_0,v_0\in{\mathbb F}$ are constant values and $f\in \mathcal{{\mathcal A}{\mathcal C}}_g([0,T];{\mathbb F})$. Since the coefficients
are constant, it will be the regularity of the term $f$ that determines the additional
regularity of the solution. As in the previous sections, we will see that it is possible to prove the uniqueness of solution when we consider the solution
in the space $\mathcal{BC}_g^2([0,T];{\mathbb F})$.

\begin{thm}\label{thmsol2o2} { Let be 
$f\in \mathcal{BC}^n_g([0,T];{\mathbb F})$ and assume $1+\lambda \, \Delta^+g(t)\neq 0$, for all
	$t\in [0,T)\cap D_g$} and $\lambda\in {\mathbb F}$ such that $\lambda^2+P\,\lambda +
	Q =0$. Then problem~\eqref{eq:secondordernh} has a unique solution $v\in\mathcal{BC}_g^{n+2}([0,T],{\mathbb F})$ given by
	\begin{equation} \label{eq:gen1}
		\begin{aligned}
			v(t)= &
			x_0\,\exp_g(\lambda_2;0,t)+
			\left( v_0-\lambda_2x_0\right)\, \exp_g(\lambda_2;0,t)\\ &\cdot \int_{[0,t)} \frac{\exp_g(\lambda_2;0,s)^{-1}}{1+\lambda_2 \Delta^+g(s)}\,
			\exp_g(\lambda_1;0,s)\, \operatorname{d} \mu_g(s)  \\
			&  + \exp_g(\lambda_2;0,t)\, \int_{[0,t)}
			\frac{\exp_g(\lambda_2;0,s)^{-1}}{1+\lambda_2 \Delta^+g(s)}\,
			\exp_g(\lambda_1;0,s)\\ &\cdot\left(
			\int_{[0,s)} \frac{\exp_g(\lambda_1;0,r)^{-1}}{1+\lambda_1 \Delta^+g(r)}\, f(r)\, \operatorname{d} \mu_g(r)\right)\, \operatorname{d} \mu_g(s).
 \end{aligned}
\end{equation}
\end{thm}
\begin{proof}Let be $\lambda_1,\lambda_2$ the two complex eigenvalues of the characteristic polynomial $x^2+Px+Q=0$. Assume $v_2\in{\mathcal C}_g^{n+2}([0,T],{\mathbb F})$ is a solution of problem~\eqref{eq:secondordernh}. Observe that, if we define $v_1=(v_2)_g'-\lambda_2v_2$, it is clear that $v_1\in{\mathcal C}_g^{n+1}([0,T],{\mathbb F})$  and $(v_1)_g'-\lambda_1v_1=f$, $v_1(0)=(v_2)_g'(0)+\lambda_2v_2(0)=v_0-\lambda_2x_0$, so $v_1$ has to solve the problem
	\begin{equation}\label{pa1}\left\{\begin{array}{l}
			(v_1)'_g(t)=  \lambda_1v_1(t)+f(t),\; g-a.e. \, t \in [0,T)\\
			v(0)=  v_0-\lambda_2x_0.
		\end{array}\right.\end{equation}
		By Corollary~\ref{corsolregnh}, problem~\eqref{pa1} has a unique solution in ${\mathcal C}_g^{n+1}([0,T],{\mathbb F})$, so $v_1$ is that unique solution. Furthermore, by definition, $(v_2)_g'+\lambda_2v_2=v_1$ and $v_2(0)=x_0$, so $v_2$ solves the problem
	\begin{equation}\label{pa2}\left\{\begin{array}{l}
			(v_2)'_g(t)=  \lambda_2v_2(t)+v_1(t),\; g-a.e. \, t \in [0,T)\\
			v(0)=  x_0.
		\end{array}\right.\end{equation}
	By	Corollary~\ref{corsolregnh}, problem~\eqref{pa2} has a unique solution in ${\mathcal C}_g^{n+2}([0,T],{\mathbb F})$, so $v_2$ is that unique solution.
	This implies that, if a solution in ${\mathcal C}_g^{n+2}([0,T],{\mathbb F})$ of problem~	\eqref{eq:secondordernh} exists, it has to be unique.

	In order to obtain that unique solution it is enough to retrace the steps we have taken to prove the uniqueness. Let $v_1$ be the unique solution of problem~\eqref{pa1} in ${\mathcal C}_g^{n+1}([0,T],{\mathbb F})$ and let $v_2$ be the unique solution of problem~\eqref{pa2} in ${\mathcal C}_g^{n+2}([0,T],{\mathbb F})$. Clearly $v_2$ is a solution of problem~\eqref{eq:secondordernh}.

	In order to obtain the explicit expression of the solution, observe that, by Proposition~\ref{nonhomoexpprop}, $v_1$ is of the form
	\begin{equation*}
		v_1(t)=\left( v_0-\lambda_2x_0\right) \,\exp_g(\lambda_1;0,t) + \exp_g(\lambda_1;0,t)\, \int_{[0,t)} \exp_g(\lambda_1;0,s)^{-1}\, \frac{f(s)}{1+\lambda_1 \Delta^+g(s)}\, \operatorname{d} \mu_g(s),
	\end{equation*}
	and $v_2$ of the form
	\begin{equation}%
		\begin{aligned}
			&	v_2(t)=  x_0\exp_g(\lambda_2;0,t) + \exp_g(\lambda_2;0,t)\, \int_{[0,t)} \exp_g(\lambda_2;0,s)^{-1}\, \frac{v_1(s)}{1+\lambda_2 \Delta^+g(s)}\, \operatorname{d} \mu_g(s)\\ = &  x_0\exp_g(\lambda_2;0,t) +\left( v_0-\lambda_2x_0\right)  \,
			\exp_g(\lambda_2;0,t)\, \int_{[0,t)} \exp_g(\lambda_2;0,s)^{-1}\, \frac{\exp_g(\lambda_1;0,s)}{1+\lambda_2 \Delta^+g(s)}\, \operatorname{d} \mu_g(s)\\
			& + \exp_g(\lambda_2;0,t)\, \int_{[0,t)} \exp_g(\lambda_2;0,s)^{-1}  \frac{\exp_g(\lambda_1;0,s)}{1+\lambda_2 \Delta^+g(s)}\\
			& \cdot
			\bigg(\int_{[0,s)} \exp_g(\lambda_1;0,r)^{-1}\, \frac{f(r)}{1+\lambda_1 \Delta^+g(r)}\, \operatorname{d} \mu_g(r)\bigg)\, \operatorname{d} \mu_g(s).
	\end{aligned}\end{equation}
	Now, thanks to Proposition~\ref{gexpprog},
	\begin{equation}%
		\begin{aligned}
			\exp_g(\lambda_2;0,t)^{-1}\,\exp_g(\lambda_1;0,t) =&
			\exp_g(-\lambda_2/(1+\lambda_2\,\Delta^+g(t));0,t)\, \exp_g(\lambda_1;0,t)  \\
			=& \exp_g((\lambda_1-\lambda_2)/(1+\lambda_2\,\Delta^+g(t));0,t),
	\end{aligned}\end{equation}
	thus,
	\begin{equation}%
		\begin{aligned}
			v_2=&x_0\exp_g(\lambda_2;0,t) \\&+
			\left( v_0-\lambda_2x_0\right) \,\exp_g(\lambda_2;0,t)\, \int_{[0,t)} \frac{1}{1+\lambda_2 \Delta^+g(s)}\, \exp_g
			\left( \frac{\lambda_1-\lambda_2}{1+\lambda_2\, \Delta^+g};0,s\right) \, \operatorname{d} \mu_g(s)\\
			&+\exp_g(\lambda_2;0,t)\, \int_{[0,t)} \frac{1}{1+\lambda_2 \Delta^+g(s)}\, \exp_g
			\left( \frac{\lambda_1-\lambda_2}{1+\lambda_2\, \Delta^+g};0,s\right)   \\
			& \cdot
			\bigg(\int_{[0,s)} \exp_g(\lambda_1;0,r)^{-1}\, \frac{f(r)}{1+\lambda_1 \Delta^+g(r)}\, \operatorname{d} \mu_g(r)\bigg)\, \operatorname{d} \mu_g(s).
	\end{aligned}\end{equation}
\end{proof}

\begin{rem}
	Note that, for $\lambda\in {\mathbb F}$ such that $\lambda^2+P\,\lambda + Q =0$,  the condition $1+\lambda \, \Delta^+g(t)= 0$ can only happen for a finite number of $t\in[0,T)\cap D_g$.
\end{rem}		

\begin{rem}
	From the previous expression we can derive the Green's function of  problem~\eqref{eq:secondordernh} just by equating
	\begin{equation} \label{eq:green0}
		\begin{aligned}
			& \int_{{\mathbb R}} G(t,r)f(r)\operatorname{d} \mu_g(r)\\
			= & \exp_g(\lambda_2;0,t)\, \int_{[0,t)} \frac{1}{1+\lambda_2 \Delta^+g(s)}\, \exp_g
			\left( \frac{\lambda_1-\lambda_2}{1+\lambda_2\, \Delta^+g};0,s\right)   \\
			& \cdot
			\bigg(\int_{[0,s)} \exp_g(\lambda_1;0,r)^{-1}\, \frac{f(r)}{1+\lambda_1 \Delta^+g(r)}\, \operatorname{d} \mu_g(r)\bigg)\, \operatorname{d} \mu_g(s).
	\end{aligned}\end{equation}
	Now, if we consider the product measure space $([0,T],\mathcal{M}_g\cdot
	\mathcal{M}_g,\mu_g \cdot \mu_g)$, we have by Fubini's Theorem~
	\cite[Theorem~10.10]{Bartle1995}:
	\begin{equation}%
		\begin{aligned}
			& \int_{{\mathbb R}} G(t,r)f(r)\operatorname{d} \mu_g(r)\\
			=& \exp_g(\lambda_2;0,t)\, \int_{[0,T]\cdot [0,T]}
			\frac{1}{1+\lambda_2 \Delta^+g(s)}\, \exp_g
			\left( \frac{\lambda_1-\lambda_2}{1+\lambda_2\, \Delta^+g};0,s\right)  \\
			& \cdot  \exp_g(\lambda_1;0,r)^{-1}\, \frac{f(r)}{1+\lambda_1 \Delta^+g(r)}\,\chi_{[0,s)}(r)\, \chi_{[0,t)}(s)
			\, \operatorname{d} \mu_g\cdot \operatorname{d} \mu_g \\
			= &\exp_g(\lambda_2;0,t)\, \int_{[0,T]} \int_{[0,T]}
			\exp_g(\lambda_1;0,r)^{-1}\,  \exp_g
			\left( \frac{\lambda_1-\lambda_2}{1+\lambda_2\, \Delta^+g};0,s\right)  \\
			& \cdot (1+\lambda_1 \Delta^+g(r))^{-1}\,(1+\lambda_2 \Delta^+g(s))^{-1}\, f(r)\,
			\chi_{(r,t)}(s) \, \chi_{[0,t)}(r)\, \operatorname{d} \mu_g(s)\, \operatorname{d} \mu_g(r) \\
			=&\exp_g(\lambda_2;0,t)\,\int_{[0,t)} \exp_g(\lambda_1;0,r)^{-1}\, \frac{f(r)}{1+\lambda_1\,
				\Delta^+g(r)} \\
			& \cdot \left(\int_{(r,t)} \frac{1}{1+\lambda_2\, \Delta^+g(s)}\,
			\exp_g
			\left( \frac{\lambda_1-\lambda_2}{1+\lambda_2\, \Delta^+g};0,s\right) \, \operatorname{d} \mu_g(s)\right) \, \operatorname{d} \mu_g(r).
	\end{aligned}\end{equation}
	Therefore, for $t,r\in[0,T]$,
	\begin{equation}\label{eq:green}
		\begin{aligned}
			G(t,r)= & \exp_g(\lambda_2;0,t) \, \exp_g(\lambda_1;0,r)^{-1}\,(1+\lambda_1 \Delta^+g(r))^{-1}\, \chi_{[0,t)}(r)\\
			&\cdot  \int_{(r,t)} \frac{1}{1+\lambda_2\, \Delta^+g(s)}\,
			\exp_g
			\left( \frac{\lambda_1-\lambda_2}{1+\lambda_2\, \Delta^+g};0,s\right) \, \operatorname{d} \mu_g(s)\\
			=& \exp_g(\lambda_2;0,t) \, \exp_g(\lambda_1;0,r)^{-1}\,(1+\lambda_1 \Delta^+g(r))^{-1}\, \chi_{[0,t)}(r)\\
			&\cdot  \int_{[r,t)} \frac{1}{1+\lambda_2\, \Delta^+g(s)}\,
			\exp_g
			\left( \frac{\lambda_1-\lambda_2}{1+\lambda_2\, \Delta^+g};0,s\right) \, \operatorname{d} \mu_g(s) \\
			&-\exp_g(\lambda_2;0,t)\, \exp_g(\lambda_2;0,r)^{-1}\,(1+\lambda_1 \Delta^+g(r))^{-1}\,
			(1+\lambda_2 \Delta^+g(r))^{-1}\, \Delta^+g(r)\, \chi_{[0,t)}(r).
	\end{aligned}\end{equation}
	Observe that:
	\begin{itemize}
		\item If $\lambda_1\neq \lambda_2$,
		\begin{equation}%
			v(s)=\exp_g
			\left( \frac{\lambda_1-\lambda_2}{1+\lambda_2\, \Delta^+g};0,s\right)  \in \mathcal{AC}_g([r,T];{\mathbb F})\end{equation}
		is the solution of:
		\begin{equation}%
			\left\{\begin{aligned}
				&v_g^{\prime}(s)=\frac{\lambda_1-\lambda_2}{1+\lambda_2\, \Delta^+g(s)} \, v(s),\; g-a.e. \, s \in [r,T), \\
				&v(r)=\exp_g
				\left( \frac{\lambda_1-\lambda_2}{1+\lambda_2\, \Delta^+g};0,r\right) .
			\end{aligned}\right.\end{equation}
		Therefore:
		\begin{equation}%
			\begin{aligned}
				&\int_{[r,t)}
				\frac{\lambda_1-\lambda_2}{1+\lambda_2\, \Delta^+g(s)}
				\exp_g
				\left( \frac{\lambda_1-\lambda_2}{1+\lambda_2\, \Delta^+g};0,s\right)  \, \operatorname{d} \mu_g(s) \\
				=&\exp_g
				\left( \frac{\lambda_1-\lambda_2}{1+\lambda_2\, \Delta^+g};0,t\right) -\exp_g
				\left( \frac{\lambda_1-\lambda_2}{1+\lambda_2\, \Delta^+g};0,r\right)  \\
				=& \exp_g(\lambda_1;0,t)\, \exp_g(\lambda_2;0,t)^{-1} -
				\exp_g(\lambda_1;0,r)\, \exp_g(\lambda_2;0,r)^{-1}.
		\end{aligned}\end{equation}
		Thus the Green function in the case $\lambda_1\neq \lambda_2$ has the
		following expression:
		\begin{equation} \label{eq:green1}
			\begin{aligned}
				&G(t,r)= \exp_g(\lambda_2;0,t) \, \exp_g(\lambda_1;0,r)^{-1}\,(1+\lambda_1 \Delta^+g(r))^{-1}\, \chi_{[0,t)}(r) \\
				&\cdot (\lambda_1-\lambda_2)^{-1} \left(
				\exp_g(\lambda_1;0,t)\, \exp_g(\lambda_2;0,t)^{-1} -
				\exp_g(\lambda_1;0,r)\, \exp_g(\lambda_2;0,r)^{-1}
				\right) \\
				&-\exp_g(\lambda_2;0,t)\, \exp_g(\lambda_2;0,r)^{-1}\,(1+\lambda_1 \Delta^+g(r))^{-1}\,
				(1+\lambda_2 \Delta^+g(r))^{-1}\, \Delta^+g(r)\, \chi_{[0,t)}(r) \\
				=&+(\lambda_1-\lambda_2)^{-1}\,\exp_g(\lambda_1;0,t)\,\exp_g(\lambda_1;0,r)^{-1}\,(1+\lambda_1 \Delta^+g(r))^{-1}\, \chi_{[0,t)}(r)
				\\
				&-(\lambda_1-\lambda_2)^{-1}\,\exp_g(\lambda_2;0,t)\,\exp_g(\lambda_2;0,r)^{-1}\,(1+\lambda_1 \Delta^+g(r))^{-1}\, \chi_{[0,t)}(r)
				\\
				& -\exp_g(\lambda_2;0,t)\, \exp_g(\lambda_2;0,r)^{-1}\,(1+\lambda_1 \Delta^+g(r))^{-1}\,
				(1+\lambda_2 \Delta^+g(r))^{-1}\, \Delta^+g(r)\, \chi_{[0,t)}(r) \\
				=&+(\lambda_1-\lambda_2)^{-1}\,\exp_g(\lambda_1;0,t)\,\exp_g(\lambda_1;0,r)^{-1}\,(1+\lambda_1 \Delta^+g(r))^{-1}\, \chi_{[0,t)}(r)
				\\
				&-(\lambda_1-\lambda_2)^{-1}\,\exp_g(\lambda_2;0,t)\,\exp_g(\lambda_2;0,r)^{-1}\,(1+\lambda_2 \Delta^+g(r))^{-1}\, \chi_{[0,t)}(r).
		\end{aligned}\end{equation}
		\item If $\lambda_1=\lambda_2$, we have the following expression for
		the Green function~\eqref{eq:green}:
		\begin{equation} \label{eq:green2}
			\begin{aligned}
				& G(r,t) \\ = & \exp_g(\lambda;0,t) \, \exp_g(\lambda;0,r)^{-1}\,(1+\lambda\, \Delta^+g(r))^{-1}\, \chi_{[0,t)}(r)\,
				\\ \cdot & \int_{(r,t)} \frac{1}{1+\lambda \, \Delta^+g(s)}\, \operatorname{d} \mu_g(s) \\
				=&\exp_g(\lambda;0,t) \, \exp_g(\lambda;0,r)^{-1}\,(1+\lambda\, \Delta^+g(r))^{-1}\, \chi_{[0,t)}(r) \\
				& \cdot \left(\int_{[0,t)} \frac{1}{1+\lambda \, \Delta^+g(s)}\, \operatorname{d} \mu_g(s) -
				\int_{[0,r)} \frac{1}{1+\lambda \, \Delta^+g(s)}\, \operatorname{d} \mu_g(s) -
				\frac{\Delta^+g(r)}{1+\lambda\, \Delta^+g(r)}
				\right) \\
				=&+ \exp_g(\lambda;0,t)\, \left(
				\int_{[0,t)} \frac{1}{1+\lambda \, \Delta^+g(s)}\, \operatorname{d} \mu_g(s)\right)\,
				\exp_g(\lambda;0,r)^{-1}\,(1+\lambda\, \Delta^+g(r))^{-1}\, \chi_{[0,t)}(r) \\
				&-\exp_g(\lambda;0,t)\,\left(\int_{[0,r)} \frac{1}{1+\lambda \, \Delta^+g(s)}\, \operatorname{d} \mu_g(s)
				\right) \, \exp_g(\lambda;0,r)^{-1}\,(1+\lambda\, \Delta^+g(r))^{-1}\, \chi_{[0,t)}(r)\\
				&-\exp_g(\lambda;0,t) \, \exp_g(\lambda;0,r)^{-1}\,(1+\lambda\, \Delta^+g(r))^{-2}\,\Delta^+g(r)\, \chi_{[0,t)}(r).
		\end{aligned}\end{equation}
	\end{itemize}
\end{rem}
\begin{rem}
	Observe that we can arrive to expressions \eqref{eq:green1} and
	\eqref{eq:green2} using an integration by parts argument in
	formula~\eqref{eq:green0}. Indeed,
	given two elements $h_1,h_2\in \mathcal{AC}_g([0,T];{\mathbb F})$ we have that
	$h_1\, h_2\in \mathcal{AC}_g([0,T];{\mathbb F})$ and:
	\begin{equation}
		(h_1 \, h_2)^{\prime}_g(t)=(h_1)_g^{\prime}(t) \, h_2(t)+ h_1(t)\, (h_2)^{\prime}_g(t)+
		(h_1)^{\prime}_g(t) \, (h_2)^{\prime}(t) \, \Delta^+ g(t),\; g-a.e.\, t \in [0,T].\end{equation}
	In particular:
	\begin{equation} \label{eq:partes}
		\begin{array}{c}
			\displaystyle
			h_1(t) \, h_2(t)-h_1(0)\, h_2(0)=
			\int_{[0,t)} (h_1)_g^{\prime}(s) \, h_2(s)\, \operatorname{d} \mu_g(s)
			+\int_{[0,t)} h_1(s)\, (h_2)^{\prime}_g(s)\, \operatorname{d} \mu_g(s)  \\
			\displaystyle
			+\int_{[0,t)} (h_1)^{\prime}_g(s) \, (h_2)^{\prime}(s) \, \Delta^+ g(s)\, \operatorname{d}
			\mu_g(s),\; \forall t \in [0,T].
	\end{array}\end{equation}
	Now we study two cases:
	\begin{itemize}
		\item Case $\lambda_1\neq \lambda_2$. Let us consider:
		\begin{equation}%
			\begin{array}{rcl}
				h_1(t)&=& \displaystyle
				(\lambda_1-\lambda_2)^{-1} \exp_g(\lambda_1;0,t)\, \exp_g(\lambda_2;0,t)^{-1},  \\
				h_2(t)&=& \displaystyle
				\int_{[0,t)} \frac{\exp_g(\lambda_1;0,r)^{-1}}{1+\lambda_1 \Delta^+g(r)}\, f(r)\, \operatorname{d} \mu_g(r),
		\end{array}\end{equation}
		we have that $h_1,h_2\in \mathcal{AC}_g([0,T];{\mathbb F})$, so:
		\begin{equation}\label{eq:int2}
			\begin{array}{c}
				\displaystyle  \int_{[0,t)}
				\frac{\exp_g(\lambda_2;0,s)^{-1}}{1+\lambda_2 \Delta^+g(s)}\,
				\exp_g(\lambda_1;0,s)\,\left(
				\int_{[0,s)} \frac{\exp_g(\lambda_1;0,r)^{-1}}{1+\lambda_1 \Delta^+g(r)}\, f(r)\, \operatorname{d} \mu_g(r)\right)\, \operatorname{d} \mu_g(s)
				\\
				\displaystyle =
				(\lambda_1-\lambda_2)^{-1} \exp_g(\lambda_1;0,t)\, \exp_g(\lambda_2;0,t)^{-1}
				\int_{[0,t)} \frac{\exp_g(\lambda_1;0,r)^{-1}}{1+\lambda_1 \Delta^+g(r)}\, f(r)\, \operatorname{d} \mu_g(r)  \\
				\displaystyle -
				\int_{[0,t)} (\lambda_1-\lambda_2)^{-1} \exp_g(\lambda_1;0,s)\, \exp_g(\lambda_2;0,s)^{-1}
				\frac{\exp_g(\lambda_1;0,s)^{-1}}{1+\lambda_1 \Delta^+g(s)}\, f(s)\, \operatorname{d} \mu_g(s)  \\
				\displaystyle -
				\int_{[0,t)}  \frac{\exp_g(\lambda_2;0,s)^{-1}}{1+\lambda_2 \Delta^+g(s)}\,
				\exp_g(\lambda_1;0,s) \,  \frac{\exp_g(\lambda_1;0,s)^{-1}}{1+\lambda_1 \Delta^+g(s)}\, f(s)\,
				\Delta^+ g(s)\, \operatorname{d} \mu_g(s),
		\end{array}\end{equation}
		and we recover the Green's function \eqref{eq:green1}. Observe that by
		substituting \eqref{eq:int2} in \eqref{eq:gen1} we obtain:
		\begin{equation}\label{eq:gen2}
			\begin{aligned}
				\displaystyle v_2(t)=& \left(
				\frac{v_0-\lambda_2x_0}{ \lambda_1-\lambda_2}\right) \, \exp_g(\lambda_1;0,t)-
				\left(\frac{v_0-\lambda_1 x_0}{\lambda_1-\lambda_2}\right)\,
				\exp_g(\lambda_2;0,t)
				\\
				& \displaystyle+
				(\lambda_1-\lambda_2)^{-1} \exp_g(\lambda_1;0,t) \int_{[0,t)} \frac{\exp_g(\lambda_1;0,s)^{-1}}{1+\lambda_1 \, \Delta^+ g(s)}\,
				f(s)\, \operatorname{d} \mu_g
				\\
				&\displaystyle-(\lambda_1-\lambda_2)^{-1} \exp_g(\lambda_2;0,t) \int_{[0,t)}
				\frac{\exp_g(\lambda_2;0,s)^{-1}}{1+\lambda_2 \, \Delta^+g(s)}\, f(s)\, \operatorname{d} \mu_g.
		\end{aligned}\end{equation}
		We have that
		\begin{equation}%
			v_h(t)= \left(
			\frac{v_0-\lambda_2x_0}{ \lambda_1-\lambda_2}\right) \, \exp_g(\lambda_1;0,t)-
			\left(\frac{v_0-\lambda_1 x_0}{\lambda_1-\lambda_2}\right)\,
			\exp_g(\lambda_2;0,t) \in \mathcal{BC}_g^{\infty}([0,T];{\mathbb F})\end{equation}
		is the solution of the homogeneous equation \eqref{eq:secondorder} and
		\begin{equation}\label{eq:particular1}
			\begin{aligned}
				v_p(t)=&+ (\lambda_1-\lambda_2)^{-1} \exp_g(\lambda_1;0,t) \int_{[0,t)} \frac{\exp_g(\lambda_1;0,s)^{-1}}{1+\lambda_1 \, \Delta^+ g(s)}\,
				f(s)\, \operatorname{d} \mu_g
				\\
				&\displaystyle-(\lambda_1-\lambda_2)^{-1} \exp_g(\lambda_2;0,t) \int_{[0,t)}
				\frac{\exp_g(\lambda_2;0,s)^{-1}}{1+\lambda_2 \, \Delta^+g(s)}\, f(s)\, \operatorname{d} \mu_g.
		\end{aligned}\end{equation}
		is a particular solution in the space $\mathcal{BC}_g^{n+2}([0,T],{\mathbb F})$ of the non
		homogeneous equation \eqref{eq:secondordernh} that satisfies $v_p(0)=(v_p)'_g(0)=0$.
		\item Case $\lambda_1=\lambda_2$. We have that
		expression~\eqref{eq:gen1} reduces to:
		\begin{equation} \label{eq:gen3}
			\begin{array}{c}
				\displaystyle v_2=x_0\,\exp_g(\lambda;0,t)+
				\left( v_0-\lambda x_0\right) \, \exp_g(\lambda;0,t)\, \int_{[0,t)} \frac{1}{1+\lambda \Delta^+g(s)}\,
				\operatorname{d} \mu_g(s)  \\
				\displaystyle
				+ \exp_g(\lambda;0,t)\, \int_{[0,t)}
				\frac{1}{1+\lambda \Delta^+g(s)}\,\left(
				\int_{[0,s)} \frac{\exp_g(\lambda;0,r)^{-1}}{1+\lambda \Delta^+g(r)}\, f(r)\, \operatorname{d} \mu_g(r)\right)\, \operatorname{d} \mu_g(s).
		\end{array}\end{equation}
		Denote by:
		\begin{equation}%
			\begin{array}{rcl}
				h_1(t)&=&\displaystyle
				\int_{[0,t)} \frac{1}{1+\lambda \, \Delta^+g(s)}\, \operatorname{d} \mu_g(s),  \\
				h_2(t)&=&\displaystyle
				\int_{[0,t)} \frac{\exp_g(\lambda;0,s)^{-1}}{1+\lambda \Delta^+g(s)}\, f(s)\, \operatorname{d} \mu_g(s).
		\end{array}\end{equation}
		We have that $h_1,h_2\in \mathcal{AC}_g([0,T];{\mathbb F})$. Hence, by formula~\eqref{eq:partes}:
		\begin{equation}\label{eq:int3}
			\begin{aligned}
				\displaystyle
				&\int_{[0,t)}
				\frac{1}{1+\lambda \Delta^+g(s)}\,\left(
				\int_{[0,s)} \frac{\exp_g(\lambda;0,r)^{-1}}{1+\lambda \Delta^+g(r)}\, f(r)\, \operatorname{d} \mu_g(r)\right)\, \operatorname{d} \mu_g(s)
				\\
				\displaystyle = &
				\int_{[0,t)} \frac{1}{1+\lambda \, \Delta^+g(s)}\, \operatorname{d} \mu_g(s) \,
				\int_{[0,t)} \frac{\exp_g(\lambda;0,s)^{-1}}{1+\lambda \Delta^+g(s)}\, f(s)\, \operatorname{d} \mu_g(s)
				\\
				\displaystyle  & - \int_{[0,t)} \left(\int_{[0,s)} \frac{1}{1+\lambda \, \Delta^+g(r)}\, \operatorname{d} \mu_g(r) \right)\,
				\frac{\exp_g(\lambda;0,s)^{-1}}{1+\lambda \Delta^+g(s)}\, f(s)\, \operatorname{d} \mu_g(s)  \\
				\displaystyle & - \int_{[0,t)}
				\frac{\exp_g(\lambda;0,s)^{-1}}{(1+\lambda \Delta^+g(s))^2}\, f(s)\, \Delta^+g(s)\, \operatorname{d} \mu_g(s).
		\end{aligned}\end{equation}
		Substituting expression~\eqref{eq:int3} in~\eqref{eq:gen3} we obtain:
		\begin{equation} \label{eq:gen4}
			\begin{aligned}
				v_2(t)=&x_0\,\exp_g(\lambda;0,t)+
				\left( v_0-\lambda x_0\right) \, \exp_g(\lambda;0,t)\, \int_{[0,t)} \frac{1}{1+\lambda \Delta^+g(s)}\,
				\operatorname{d} \mu_g(s)  \\&
				+ \exp_g(\lambda;0,t)\, \int_{[0,t)} \frac{1}{1+\lambda \, \Delta^+g(s)}\, \operatorname{d} \mu_g(s) \,
				\int_{[0,t)} \frac{\exp_g(\lambda;0,s)^{-1}}{1+\lambda \Delta^+g(s)}\, f(s)\, \operatorname{d} \mu_g(s)
				\\&
				- \exp_g(\lambda;0,t)\,  \int_{[0,t)} \left(\int_{[0,s)} \frac{1}{1+\lambda \, \Delta^+g(r)}\, \operatorname{d} \mu_g(r) \right)\,
				\frac{\exp_g(\lambda;0,s)^{-1}}{1+\lambda \Delta^+g(s)}\, f(s)\, \operatorname{d} \mu_g(s)  \\&
				- \exp_g(\lambda;0,t)\,  \int_{[0,t)}
				\frac{\exp_g(\lambda;0,s)^{-1}}{(1+\lambda \Delta^+g(s))^2}\, f(s)\, \Delta^+g(s)\, \operatorname{d} \mu_g(s).
		\end{aligned}\end{equation}
		Observe that
		\begin{equation}%
			v_h(t)=x_0\,\exp_g(\lambda;0,t)+
			\left( v_0-\lambda x_0\right) \, \exp_g(\lambda;0,t)\, \int_{[0,t)} \frac{1}{1+\lambda \Delta^+g(s)}\,
			\operatorname{d} \mu_g(s) \in \mathcal{BC}_g^{\infty}([0,T];{\mathbb F})\end{equation}
		is the solution of the homogeneous equation~\eqref{eq:secondordernh} and
		\begin{equation}\label{eq:particular2}
			\begin{aligned}
				v_p(t)=&\exp_g(\lambda;0,t)\, \int_{[0,t)} \frac{1}{1+\lambda \, \Delta^+g(s)}\, \operatorname{d} \mu_g(s) \,
				\int_{[0,t)} \frac{\exp_g(\lambda;0,s)^{-1}}{1+\lambda \Delta^+g(s)}\, f(s)\, \operatorname{d} \mu_g(s)
				\\&
				- \exp_g(\lambda;0,t)\,  \int_{[0,t)} \left(\int_{[0,s)} \frac{1}{1+\lambda \, \Delta^+g(r)}\, \operatorname{d} \mu_g(r) \right)\,
				\frac{\exp_g(\lambda;0,s)^{-1}}{1+\lambda \Delta^+g(s)}\, f(s)\, \operatorname{d} \mu_g(s)  \\&
				- \exp_g(\lambda;0,t)\,  \int_{[0,t)}
				\frac{\exp_g(\lambda;0,s)^{-1}}{(1+\lambda \Delta^+g(s))^2}\, f(s)\, \Delta^+g(s)\, \operatorname{d} \mu_g(s)
		\end{aligned}\end{equation}
		is a particular solution in the space $\mathcal{BC}_g^{n+2}([0,T],{\mathbb F})$ of the non
		homogeneous equation~\eqref{eq:secondordernh} that satisfies $v_p(0)=(v_p)'_g(0)=0$.
	\end{itemize}
\end{rem}

\section{The Stieltjes harmonic oscillator} \label{SSHO}

In this section we present an application related to the real solution of
the Stieltjes harmonic oscillator ($g$-harmonic oscillator). 
{ Let 
$g:\mathbb{R}\to \mathbb{R}$ be a derivator 
such that $0\notin N_g^-$ and $T\notin N_g^+\cup D_g 
\cup C_g$ and denote by $g^C$ its continuous part. We consider the following equation:}
\begin{equation}\label{eq:gharmonic}
	\left\{\begin{aligned}
		&v_g^{\prime \prime}(t)+2\, \zeta\, \omega_0\, v_g^{\prime}(t)+ \omega_0^2 \, v(t)= 0, \;
		g-a.e.\, t \in [0,T), \\
		& v(0)=x_0,\\
		& v_g^{\prime}(0)=v_0,
	\end{aligned}\right.\end{equation}
where $x_0$ and $v_0$ are real numbers and:
\begin{itemize}
	\item $\omega_0$ is the undamped angular frequency of the oscillator:
	\begin{equation}%
		\omega_0=\sqrt{\frac{k}{m}},\end{equation}
	being $m$ the mass of the oscillator and $k$ a measure of the stiffness of the spring,
	two positive constants.
	\item $\zeta$ is the damping ratio:
	\begin{equation}%
		\zeta=\frac{c}{2\sqrt{mk}}\end{equation}
	being $c>0$ the viscous damping coefficient (resistance of the medium).
	If $\zeta>1$ we have an overdamped
	oscillator, if $\zeta=1$ the oscillator is critically damped and, if $\zeta<1$, the oscillator
	is underdamped. Observe that the solutions of the characteristic equation are given
	by:
	\begin{equation}%
		\begin{aligned}
			\lambda=&\frac{1}{2}\left(-2\, \zeta\, \omega_0 \pm \sqrt{4\,\zeta^2 \, \omega_0^2-
				4\, \omega_0^2}\right)=-\zeta \, \omega_0 \pm \omega_0 \, \sqrt{\zeta^2-1}.
	\end{aligned}\end{equation}
	Assume that $1+\lambda \, \Delta^+g(t)\neq 0$ for all $t\in [0,T)\cap D_g$ and for all
	$\lambda$ solution of the characteristic equation.
\end{itemize}
We have the following real solution of the $g$-harmonic oscillator in terms of
the damping ratio:
\begin{itemize}
	\item If $\zeta>1$ we have two real solution of the characteristic equation:
	\begin{equation}%
		\begin{aligned}
			\lambda_1 =& -\zeta \, \omega_0 - \omega_0 \, \sqrt{\zeta^2-1},\\
			\lambda_2 =& -\zeta \, \omega_0 + \omega_0 \, \sqrt{\zeta^2-1},
	\end{aligned}\end{equation}
	thus, the solution of~\eqref{eq:gharmonic} is given by:
	\begin{equation}%
		\begin{aligned}
			v(t)=&\left(\frac{v_{0}-\lambda _{2}\,x_{0}}{\lambda _{1}-\lambda _{2}}\right)\, \exp_g(\lambda_1;0,t)
			-\left(\frac{v_{0}-\lambda _{1}\,x_{0}}{\lambda _{1}-\lambda _{2}}\right) \, \exp_g(\lambda_2;0,t),
			\\
			=& \left(\frac{v_{0}-\lambda _{2}\,x_{0}}{\lambda _{1}-\lambda _{2}}\right) \,
			\exp(\lambda_1\, g_C(t))\, \prod_{s \in [0,t)\cap D_g} \left(
			1+\lambda_1\, \Delta^+g(s) \right)  \\
			&-\left(\frac{v_{0}-\lambda _{1}\,x_{0}}{\lambda _{1}-\lambda _{2}}\right) \,
			\exp_g(\lambda_2\, g_C(t)) \, \prod_{s \in [0,t)\cap D_g} \left(
			1+\lambda_2 \, \Delta^+g(s)\right) .
	\end{aligned}\end{equation}
	\item If $\zeta=1$ we have one real solution of the characteristic equation
	\begin{equation}%
		\lambda= -\zeta\, \omega_0,\end{equation}
	thus, the solution of~\eqref{eq:gharmonic} is given by:
	\begin{equation}%
		\begin{aligned}
			v(t)=&x_0\, \exp_g(\lambda;0,t)+(v_0-\lambda\, x_0)\, \exp_g(\lambda;0,t)\, \int_{[0,t)}
			\frac{1}{1+\lambda\, \Delta^+g(s)}\, \operatorname{d} \mu_g,\\
			=& \exp(\lambda\, g_C(t))\, \prod_{s\in [0,t)\cap D_g}
			\left( 1+\lambda\, \Delta^+g(s)\right)  \\ & \cdot \bigg[ x_0 +
			(v_0-\lambda\, x_0) \, \bigg(g_C(t)+
			\sum_{s \in [0,t)\cap D_g} \frac{\Delta^+g(s)}{1+\lambda\, \Delta^+g(s)}\bigg)\bigg].
	\end{aligned}\end{equation}
	\item If $\zeta<1$ we have a conjugate complex solution:
	\begin{equation}%
		\begin{aligned}
			\lambda_1=&-\zeta\, \omega_0 + i\, \omega_0\, \sqrt{1-\zeta^2}, \\
			\lambda_2 =&-\zeta\, \omega_0 - i\, \omega_0\, \sqrt{1-\zeta^2},
	\end{aligned}\end{equation}
	We will denote by $a=-\zeta\, \omega_0$ and $b=\omega_0\, \sqrt{1-\zeta^2}$, we have
	that the solution is given by the following formula:
	\begin{equation}%
		\begin{aligned}
			v(t)=&\left(
			\frac{v_0-a\, x_0+ b\,i\,x_0}{ 2\, b \, i}\right) \, \exp_g(a+b\,i;0,t)-
			\left(\frac{v_0-a\, x_0-b\,i\, x_0}{2\, b\, i}\right)\,
			\exp_g(a-b\,i;0,t).
	\end{aligned}\end{equation}
	If we take into account expression~\eqref{eq:expgab}:
	\begin{equation}%
		\begin{aligned}
			v(t)=&\exp_g(a;0,t)\, \bigg[\left( \frac{v_0-a\, x_0}{b}\right) \, \frac{1}{2\,i}\,
			\bigg(\exp_g\left( \frac{i\, b}{1+a\, \Delta^+g};0,t\right) -\exp_g\left( \frac{-i\, b}{1+a\, \Delta^+g};0,t\right)
			\bigg)\\
			&+x_0 \, \frac{1}{2} \bigg(
			\exp_g\left( \frac{i\, b}{1+a\, \Delta^+g};0,t\right) +\exp_g\left( \frac{-i\, b}{1+a\, \Delta^+g};0,t\right)
			\bigg)\bigg] \\
			=& \exp_g(a;0,t) \,\bigg[ \left( \frac{v_0-a\, x_0}{b}\right) \,\sin_g\left( \frac{b}{1+a\, \Delta^+g};0,t\right)
			+ x_0 \, \cos_g\left( \frac{b}{1+a\, \Delta^+g};0,t\right) \bigg].
	\end{aligned}\end{equation}
\end{itemize}

\begin{exa} \label{example1} Let us analyze the behavior of the $g$-harmonic
	oscillator in the particular
	case $\omega_0=2$, $x_0=v_0=1$. We consider the following derivators:
	\begin{equation}%
		\begin{aligned}
			g_1(t)=& g_1^C(t)+g^B(t), \\
			g_2(t)=&g_2^C(t)+g^B(t),
	\end{aligned}\end{equation}
	where $g_1^C(t)=t$,
	\begin{equation}%
		g_2(t)=\begin{dcases}
			\frac{1}{2}+\frac{[x]}{2},\; & [x]+1 = 2\, k,\; k \in \mathbb{N}, \\
			x-\frac{[x]}{2},\; & [x]+1 = 2\, k -1, \; k \in \mathbb{N},
	\end{dcases}\end{equation}
	$[\cdot]$ denotes the floor function and
	\begin{equation}%
		g^B(t)=\sum_{s\in[0,t)\cap D_g} \Delta^+ g(s),\end{equation}
	being $D_g=\{t \in [0,\infty):\; t = k\, \pi/4,\; k \in \mathbb{N}\}$ and
	$\Delta^+g(t)=l$ for some fixed value $l$ for all $t\in D_g$. In order to compare the effect
	that discontinuities in the derivator have on the solution, we have considered
	the cases where $l\in\{0,1/3^3,1/3^2,1/3^1,1/3^0\}$. Observe that the solution
	associated to the derivator $g_1$ and $l=0$ corresponds to
	the classical solution of the harmonic oscillator. In Figure~\ref{fig1} we
	can see a graphical representation of the derivators considered.
	\begin{figure}[h]
		\begin{subfigure}{.5\textwidth}
			\centering
			\includegraphics[width=1\linewidth]{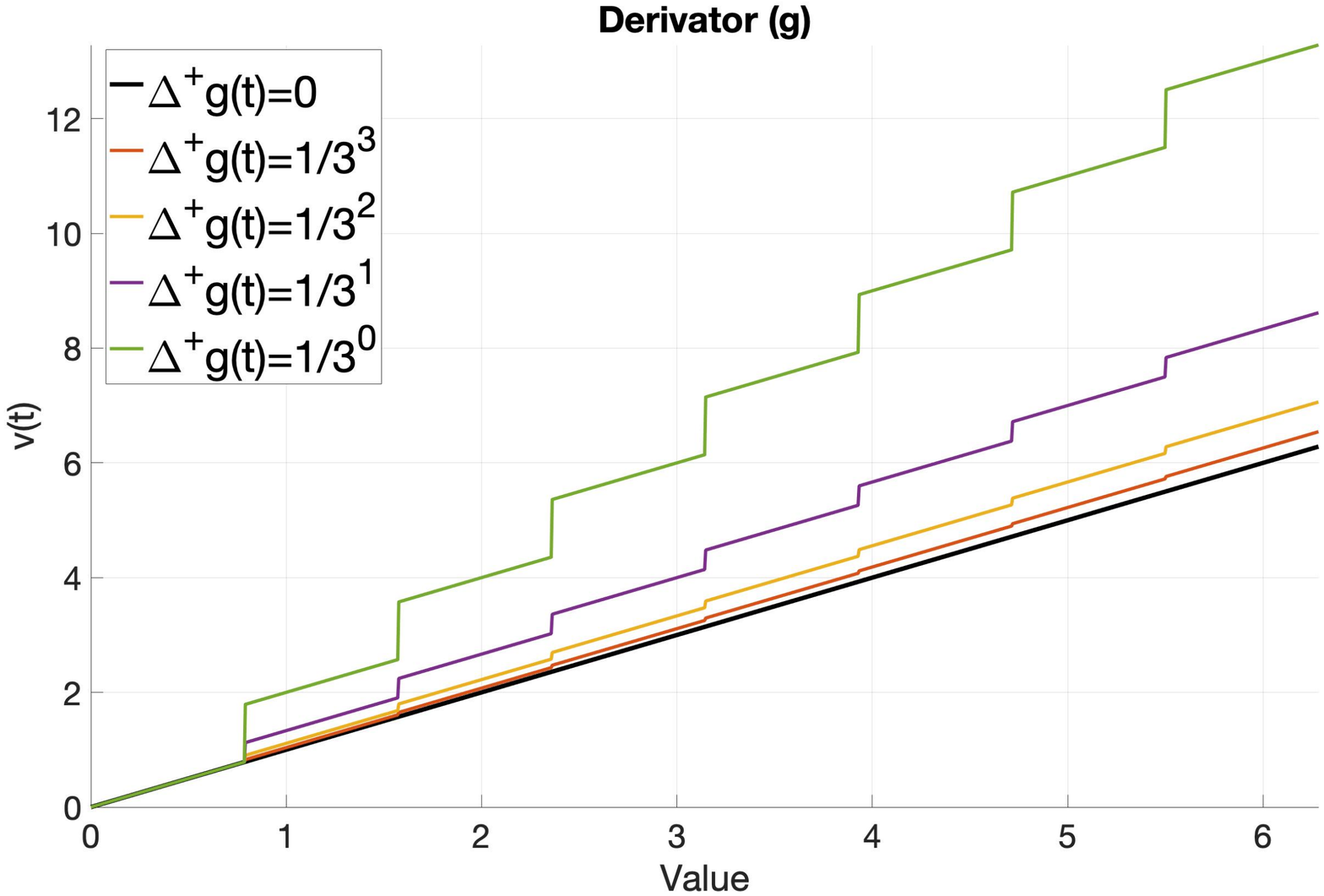}
			\caption{Derivator $g_1$.}
		\end{subfigure}
		\begin{subfigure}{.5\textwidth}
			\centering
			\includegraphics[width=1\linewidth]{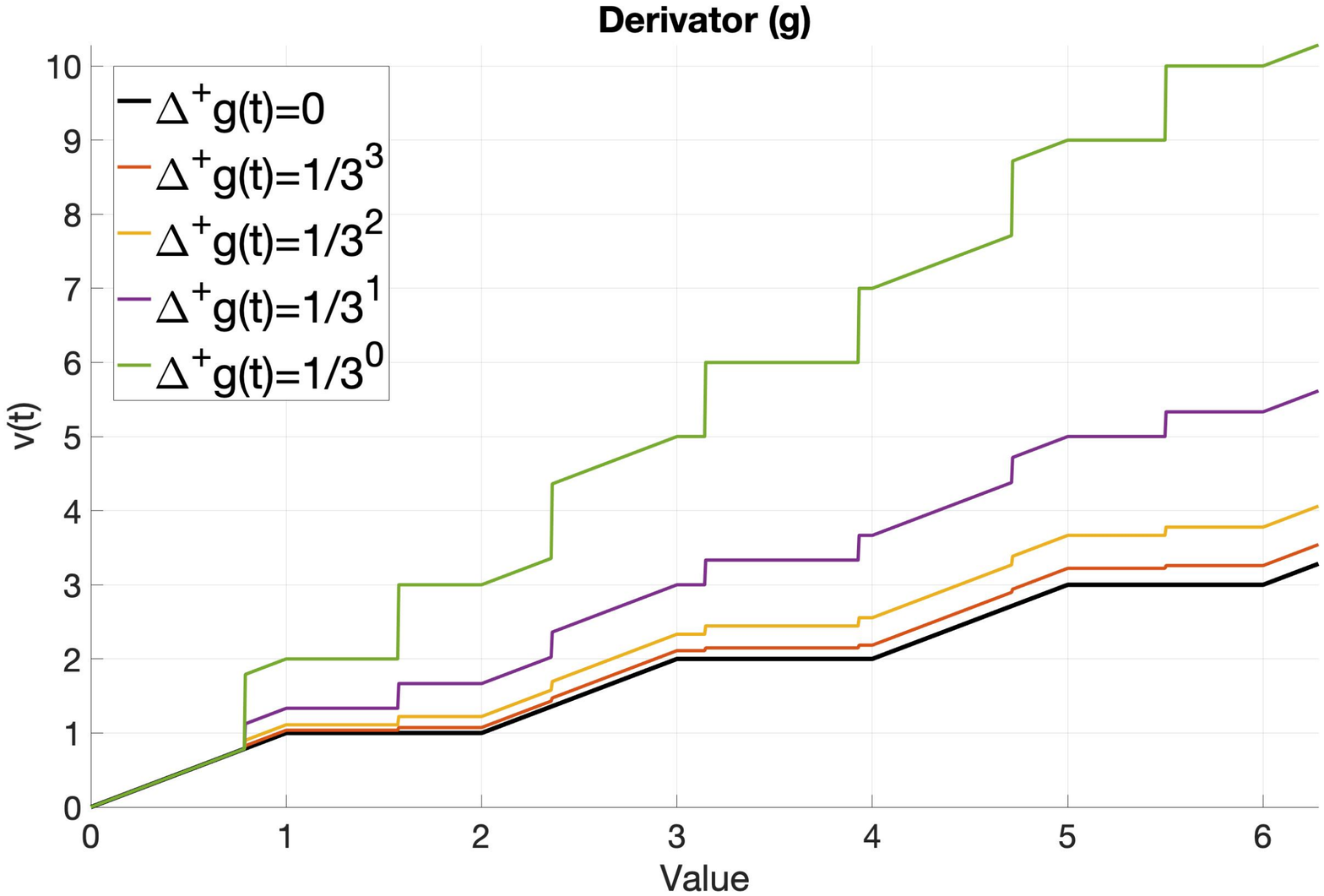}
			\caption{Derivator $g_2$.}
		\end{subfigure}
		\vspace{+6pt}
		\caption{Derivators $g_1$ and $g_2$ associated to $l\in\{0,1/3^3,1/3^2,1/3^1,1/3^0\}$
			{(vertical lines have to be understood as jumps and not as a multivalued function)}.}
		\label{fig1}
	\end{figure}
	Now, we present a graphical representation of the solution associated to
	$\zeta=1.5$ (overdamped oscillator, Figure~\ref{fig2}), $\zeta=1$
	(critically damped oscillator, Figure~\ref{fig3}) and $z=0.5$
	(underdamped oscillator, Figure~\ref{fig4}).

	\begin{figure}[h]
		\begin{subfigure}{.5\textwidth}
			\centering
			\includegraphics[width=1\linewidth]{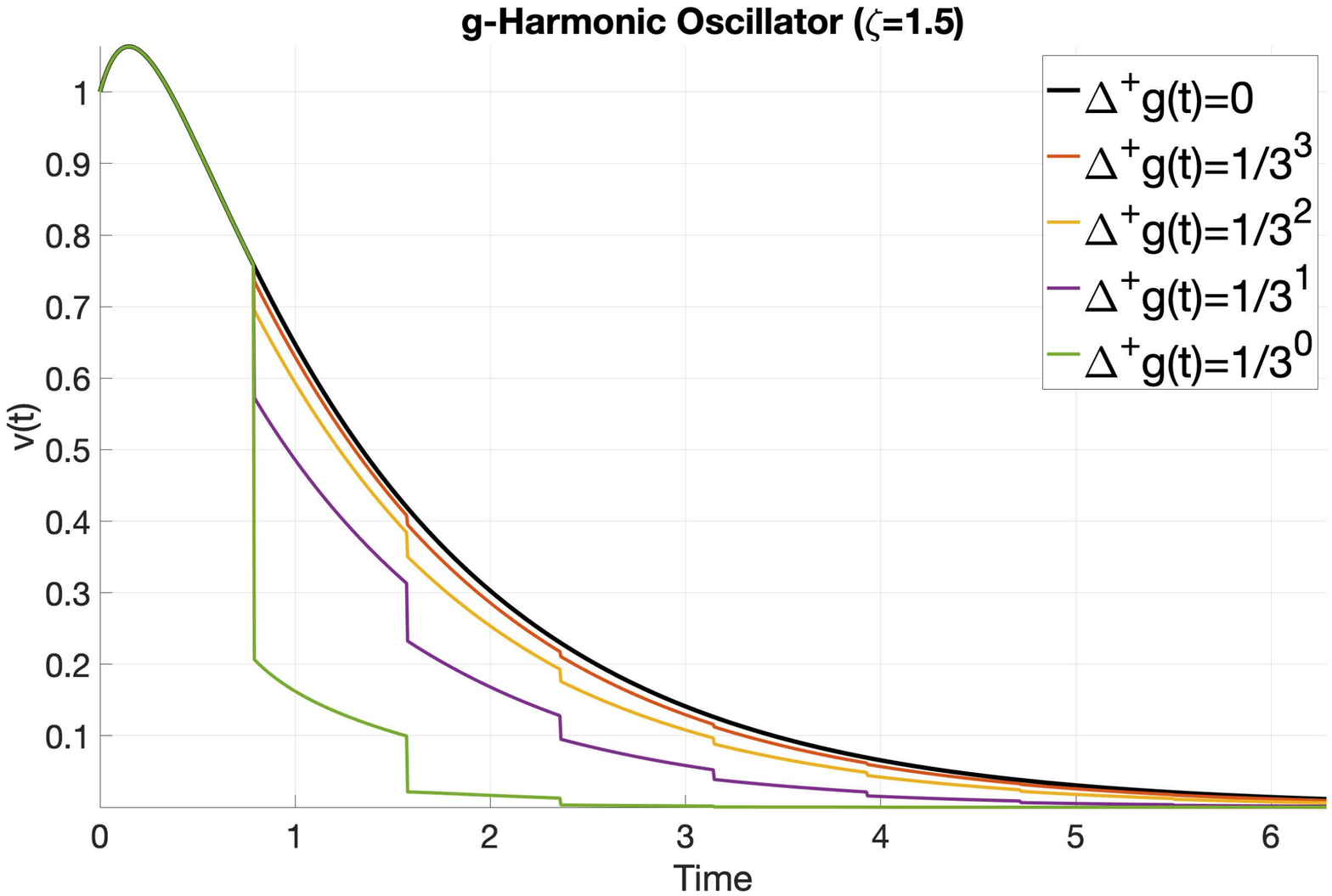}
			\caption{$\zeta=1.5$, overdamped oscillator ($g_1$).}
		\end{subfigure}
		\begin{subfigure}{.5\textwidth}
			\centering
			\includegraphics[width=1\linewidth]{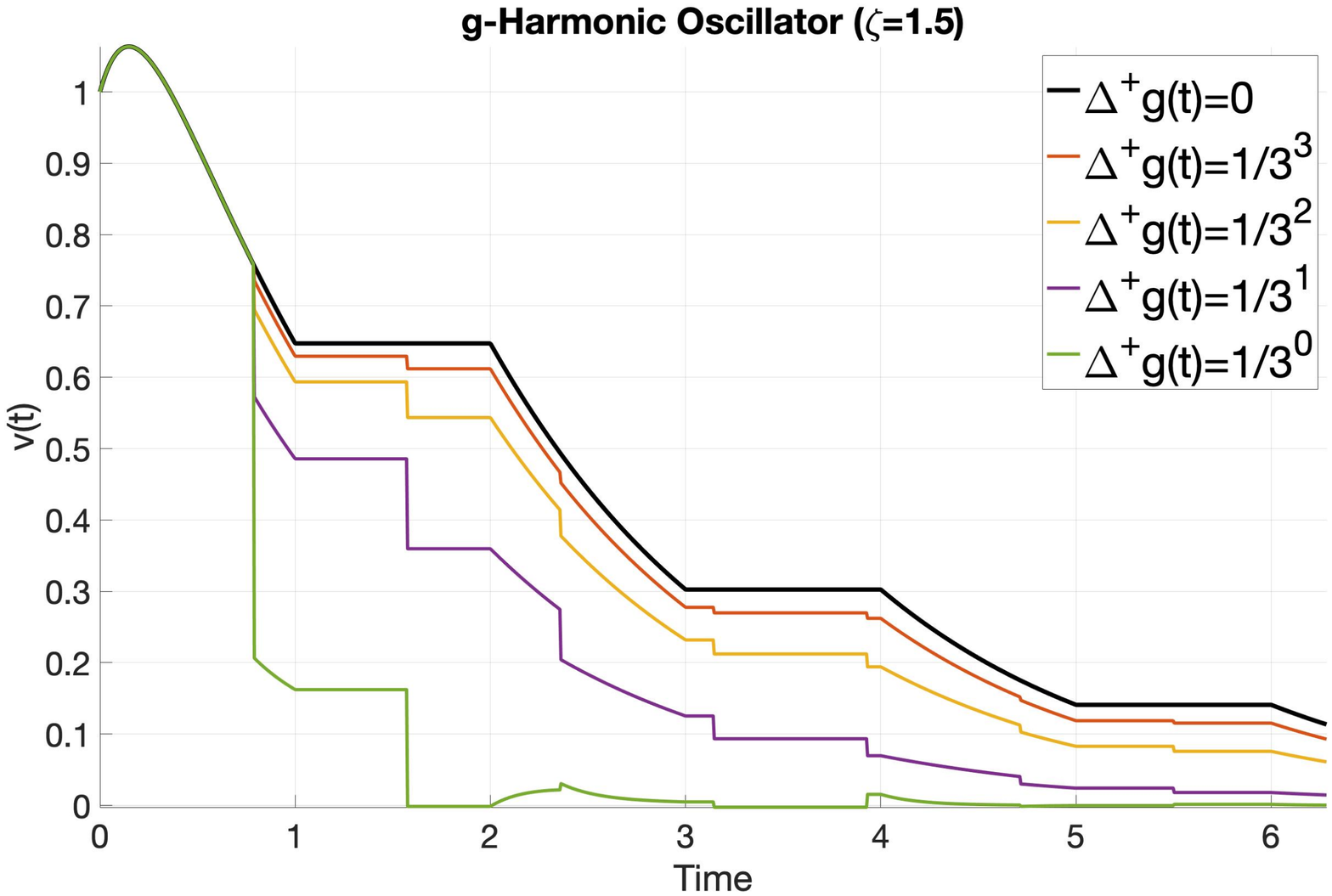}
			\caption{$\zeta=1.5$, overdamped oscillator ($g_2$).}
		\end{subfigure}
		\vspace{+6pt}
		\caption{solution of the $g$-harmonic oscillator for $\zeta=1.5$ (overdamped oscillator)
			associated to $g_1$ and $g_2$ {(vertical lines have to be understood as jumps and not as a multivalued function)}.}
		\label{fig2}
	\end{figure}

	\begin{figure}[h]
		\begin{subfigure}{.5\textwidth}
			\centering
			\includegraphics[width=1\linewidth]{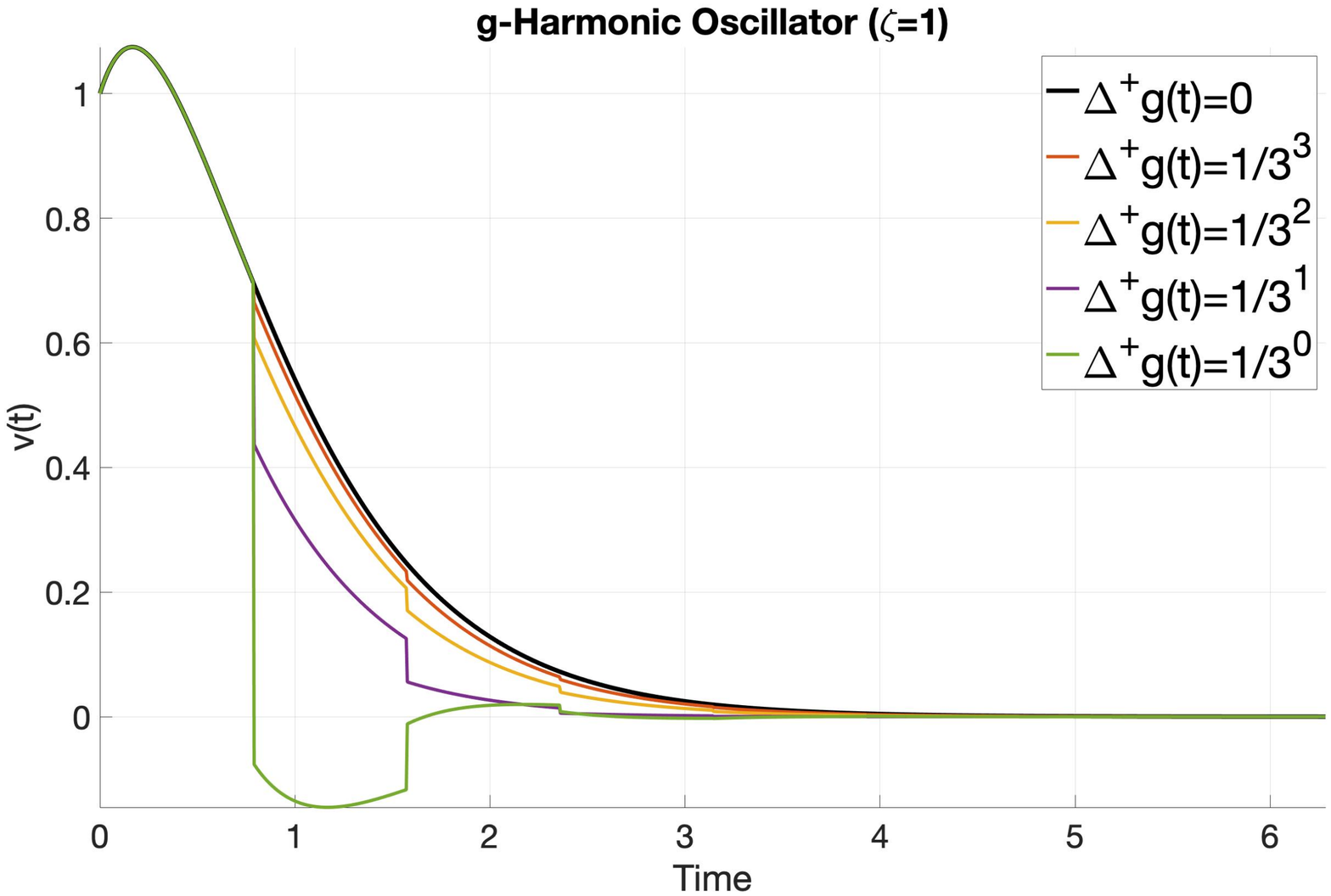}
			\caption{$\zeta=1$, critically damped oscillator ($g_1$).}
		\end{subfigure}
		\begin{subfigure}{.5\textwidth}
			\centering
			\includegraphics[width=1\linewidth]{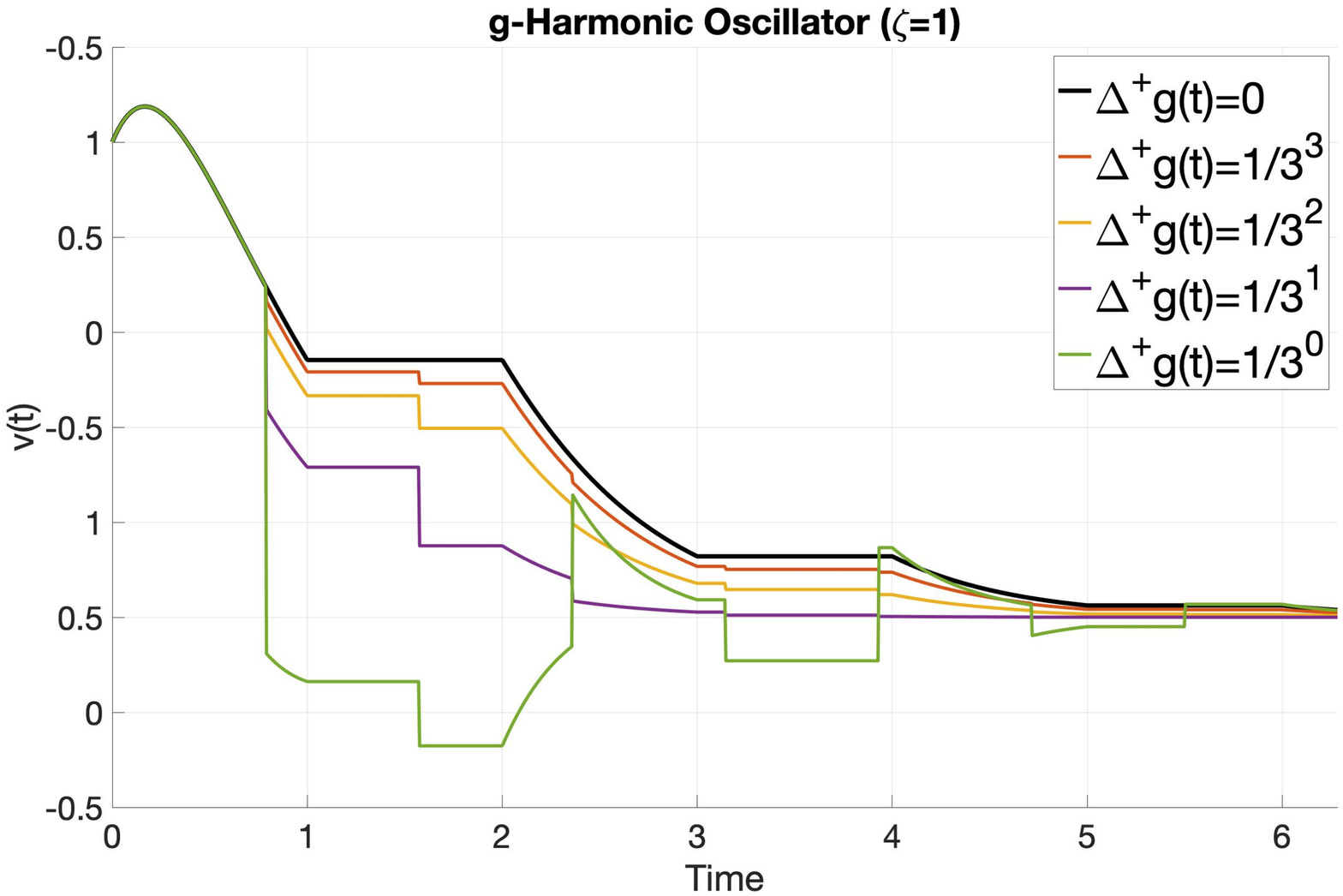}
			\caption{$\zeta=1$, critically damped oscillator ($g_2$).}
		\end{subfigure}
		\vspace{+6pt}
		\caption{solution of the $g$-harmonic oscillator for $\zeta=1$ (critically damped oscillator)
			associated to $g_1$ and $g_2$ {(vertical lines have to be understood as jumps and not as a multivalued function)}.}
		\label{fig3}
	\end{figure}

	\begin{figure}[h]
		\begin{subfigure}{.5\textwidth}
			\centering
			\includegraphics[width=1\linewidth]{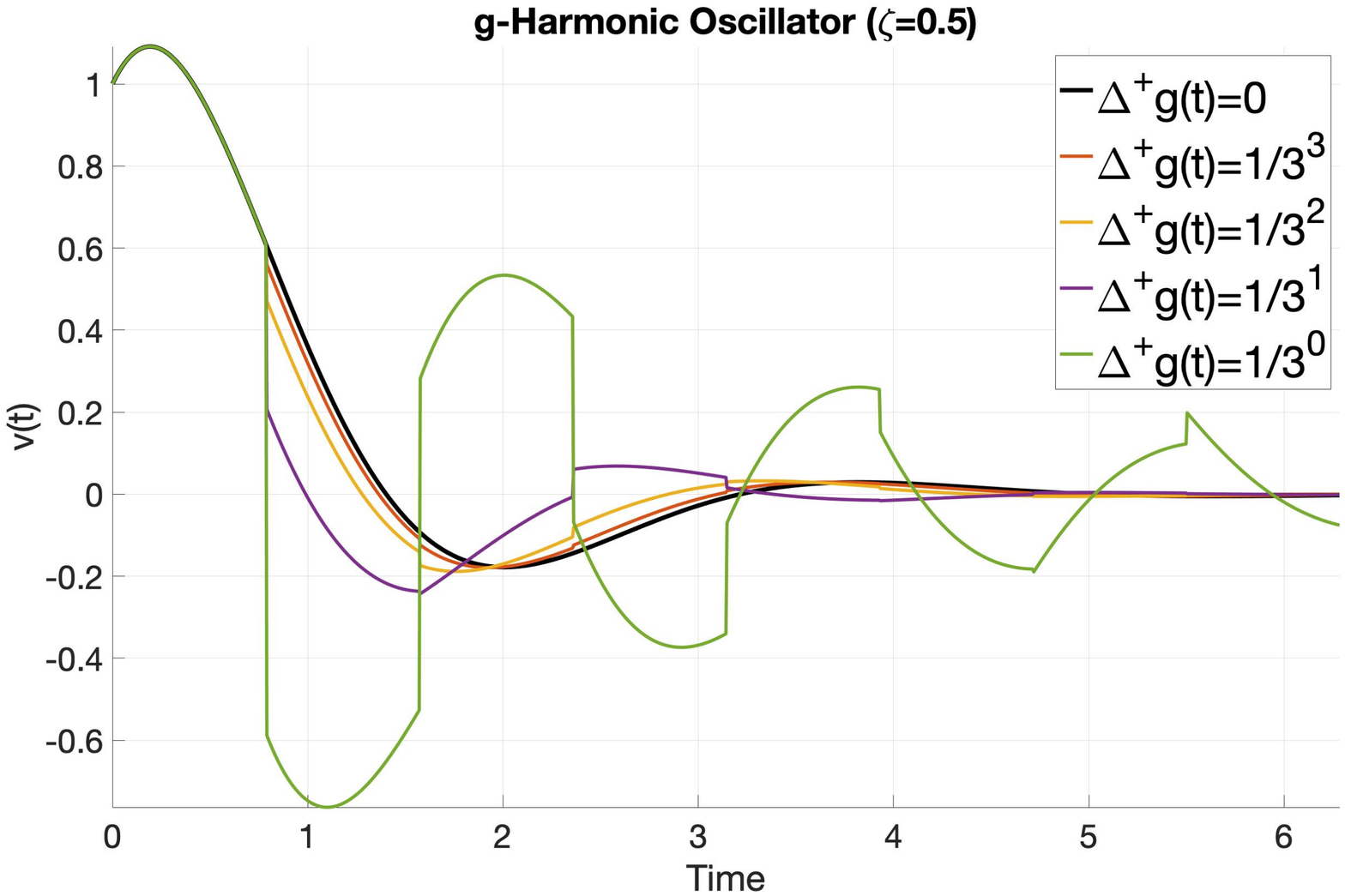}
			\caption{$\zeta=0.5$, underdamped oscillator ($g_1$).}
		\end{subfigure}
		\begin{subfigure}{.5\textwidth}
			\centering
			\includegraphics[width=1\linewidth]{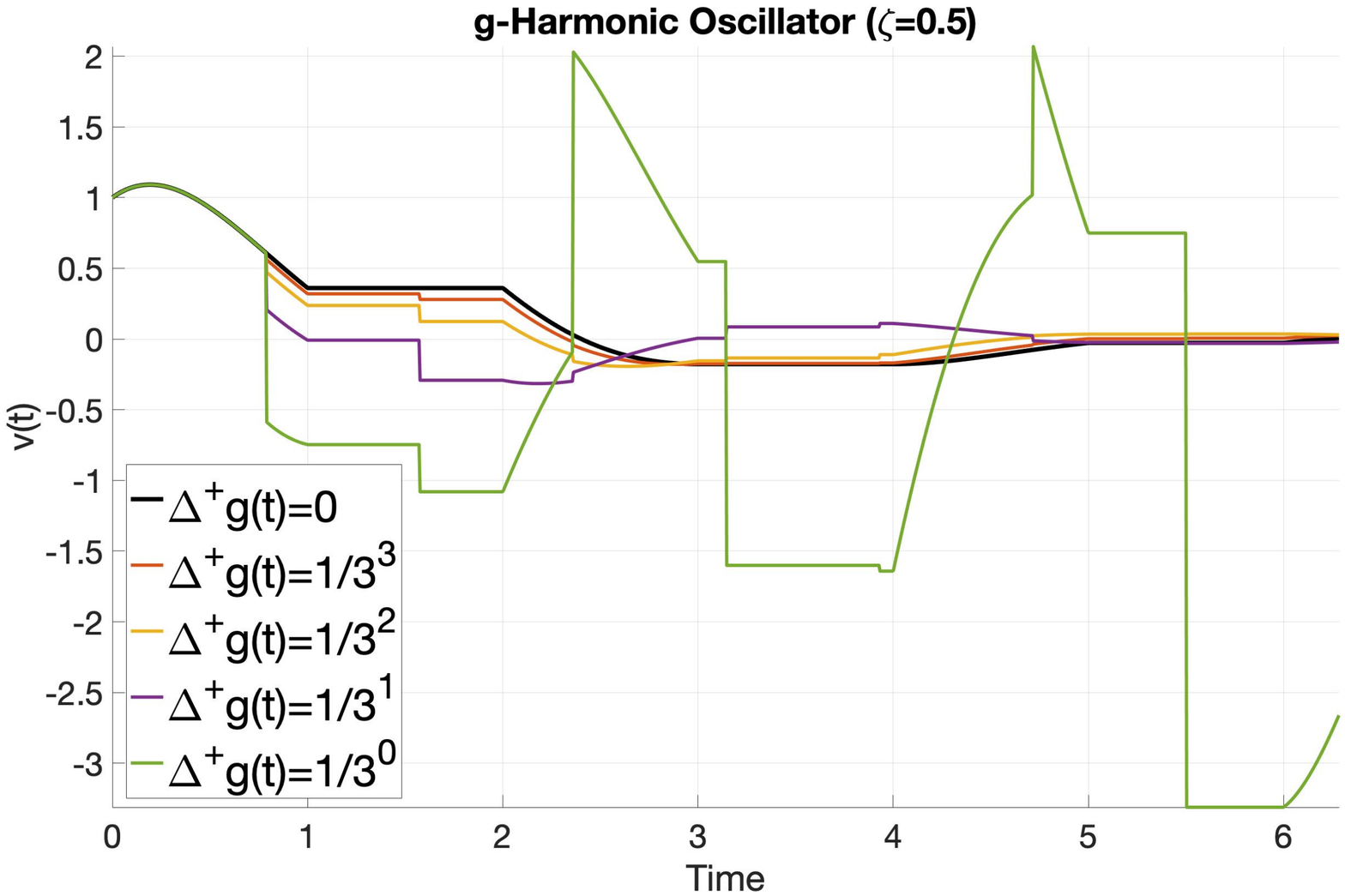}
			\caption{$\zeta=0.5$, underdamped oscillator ($g_2$).}
		\end{subfigure}
		\vspace{+6pt}
		\caption{solution of the $g$-harmonic oscillator for $\zeta=0.5$ (underdamped oscillator)
			associated to $g_1$ and $g_2$ {(vertical lines have to be understood as jumps and not as a multivalued function)}.}
		\label{fig4}
	\end{figure}
\end{exa}

Now we will study the effect of a $g$-periodic source term with the same frequency
as the natural frequency of the oscillator. We will consider the non-homogeneous
$g$-harmonic oscillator with $c=0$ and we will consider as a source
term $f(t)=\cos_g(\omega_0;0,t)$.
\begin{equation}\label{eq:gharmonic2}
	\left\{\begin{aligned}
		&v_g^{\prime \prime}(t)+ \omega_0^2 \, v(t)= \cos_g(\omega_0;0,t), \;
		g-a.e.\, t \in [0,T), \\
		& v(0)=x_0,\\
		& v_g^{\prime}(0)=v_0,
	\end{aligned}\right.\end{equation}
Observe that in the case $g(t)=t$ we recover the classical resonance effect. It is reasonable to expect that in the case of having a generic derivator
the amplitude of the oscillations increases with $t$. Indeed we have that
the solution $v$ is given by formula~\eqref{eq:gen2}:
\begin{equation}%
	v(t)=v_h(t)+v_p(t),\end{equation}
where $v_h$ is the solution of the homogeneous equation:
\begin{equation}%
	\begin{aligned}
		v_h(t)=&\left(
		\frac{v_0+i\, \omega_0\,x_0}{ 2\,i \, \omega_0} \right) \, \exp_g(i\,\omega_0;0,t)
		- \left(\frac{v_0-i\,\omega_0\, x_0}{2\,i\, \omega_0}\right)\,
		\exp_g(-i\, \omega_0;0,t) \\
		=&x_0\, \cos_g(\omega_0;0,t)+ \frac{v_0}{\omega_0}\,\sin_g(\omega_0;0,t),
\end{aligned}\end{equation}
and $v_p$ is a particular solution that satisfies $v_p(0)=(v_p)^{\prime}_g(t)=0$ given by
\begin{equation}%
	\begin{aligned}
		v_p(t)=&+ (2\,i\, \omega_0)^{-1} \exp_g(+i\,\omega_0;0,t)
		\int_{[0,t)} \frac{\exp_g(+i\, \omega_0 ;0,s)^{-1}}{1+i\, \omega_0 \, \Delta^+ g(s)}\,
		\cos_g(\omega_0;0,s)\, \operatorname{d} \mu_g
		\\
		&\displaystyle-(2\,i\,\omega_0)^{-1} \exp_g(-i\, \omega_0;0,t) \int_{[0,t)}
		\frac{\exp_g(-i\, \omega_0;0,s)^{-1}}{1-i\,\omega_0 \, \Delta^+g(s)}\,
		\cos_g(\omega_0;0,s)\, \operatorname{d} \mu_g.
\end{aligned}\end{equation}
Now,
\begin{equation}%
	\begin{aligned}
		&\int_{[0,t)} \frac{\exp_g(+i\, \omega_0 ;0,s)^{-1}}{1+i\, \omega_0 \, \Delta^+ g(s)}\,
		\cos_g(\omega_0;0,s)\, \operatorname{d} \mu_g \\
		=&\frac{1}{2}\int_{[0,t)} \frac{\exp_g(+i\, \omega_0 ;0,s)^{-1}}{1+i\, \omega_0 \, \Delta^+ g(s)}\,
		\left( \exp_g(i\,\omega_0;0,s)+\exp_g(-i\,\omega_0;0,s)\right) \, \operatorname{d} \mu_g\\
		=&\frac{1}{2}\, \int_{[0,t)} \frac{1}{1+i\, \omega_0 \, \Delta^+ g(s)}\, \operatorname{d} \mu_g+
		\frac{1}{2}\,\int_{[0,t)} \exp_g(-i\,\omega_0;0,s) \,
		\frac{\exp_g(+i\, \omega_0 ;0,s)^{-1}}{1+i\, \omega_0 \, \Delta^+ g(s)}\,\operatorname{d} \mu_g \\
		=&\frac{1}{2}\, \int_{[0,t)} \frac{1}{1+i\, \omega_0 \, \Delta^+ g(s)}\, \operatorname{d} \mu_g+
		\frac{1}{2}\, \int_{[0,t)} \frac{1}{1+i\, \omega_0\, \Delta^+g(s)} \,
		\exp_g\left( \frac{-2\,i \, \omega_0}{1+i\, \omega_0\, \Delta^+g};0,s\right) \, \operatorname{d} \mu_g \\
		=& \frac{1}{2}\, \int_{[0,t)} \frac{1}{1+i\, \omega_0 \, \Delta^+ g(s)}\, \operatorname{d} \mu_g-
		\frac{1}{4\,i\,\omega_0}\,\bigg[
		\exp_g\left( \frac{-2\,i \, \omega_0}{1+i\, \omega_0\, \Delta^+g};0,t\right) -1\bigg].
\end{aligned}\end{equation}
On the other hand,
\begin{equation}%
	\begin{aligned}
		& \int_{[0,t)}
		\frac{\exp_g(-i\, \omega_0;0,s)^{-1}}{1-i\,\omega_0 \, \Delta^+g(s)}\,
		\cos_g(\omega_0;0,s)\, \operatorname{d} \mu_g \\
		=&\frac{1}{2} \, \int_{[0,t)}
		\frac{\exp_g(-i\, \omega_0;0,s)^{-1}}{1-i\,\omega_0 \, \Delta^+g(s)}\,
		\left( \exp_g(+i\,\omega_0;0,s)+\exp_g(-i\,\omega_0;0,s)\right) \, \operatorname{d} \mu_g\\
		=& \frac{1}{2} \, \int_{[0,t)} \exp_g(+i\,\omega_0;0,s)\,
		\frac{\exp_g(-i\, \omega_0;0,s)^{-1}}{1-i\,\omega_0 \, \Delta^+g(s)}\,\operatorname{d} \mu_g
		+\frac{1}{2}\, \int_{[0,t)} \frac{1}{1-i\,\omega_0\, \Delta^+g(s)}\, \operatorname{d} \mu_g \\
		=& \frac{1}{2} \int_{[0,t)} \frac{1}{1-i\,\omega_0 \, \Delta^+g(s)}\,
		\exp_g\left( \frac{2\,i\, \omega_0}{1-i\,\omega_0\, \Delta^+g};0,s \right) \, \operatorname{d} \mu_g
		+\frac{1}{2}\, \int_{[0,t)} \frac{1}{1-i\,\omega_0\, \Delta^+g(s)}\, \operatorname{d} \mu_g \\
		=&\frac{1}{4\,i\, \omega_0}\bigg[
		\exp_g\left( \frac{2\,i\, \omega_0}{1-i\,\omega_0\, \Delta^+g};0,t \right) -1
		\bigg]+\frac{1}{2}\, \int_{[0,t)} \frac{1}{1-i\,\omega_0\, \Delta^+g(s)}\, \operatorname{d} \mu_g.
\end{aligned}\end{equation}
Thus,
\begin{equation}%
	\begin{aligned}
		v_p(t)=&+ \frac{1}{4\,i\, \omega_0} \, \exp_g(+i\,\omega_0;0,t) \,
		\int_{[0,t)} \frac{1}{1+i\, \omega_0 \, \Delta^+ g(s)}\, \operatorname{d} \mu_g\\
		&+\frac{1}{8\, \omega_0^2} \, \exp_g(+i\,\omega_0;0,t) \,\bigg[
		\exp_g\left( \frac{-2\,i \, \omega_0}{1+i\, \omega_0\, \Delta^+g};0,t\right) -1\bigg]\\
		& - \frac{1}{4\,i\,\omega_0} \, \exp_g(-i\, \omega_0;0,t)\,
		\int_{[0,t)} \frac{1}{1-i\,\omega_0\, \Delta^+g(s)}\, \operatorname{d} \mu_g
		\\
		& + \frac{1}{8\, \omega_0^2}\,\exp_g(-i\, \omega_0;0,t)\, \bigg[
		\exp_g\left( \frac{2\,i\, \omega_0}{1-i\,\omega_0\, \Delta^+g};0,t \right) -1
		\bigg].
\end{aligned}\end{equation}
Now, taking into account that:
\begin{equation}%
	\begin{aligned}
		\exp_g\left( \frac{-2\,i \, \omega_0}{1+i\, \omega_0\, \Delta^+g};0,t\right)  =&
		\exp_g(-i\, \omega_0;0,t) \, \exp_g(+i\, \omega_0;0,t)^{-1} ,\\
		\exp_g\left( \frac{2\,i\, \omega_0}{1-i\,\omega_0\, \Delta^+g};0,t \right)  =&
		\exp_g(+i\,\omega_0;0,s)\, \exp_g(-i\, \omega_0;0,s)^{-1},
\end{aligned}\end{equation}
we obtain the following expression for the particular solution:
\begin{equation}%
	\begin{aligned}
		v_p(t)=&+ \frac{1}{4\,i\, \omega_0} \, \exp_g(+i\,\omega_0;0,t) \,
		\int_{[0,t)} \frac{1}{1+i\, \omega_0 \, \Delta^+ g(s)}\, \operatorname{d} \mu_g\\&
		- \frac{1}{4\,i\,\omega_0} \, \exp_g(-i\, \omega_0;0,t)\,
		\int_{[0,t)} \frac{1}{1-i\,\omega_0\, \Delta^+g(s)}\, \operatorname{d} \mu_g \\
		& +\frac{1}{8\,\omega_0^2} \, \exp_g(-i\, \omega_0;0,t) -
		\frac{1}{8\,\omega_0^2} \, \exp_g(+i\,\omega_0;0,t) \\
		& + \frac{1}{8\, \omega_0^2}\,\exp_g(+i\,\omega_0;0,s)-
		\frac{1}{8\,\omega_0^2}\, \exp_g(-i\, \omega_0;0,t).
\end{aligned}\end{equation}
Thus,
\begin{equation}%
	\begin{aligned}
		v_p(t)=&+ \frac{1}{4\,i\, \omega_0} \, \exp_g(+i\,\omega_0;0,t) \,
		\int_{[0,t)} \frac{1}{1+i\, \omega_0 \, \Delta^+ g(s)}\, \operatorname{d} \mu_g\\&
		- \frac{1}{4\,i\,\omega_0} \, \exp_g(-i\, \omega_0;0,t)\,
		\int_{[0,t)} \frac{1}{1-i\,\omega_0\, \Delta^+g(s)}\, \operatorname{d} \mu_g.
\end{aligned}\end{equation}
In order to simplify the previous expression let us consider the following
computations:
\begin{equation}%
	\begin{aligned}
		\int_{[0,t)} \frac{1}{1+i\, \omega_0 \, \Delta^+ g(s)}\, \operatorname{d} \mu_g=&
		\int_{[0,t)} \frac{1-i\,\omega_0\, \Delta^+g(s)}{1+\omega_0^2\, \Delta^+ g(s)^2}\,
		\operatorname{d} \mu_g\\
		=& \int_{[0,t)} \frac{1}{1+\omega_0^2\, \Delta^+ g(s)^2}\,
		\operatorname{d} \mu_g - i\, \omega_0\, \int_{[0,t)} \frac{\Delta^+g(s)}{1+\omega_0^2\, \Delta^+ g(s)^2}\,
		\operatorname{d} \mu_g,\\
		\int_{[0,t)} \frac{1}{1-i\,\omega_0\, \Delta^+g(s)}\, \operatorname{d} \mu_g=&
		\int_{[0,t)} \frac{1+i\,\omega_0\, \Delta^+g(s)}{1+\omega_0^2\, \Delta^+ g(s)^2}\,
		\operatorname{d} \mu_g \\
		=& \int_{[0,t)} \frac{1}{1+\omega_0^2\, \Delta^+ g(s)^2}\,
		\operatorname{d} \mu_g + i\, \omega_0\, \int_{[0,t)} \frac{\Delta^+g(s)}{1+\omega_0^2\, \Delta^+ g(s)^2}\,
		\operatorname{d} \mu_g.
\end{aligned}\end{equation}
Therefore,
\begin{equation}%
	\begin{aligned}
		v_p(t)=& +\frac{1}{2\,\omega_0}\,\int_{[0,t)} \frac{1}{1+\omega_0^2\, \Delta^+g(s)^2}\, \operatorname{d} \mu_g\,
		\bigg[ \frac{1}{2\,i}\, \big(\exp_g(+i\,\omega_0;0,t)-\exp_g(-i\, \omega_0;0,t)\big) \bigg]\\
		&-\frac{1}{2}\,  \int_{[0,t)} \frac{\Delta^+g(s)}{1+\omega_0^2\, \Delta^+ g(s)^2}\,
		\operatorname{d} \mu_g\, \bigg[\frac{1}{2}\,\big(\exp_g(+i\,\omega_0;0,t)+\exp_g(-i\,\omega_0;0,t)
		\big)
		\bigg] \\
		=& \frac{1}{2\,\omega_0} \,\sin_g(\omega_0;0,t) \,
		\int_{[0,t)} \frac{1}{1+\omega_0^2\, \Delta^+ g(s)^2}\, \operatorname{d} \mu_g\\&
		-\frac{1}{2}\, \cos_g(\omega_0;0,t)\,
		\int_{[0,t)} \frac{\Delta^+g(s)}{1+\omega_0^2\, \Delta^+ g(s)^2}\,
		\operatorname{d} \mu_g.
\end{aligned}\end{equation}
Finally, the solution of~\eqref{eq:gharmonic2} is given by
\begin{equation}\label{resosol}
	\begin{aligned}
		v(t)=&x_0\, \cos_g(\omega_0;0,t)+ \frac{v_0}{\omega_0}\,\sin_g(\omega_0;0,t)+
		\frac{1}{2\,\omega_0} \,\sin_g(\omega_0;0,t) \,
		\int_{[0,t)} \frac{1}{1+\omega_0^2\, \Delta^+ g(s)^2}\, \operatorname{d} \mu_g\\&
		-\frac{1}{2}\, \cos_g(\omega_0;0,t)\,
		\int_{[0,t)} \frac{\Delta^+g(s)}{1+\omega_0^2\, \Delta^+ g(s)^2}\,
		\operatorname{d} \mu_g.
\end{aligned}\end{equation}
Observe that if $g(t)=t$, we recover the classical solution and the amplitude
of oscillations grows with $t$.

\begin{exa} Let us take the same data and derivators considered in
	Example~\ref{example1}. In Figure~\ref{fig5} we can see
	the solution of~\eqref{eq:gharmonic2} for some of the derivators
	considered in Example~\ref{example1}.

	\begin{figure}[h]
		\begin{subfigure}{.5\textwidth}
			\centering
			\includegraphics[width=1\linewidth]{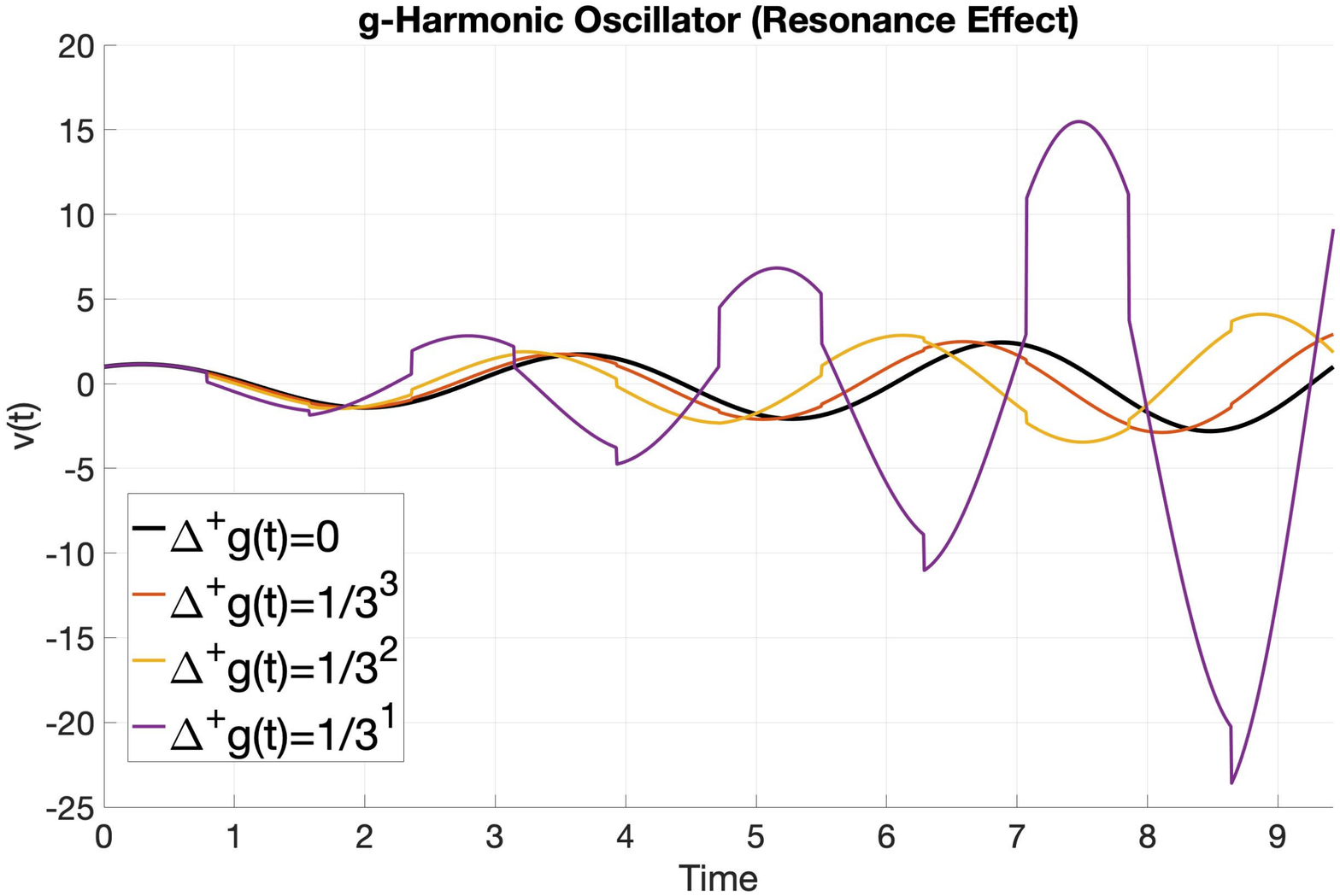}
			\caption{Resonance effect ($g_1$).}
		\end{subfigure}
		\begin{subfigure}{.5\textwidth}
			\centering
			\includegraphics[width=1\linewidth]{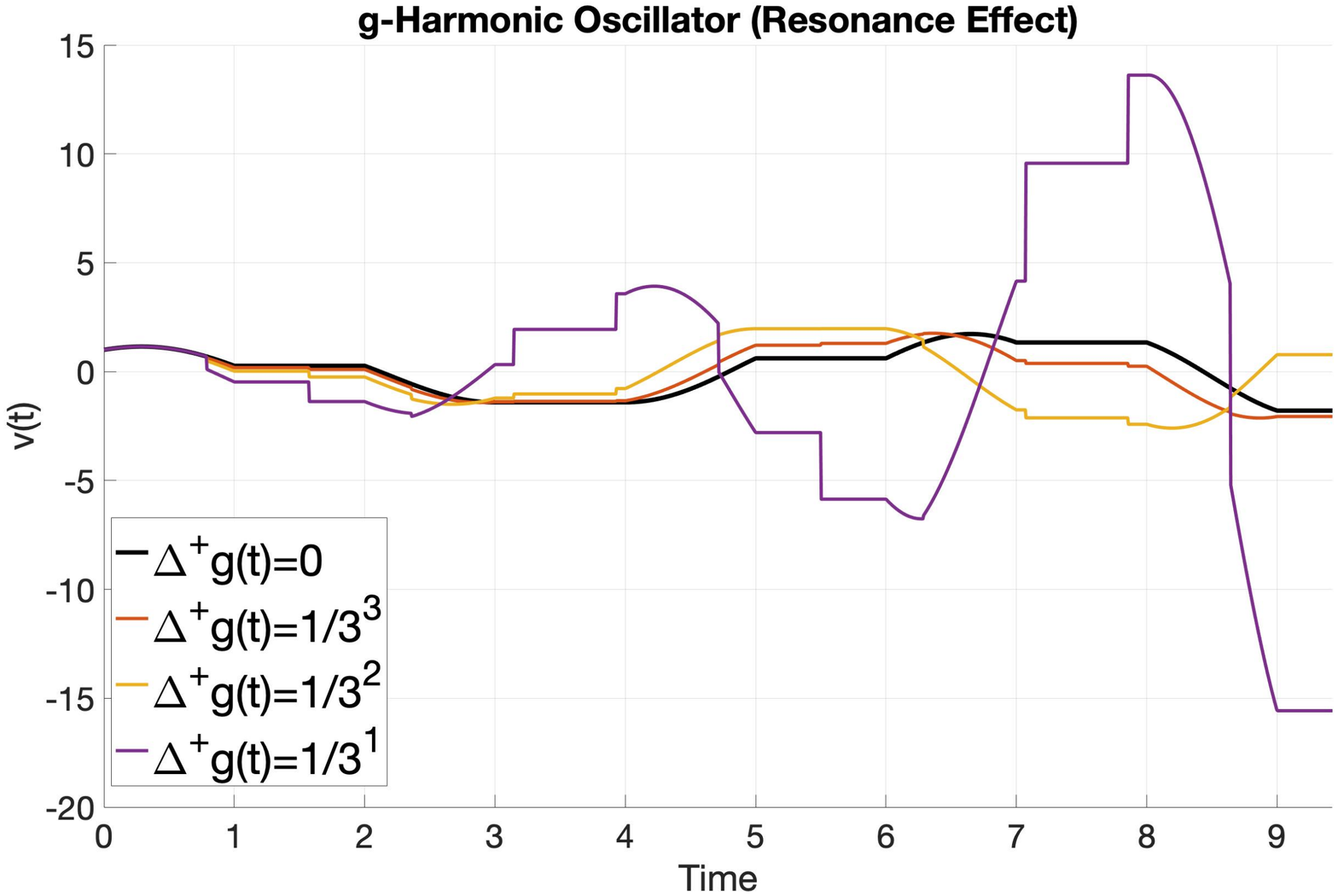}
			\caption{Resonance effect ($g_2$).}
		\end{subfigure}
		\vspace{+6pt}
		\caption{Solution of the $g$-harmonic oscillator associated to
			$g_1$ and $g_2$ (resonance effect) {(vertical lines have to be understood as jumps and not as a multivalued function)}.}
		\label{fig5}
	\end{figure}

	In order to validate the exact solution~\eqref{resosol}, let us compare the
	solution~\eqref{resosol} with the numerical solution of the system:
	\begin{equation}
		\left\{
		\begin{aligned} & 
			\left(\begin{array}{c}
				u\\v
			\end{array}\right)^{\prime}_g(t)=\left(\begin{array}{cc}
				0 & -\omega_0^2 \\
				1 & 0
			\end{array}\right) \left(\begin{array}{c}
				u(t)\\v(t)
			\end{array}\right)+\left( \begin{array}{c}
				\cos_g(\omega_0;0,t) \\
				0
			\end{array}\right), \\
		& 	u(0)=v_0,\  v(0)=x_0.
		\end{aligned}\right.\end{equation}
	The numerical approximation of the solution of
	a system of Stieltjes differential equations was introduced in \cite{FerFraTo2020}, where the authors presented a
	predictor-corrector numerical scheme to approximate the solution of a
	Stieljes dieferential equation (also for systems) from a quadrature formula for the Lebesgue
	Stieltjes integral. For this, a finite set of times $\{t_j\}_{j=0}^{N+1}\subset [0,T]$
	is considered such that $D_g \subset \{t_j\}_{j=0}^{N+1}$, $t_0=0$, $t_{N+1}=T$ and
	$t_{k+1}-t_k=h>0$, for every $k=1,\ldots,N$. The application of the numerical method to our case as follows. Given an element $\mathbf{y}_0=(v_0,x_0)$, we compute
	$\{(\mathbf{y}^+_{j-1},\mathbf{y}^*_j,\mathbf{y}_j)\}_{j=1}^{N+1}$ as:
	\begin{equation} \label{eq:scheme1}
		\left\{\begin{array}{rcl}
			y_{i,j}^+ &=&\displaystyle y_{i,k}+F_{i}(t_j,\mathbf{y}_j)\, \Delta^+ g(t_j), \vspace{0.2cm} \\
			y_{i,j+1}^*&=&\displaystyle y_{i,k}^+ + F_{i} (t_j^+, \mathbf{y}_j^+ )\,(g(t_{j+1})-g(t_j^+)),
			\vspace{0.2cm} \\
			y_{i,j+1} &=& \displaystyle y_{i,j}^+ + \frac{1}{2} \left(F_{i}(t_j^+, \mathbf{y}_j^+ )+
			F_{i}(t_{j+1}^-,\mathbf{y}_{j+1}^*)\right)\,(g(t_{j+1})-g_s(t_j^+)),
		\end{array}\right.\end{equation}
	for every $j=0,\ldots,N$ and $i=1,2$, being $\mathbf{y}_j=(x_{1,j},x_{2,j})$ and
	\begin{equation}%
		\mathbf{F}(t,\mathbf{y})=\left(\begin{array}{l}
			\cos_g(\omega_0;0,t)-\omega_0^2\, y_2\\
			y_1
		\end{array}\right).\end{equation}
	In Table~\ref{table1} we can see the numerical errors $e_h=\max\{|v(t_j)-y_{2,j}|:\; j=0,\ldots,N+1\}$
	for different values of $h$, taking $g(t)=g_2^C(t)+g^B(t)$, with $\Delta^+g=1/3$ (see
	Example~\ref{example1} for the definitions of $g_2^C$ and $g^B$).
	\begin{table}[h]
		\centering
		\noindent\resizebox{.9\textwidth}{!}{
			\begin{tabular}{|l|l|l|l|l|l|l|}
				\hline
				$h$ & $1.e-1$ & $1.e-2$ & $1.e-3$ & $1.e-4$ & $1.e-5$ & $1.e-6$\\
				\hline
				$e_h$ & $4.5260e-01$ & $3.8906e-03$ & $3.8335e-05$ & $3.8274e-07$ & $3.8273e-09$ &
				$3.6102e-11$\\
				\hline
		\end{tabular}}
		\caption{Numerical errors $e_h=\max\{|v(t_j)-y_{2,j}|:\; j=0,\ldots,N+1\}$ for different values of $h$.}
		\label{table1}
	\end{table}
	Finally, in Figure~\ref{fig6}, we can see the comparison between the exact solution
	and the numerical approximation for $h=1.e-1$ and $h=1.e-2$.

	\begin{figure}[h]
		\begin{subfigure}{.5\textwidth}
			\centering
			\includegraphics[width=1\linewidth]{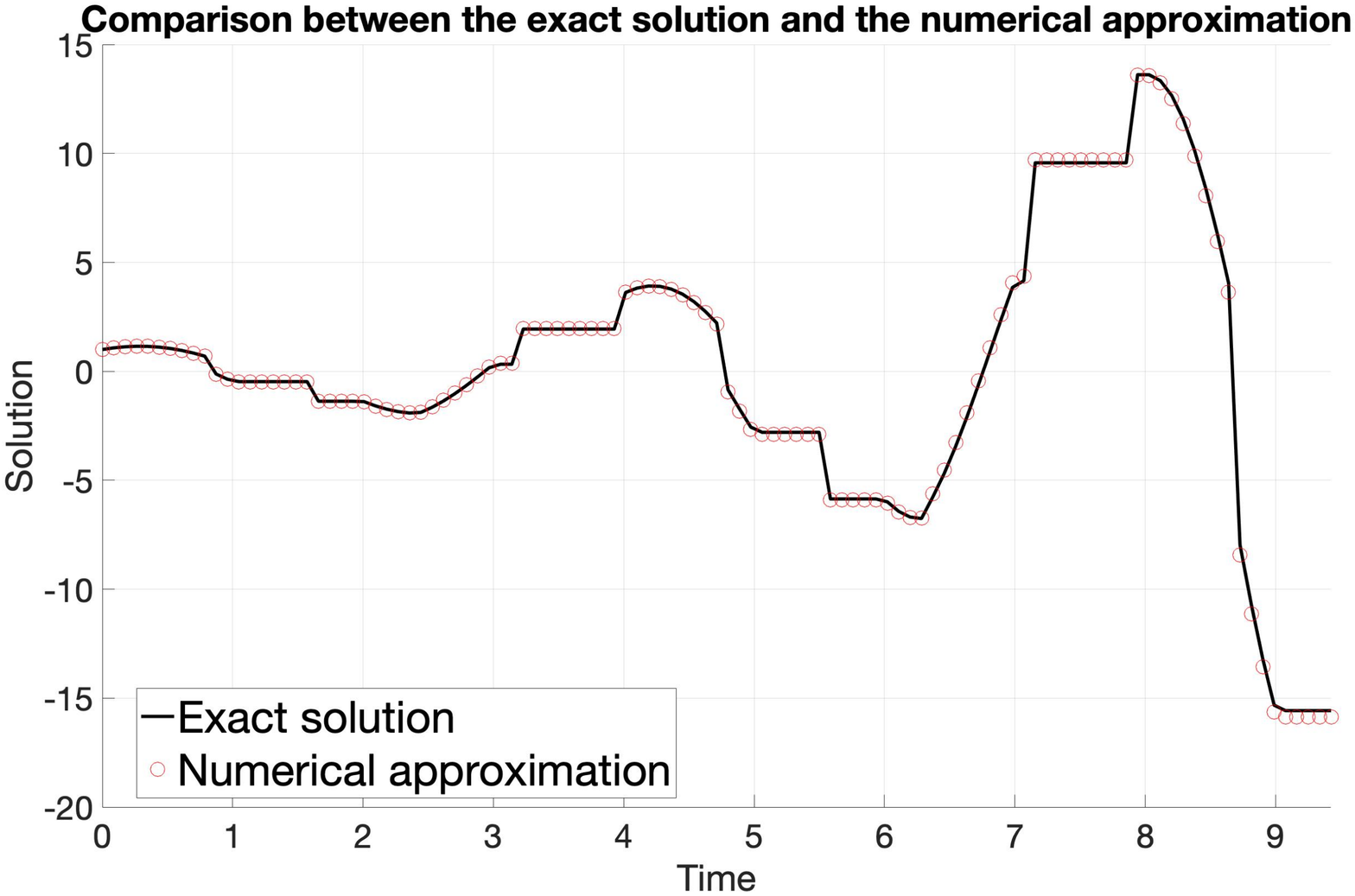}
			\caption{Exact solution vs. numerical approximation \\left( $h=1.e-1$).}
		\end{subfigure}
		\begin{subfigure}{.5\textwidth}
			\centering
			\includegraphics[width=1\linewidth]{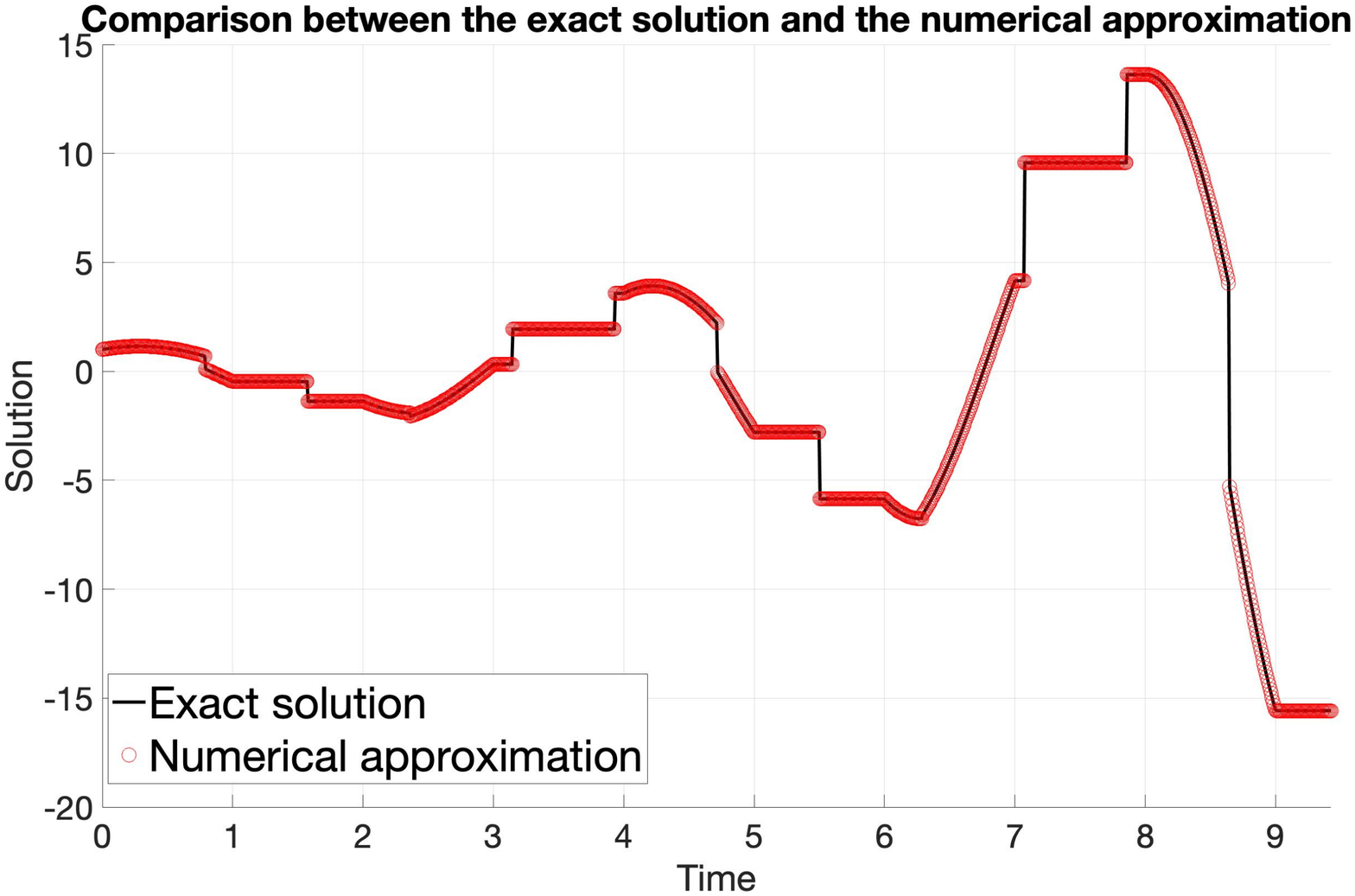}
			\caption{Exact solution vs. numerical approximation \\left( $h=1.e-2$).}
		\end{subfigure}
		\vspace{+6pt}
		\caption{Comparison between the exact solution and the numerical approximation
			(vertical lines have to be understood as jumps and not as a multivalued function).}
		\label{fig6}
	\end{figure}
\end{exa}

\section*{Acknowledgments}
The authors were partially supported by Xunta de Galicia, project ED431C 2019/02, and by the Agencia Estatal de Investigaci\'on (AEI) of Spain under grant MTM2016-75140-P, co-financed by the European Community fund FEDER.
%\section{Conclusions}

\bibliography{refs-gode}

\end{document}